\documentclass{amsproc}
\usepackage{color}

\newtheorem{theorem}{Theorem}[section]
\newtheorem{lemma}[theorem]{Lemma}
\newtheorem{proposition}[theorem]{Proposition}
\newtheorem{corollary}[theorem]{Corollary}
\theoremstyle{definition}
\newtheorem{definition}[theorem]{Definition}
\newtheorem{example}[theorem]{Example}

\newtheorem{assumption}[theorem]{Assumption}

\theoremstyle{remark}
\newtheorem{remark}[theorem]{Remark}
\newtheorem{notations}[theorem]{Notations}

\numberwithin{equation}{section}

\newcommand{\abs}[1]{\lvert#1\rvert}
\newcommand{\Abs}[1]{\Vert\lvert#1\rvert\Vert}
\newcommand{\blankbox}[2]{%
  \parbox{\columnwidth}{\centering

    \setlength{\fboxsep}{0pt}%
    \fbox{\raisebox{0pt}[#2]{\hspace{#1}}}%
  }%
}

\newcommand{\s}{{\rm supp}}
\newcommand{\R}{{\mathbb R}}
\newcommand{\RR}{\mathcal R}
\newcommand{\Z}{{\mathbb Z}}
\newcommand{\FF}{{\mathcal F}}
\newcommand{\GG}{{\mathcal G}}
\newcommand{\Int}{{{\rm Int}\,}}
\newcommand{\A}{{\mathbb A}}
\newcommand{\N}{{\mathbb N}}
\newcommand{\LL}{{\mathcal L}}
\newcommand{\IA}{{\stackrel{\circ}{\A}}}
\newcommand{\C}{{\mathbb C}}
\newcommand{\Fix}{{\rm Fix}}
\newcommand{\Moeb}{\mbox{M\"ob}_+(S^1_\C)}
\newcommand{\BB}{{\mathcal B}}
\newcommand{\PP}{{\mathcal P}}
\newcommand{\EE}{{\mathbb E}}
\newcommand{\HH}{{\mathcal H}}
\newcommand{\lev}{{\rm lev}}
\newcommand{\Stab}{{\rm Stab}}
\newcommand{\bd}{{\sc Proof}.\ }
\newcommand{\QQ}{{\mathcal Q}}
\newcommand{\VV}{{\mathbb V}}
\newcommand{\JJ}{{\mathcal J}}
\newcommand{\TT}{{\mathcal T}}
\newcommand{\II}{{\stackrel{\circ} I}}
\newcommand{\VVV}{{\mathcal V}}
\newcommand{\EEE}{{\mathcal E}}
\newcommand{\DD}{{\mathbb D}}

\begin{document}
\title[Basic partitions and combinations]
{Basic partitions and combinations of group actions on the circle:
A new approach to a theorem of Kathryn Mann}


\author{Shigenori Matsumoto}
\address{Department of Mathematics, College of
Science and Technology, Nihon University, 1-8-14 Kanda, Surugadai,
Chiyoda-ku, Tokyo, 101-8308 Japan
}
\email{matsumo@math.cst.nihon-u.ac.jp
}
\thanks{The author is partially supported by Grant-in-Aid for
Scientific Research (C) No.\ 25400096.}


\date{\today }
\begin{abstract}
Let $\Pi_g$ be the surface group of genus $g$ ($g\geq2$),
and denote by $\RR_{\Pi_g}$ the space of the homomorphisms from
$\Pi_g$ into the group of the orientation preserving homeomorphisms of
$S^1$. Let $2g-2=kl$  for some positive integers $k$ and $l$. 
Then the subset of $\RR_{\Pi_g}$ formed by those $\varphi$ which
are semiconjugate to $k$-fold lifts of some homomorphisms and which have
 Euler
number $eu(\varphi)=l$ is shown to be clopen. This leads to a new
proof of the main result of Kathryn Mann \cite{Mann} from a completely
 different approach.
\end{abstract}

\maketitle
\section{Introduction}

Let $S^1=\R/\Z$ and denote the canonical projection by $\pi:\R\to S^1$.
Denote by $T:\R\to\R$ the translation by one: $T(x)=x+1$.

\begin{notations} Let $\HH={\rm Homeo}_+(S^1)$ denote the group of the
 orientation preserving homeomorphisms of $S^1$, and for any group $G$,
$\RR_G={\rm Homo}(G,\HH)$ the set of the homomorphisms from $G$ to
 $\HH$.
\end{notations}

\begin{definition}\label{dom} A map $h:S^1\to S^1$ is called {\em degree one
 monotone} if there is a nondecreasing (not necessarily continuous) map $\tilde h:\R\to\R$ such that
$\tilde h\circ T=T\circ\tilde h$ and $\pi\circ\tilde h=h\circ \pi$.
\end{definition}

Denote 
$$\RR_G^*=\{\varphi\in\RR_G\mid \exists x\in S^1\ \mbox{such that }\ \varphi(g)(x)=x,\ \forall g\in G\}.$$ 

\begin{definition}\label{semiconjugacy}
Two homomorphisms $\varphi^1,\varphi^2\in \RR_G$ are called {\em semiconjugate},
denoted $\varphi^1\sim\varphi^2$,
if either $\varphi^1,\varphi^2\in\RR_G^*$ or $\varphi^1,\varphi^2\in \RR_G\setminus \RR_G^*$
and there is a degree one monotone map $h:S^1\to S^1$ such that
$\varphi^2(g)\circ h=h\circ\varphi^1(g)$ for any $g\in G$.
\end{definition}

The proof of the following proposition can be found in Appendix A.

\begin{proposition}\label{first}
The semiconjugacy is an equivalence relation.
\end{proposition}

\begin{definition}
Let $F^i\subset S^1$ be a $\varphi^i(G)$-invariant subset
 ($\varphi^i\in\RR_G$, $i=1,2$). A map $\xi:F^1\to F^2$ is called
{\em $(\varphi^1,\varphi^2)$-equivariant} if
 $\xi\circ\varphi^1(g)=\varphi^2(g)\circ\xi$ on $F^1$
for any $g\in G$.
\end{definition}

We have the following easy proposition.

\begin{proposition}\label{equivariant}
Let $F^i\subset S^1$ be a $\varphi^i(G)$-invariant subset
 ($\varphi^i\in\RR_G$, $i=1,2$), and assume there is a cyclic order
 preserving  $(\varphi^1,\varphi^2)$-equivariant
bijection $\xi:F^1\to F^2$. Then we have $\varphi^1\sim\varphi^2$. 
\end{proposition}

\bd 
Two homomorphisms $\varphi_1\in\RR_G^*$ and
$\varphi_2\in\RR_G\setminus\RR_G^*$ can never satisfy the condition
of the proposition.
So one may assume $\varphi^i\in\RR_G\setminus\RR_G^*$.
There is an order preserving bijection 
$\widetilde\xi:\pi^{-1}(F^1)\to\pi^{-1}(F^2)$ such that 
$\widetilde\xi\circ T=T\circ\widetilde\xi$ and 
$\xi\circ\pi=\pi\circ\widetilde\xi$. Define $\widetilde h:\R\to\R$
by
$$
\widetilde h(\widetilde x)=\inf\{\widetilde\xi(\widetilde y)\mid\widetilde
y\in[\widetilde x,\infty)\cap\pi^{-1}(F^1)\}.$$
Then $\widetilde h\circ T=T\circ \widetilde h$,
and there is a monotone degree one map $h:S^1\to S^1$ such that
$h\circ\pi=\pi\circ\widetilde h$. 
Now $(\varphi^1,\varphi^2)$-equivariance of $\xi$ implies
that $h\circ\varphi^1(g)=\varphi^2(g)\circ h$ ($\forall g\in G$). \qed

\medskip

\begin{definition} A homomorphism $\varphi\in\RR_G$ is called {\em type
 0} if there is a $\varphi(G)$-invariant probability measure on $S^1$.
\end{definition}

If there is a finite $\varphi(G)$-orbit or if the
action of $\varphi(G)$ is free, then $\varphi$ is type 0.
If $\varphi$ is type 0 and $\varphi\sim\varphi'$, then $\varphi'$ is also type 0.
If $\varphi$ is not type 0, then the minimal set of $\varphi$ is unique,
either a Cantor set or the whole $S^1$. In the latter case we say that
$\varphi$ is minimal. 

\begin{definition} For $\varphi$ not of type 0, a minimal homomorphism which
 is semiconjugate to $\varphi$ is denoted by $\varphi_\sharp$, and called
a {\em minimal model}.
\end{definition}

A minimal model $\varphi_\sharp$ always exists and is unique up to
topological conjugacy for
$\varphi$ not of type 0.
For any $k\geq2$, let $\pi_k:S^1\to S^1$ be the $k$-fold covering map,
that is, $\pi_k(x+\Z)=kx+\Z$.

\begin{definition}
For $k\in\N$,  $\psi\in\RR_G$ is called a {\em $k$-fold lift} of
$\varphi\in\RR_G$ if for any $g\in G$, it holds that
 $\varphi(g)\circ\pi_k=\pi_k\circ\psi(g)$.
\end{definition}

\begin{definition}
For $k\in\N$, a homomorphism $\varphi\in\RR_G$ is called {\em type $k$} if
 it satisfies the following conditions.
\\
(1) $\varphi$ is not type $0$.
\\
(2) A minimal model $\varphi_\sharp$ is a $k$-fold lift of some
 homomorphism in $\RR_G$.
\\
(3) $k$ is the maximal among those which satisfy (2).

For $k\geq0$, the set of type $k$ homomorphisms is denoted by
 $\RR_G(k)$.
\end{definition}

Thus type 1 homomorphisms are those homomorphisms which are not type 0
 and whose
minimal model cannot be a $k$-fold lift for any $k\geq2$.

\bigskip
The group $\HH$ is a topological group with the uniform convergence
topology, defined by the metric:
$$
d(f,h)=\sup_{x\in S^1}\abs{f(x)-h(x)}\  \mbox{ for }\ f,h\in\HH.$$
The space $\RR_G$ is equipped with the following topology. Given
$\varphi\in\RR_G$, $g\in G$ and $\varepsilon>0$, let
\begin{equation}\label{e1.1}
U(\varphi;g,\varepsilon)=\{\varphi'\in\RR_G\mid d(\varphi'(g),\varphi(g))<\varepsilon\}.
\end{equation}
The topology with subbase $U(\varphi;g,\varepsilon)$ is called
the {\em weak topology}. When the group $G$ is finitely generated, this
coincides with the usual topology of uniform convergence on generators. 
The following proposition will be proven in the next section.

\begin{proposition}\label{p1}
For any group $G$ and $k\geq1$, the subset $\RR_G(0)$ is closed and $\bigcup_{1\leq i\leq k}\RR_G(i)$ is
 open in $\RR_G$.
\end{proposition}

This is best possible, for example for free groups. However for 
groups of a special kind, one can expect that some component
of $\RR_G(k)$, $k\geq2$, is
also open. The purpose of this paper is to consider this problem for
the surface group $\Pi_g$, $g\geq2$. The group $\Pi_g$ is the
fundamental group of the closed oriented surface of genus $g$, and has
a presentation:
$$
\Pi_g=\langle A_1,B_1,\ldots,A_g,B_g\mid[A_1,B_1]\cdots[A_g,B_g]=e\rangle.$$

Given $\varphi\in\RR_{\Pi_g}$, its {\em Euler number} $eu(\varphi)\in\Z$ is defined
by
$$
[\widetilde{\varphi(A_1)},\widetilde{\varphi(B_1)}]\cdots
[\widetilde{\varphi(A_g)},\widetilde{\varphi(B_g)}]=T^{eu(\varphi)},$$
where for $f\in\HH$, $\widetilde f$ denotes an arbitrary lift of $f$ to
a homeomorphism of $\R$.
The map $eu:\RR_{\Pi_g}\to\Z$ is continuous, and thus
$eu^{-1}(i)$ is clopen in $\RR_{\Pi_g}$ for any $i\in\Z$.
We have the following
classical theorem \cite{Milnor}, \cite{Wood}, called the Milnor-Wood inequality.

\begin{theorem}\label{Milor-Wood}
 The inverse image $eu^{-1}(i)$ is nonempty if and only if $\abs{i}\leq2g-2$.
\end{theorem}

For homomorphisms with the extremal values of Euler number, we have
 the following result \cite{Matsumoto2}. (In fact, the pathwise
 connectedness below is not mentioned in that paper. But it is an easy
 consequence of the main theorem.)

\begin{theorem}\label{M2}
The inverse image $E_+=eu^{-1}(2g-2)$ is pathwise connected, 
and if $\varphi,\varphi'\in E_+$,    then $\varphi\sim\varphi'$.
The same thing holds true for $E_-=eu^{-1}(-2g+2)$.
\end{theorem}

Assume $eu(\varphi)=2g-2$ and $2g-2=kl$ for some positive integers $k,l$.
Choose an arbitary $k$-fold
lift $\widehat{\varphi(A_j)}$ (resp.\ $\widehat{\varphi(B_j)}$) of $\varphi(A_j)$
(resp.\ $\varphi(B_j)$) for $j=1,\ldots,g$. Then we have
$$
[\widehat{\varphi(A_1)},\widehat{\varphi(B_1)}]\ldots[\widehat{\varphi(A_g)},\widehat{\varphi(B_g)}]={\rm
Id}.$$
In fact, this is obtained by taking a quotient by the action of $T^l$
of the formula:
$$
[\widetilde{\varphi(A_1)},\widetilde{\varphi(B_1)}]\cdots
[\widetilde{\varphi(A_g)},\widetilde{\varphi(B_g)}]=T^{2g-2}=T^{kl}.$$

Thus we have a $k$-fold lift of $\varphi$ once we  
choose
 $k$-fold lifts of the generators arbitrarily. We shall denote the $k$-fold lifts of $\varphi$
by $\psi_j$, $1\leq j\leq k^{2g}$. The following result is immediate.

\begin{proposition}\label{p2}
We have $eu(\psi_j)=l$. \qed
\end{proposition}

The main result of the present paper is the following.

\begin{theorem}\label{main}
Assume $2g-2=kl$ for some positive integers $k$ and $l$.
Then the subset $eu^{-1}(l)\cap\RR_{\Pi_g}(k)$ is 
 clopen in $\RR_{\Pi_g}$.
\end{theorem}

The closedness of $eu^{-1}(l)\cap\RR_{\Pi_g}(k)$ follows from
Proposition \ref{p1}. In fact,  we have
 $$eu^{-1}(l)\cap\RR_{\Pi_g}(k)=eu^{-1}(l)\setminus \cup_{1\leq j\leq
 k-1}\RR_{\Pi_g}(j),$$ where $eu^{-1}(l)$ is closed and 
$\cup_{1\leq j\leq
 k-1}\RR_{\Pi_g}(j)$ is open.

 For the openness, we use the following concept.

\begin{definition} For any group $G$, a homomorpism $\varphi\in\RR_G$ is
 said to be {\em locally stable} if any homomorphism $\varphi'\in\RR_G$
 sufficiently near to $\varphi$ is semi-conjugate to $\varphi$.
\end{definition}

The openness follows from the following theorem.

\begin{theorem}\label{locally stable}
Any homomorphism of $eu^{-1}(l)\cap\RR_{\Pi_g}(k)$ is locally stable.
\end{theorem}

Let $Z_j$ be the connected component of $\RR_{\Pi_g}$ which
contains
the above lift $\psi_j$, $1\leq j\leq k^{2g}$. Then we have the following
corollary.

\begin{corollary}\label{corollary}
Any two homomorphisms of the same component $Z_j$ are mutually semi-conjugate. \qed
\end{corollary}

The same result has been obtained by K. Mann
\cite{Mann}, based upon extensive use of algorythms in \cite{Calegari}. 
This paper contains a completely different approach.
Also there is a
quite simple proof for diffeomorphisms due to J. Bowden
\cite{Bowden}. 

We shall prove Proposition \ref{p1} in Section 2, and Theorem
\ref{locally stable} in Sections 4--7. We give an outline of Sections 4
and 5
in Section 3.
It seems that our method provides a new and elementary
 proof of the main result
of \cite{Matsumoto2}, but we do not pursue it in the present paper. 
Throughout the paper, we use the following notations.

\begin{notations} \label{n1}

$\bullet$ The positive cyclic order of $S^1$ is denoted by $\prec$.

$\bullet$ Given two distinct points $a,b\in S^1$, $[a,b]=\{x\in S^1,a\prec x\prec b\}$. 
\\
For a subset $X$ of $S^1$, we denote

$\bullet$ $C\sqsubset X$ if $C$ is a connected component of $X$,

$\bullet$ $X_\sharp$ the union of the closures of the connected components of
$S^1\setminus X$,

$\bullet$ $X_*=X\cap X_\sharp$.
\\
We abbreviate

$\bullet$ BP for ``basic partition'', BC for ``basic configuration''
 and COP for ``cyclic order preserving''.
\end{notations}

{\sc Acknowledgement.} Hearty thanks are due to the referee, whose
valuable comments are helpful for the improvement of the paper.

\section{Proximal actions}

In this section, $G$ is to be an arbitrary group, countable or not.
This section is devoted to the proof of Proposition \ref{p1}.
Let us begin by showing that $\RR_G(0)$ is a closed subset of $\RR_G$.
Let $\varphi$ be any homomorphism from the closure of $\RR_G(0)$.
Let us denote by $\PP(S^1)$ the space of the probability measures on $S^1$,
equipped with the weak* topology.
In order to show $\varphi$ admits an invariant probability measure, it is
sufficient to prove that for any finite subset $\{g_i\}\subset G$,
there is a probability measure invariant by
$\varphi(g_i)_*:\PP(S^1)\to\PP(S^1 )$, thanks to the finite intersection
property of the compact set $\PP(S^1)$.
Choose 
$$\varphi_n\in \bigcap_iU(\varphi;g_i,1/n)\cap \RR_G(0),$$
where $U(\cdot)$ is introduced in (\ref{e1.1}), and let
$\mu_n\in\PP(S^1)$ be a $\varphi_n(G)$-invariant measure.
Since the maps $\varphi_n(g_i)_*$ and $\varphi(g_i)_*$ are continuous and
$\varphi_n(g_i)_*$ converges to $\varphi(g_i)_*$ pointwise, an accumulation
point of $\{\mu_n\}$ is the desired measure.

\bigskip
Now let us turn to show that $\RR_G(1)$ is an open subset of $\RR_G$.
The argument is based upon the following Theorem \ref{Ghys} due to
\'E. Ghys (p.362, \cite{Ghys1}), whose proof is included in Appendix B.
 To state it, we make a definition.

\begin{definition}
A homomorphism $\varphi\in\RR_G$ is called {\em proximal} if
for any closed interval $I\subset S^1$, $\inf_{g\in G}\abs{\varphi(g)I}=0$,
where $\abs{\cdot}$ denotes the diameter.
\end{definition}

\begin{theorem}\label{Ghys}
For any $\varphi\in\RR_G$, $\varphi\in\RR_G(1)$ if and only if a minimal
 model $\varphi_\sharp$ is
 proximal.
\end{theorem}

\begin{definition} Given $x,y\in S^1$, a sequence $\{f_n\}\subset \HH$
 is called an $(x,y)$-sequence if for any $\varepsilon>0$, there is $N$ such
 that if $n\geq N$, $f_n$ maps the complement of the
 $\varepsilon$-neighbourhood of $x$ into the $\varepsilon$-neighbourhood of
 $y$.
\end{definition}

\begin{lemma}\label{l2.1}
For any $x,y\in S^1$ and $\varphi\in\RR_G(1)$, there is an $(x,y)$-sequence
in $\varphi_\sharp(G)$. 
\end{lemma}
\bd For any $x\in S^1$, define
$$
E_x=\{y\in S^1\mid \exists\mbox{$(x,y)$-sequence in $\varphi_\sharp(G)$}\}.
$$

By Theorem \ref{Ghys}, $E_x$ is nonempty for any $x\in S^1$. On the
other hand, it is easy to show that $E_x$ is closed and
$\varphi_\sharp(G)$-invariant. Therefore we have $E_x=S^1$. \qed

\medskip
There is a bounded 2-cocycle $c$ of the group $\HH$ defined by
$$
c(f,h)=\tau(\widetilde f\circ\widetilde h)-\tau(\widetilde
f)-\tau(\widetilde h),$$
where $\tilde f$ (resp.\ $\tilde h$) is an arbitrary lift of $f$ (resp.\ $h$) to $\R$,
and $\tau(\cdot)$ stands for the translation number.
As is well known, its $L^\infty$ norm satisfies $\Vert c\Vert=1$.
For $\varphi\in\RR_G$, the pull back cocycle $\varphi^*c$ lies in the second
bounded cocycle group $Z^2_b(G)$ of $G$ and satisfies
$\Vert\varphi^*c\Vert\leq1$. It is known \cite{Matsumoto1} that $\varphi^*c=0$
if and only if $\varphi\in\RR_G(0)$. For other $\RR_G(k)$, we have the following.

\begin{lemma}\label{l2.2}
For any $\varphi\in\RR_G$ and $k\geq1$, $\varphi\in\RR_G(k)$ if and only if
$\Vert\varphi^*c\Vert=1/k$.
\end{lemma}

\bd It suffices to show only the following implication:
\begin{equation}\label{e2.1}
\varphi\in\RR_G(k)\Rightarrow\Vert\varphi^*c\Vert=1/k,\ \ \forall k\geq1,
\end{equation}
since the opposite implication follows from this.
First of all, let us show (\ref{e2.1}) for $k=1$. 
Let $\varphi_\sharp$ be a minimal model of any $\varphi\in\RR_G(1)$.
Choose four points
$x\prec y\prec z\prec u\prec x$ in $S^1$. By Lemma \ref{l2.1}, there
are a $(y,x)$-sequence $f_n$ and an $(u,z)$-sequence $h_n$ in $\varphi_\sharp(G)$.
Let $\widetilde f_n$ and $\widetilde h_n$ be the lifts of $f_n$ and
$h_n$ such that $\tau(\widetilde f_n)=\tau(\widetilde h_n)=0$.  
One can choose lifts of the four points so that
$\widetilde x<\widetilde y<\widetilde z<\widetilde u<T(\widetilde x)$.
See Figure 1 for this and the next argument.

\begin{figure}\caption{}

\unitlength 0.1in
\begin{picture}( 47.7000, 31.6000)(  9.5000,-36.9000)
%
{\color[named]{Black}{%
\special{pn 8}%
\special{pa 1270 830}%
\special{pa 5710 830}%
\special{fp}%
}}%
%
{\color[named]{Black}{%
\special{pn 8}%
\special{pa 1290 1330}%
\special{pa 5720 1330}%
\special{fp}%
}}%
%
{\color[named]{Black}{%
\special{pn 4}%
\special{sh 1}%
\special{ar 2530 830 16 16 0  6.28318530717959E+0000}%
\special{sh 1}%
\special{ar 2530 830 16 16 0  6.28318530717959E+0000}%
}}%
%
{\color[named]{Black}{%
\special{pn 4}%
\special{sh 1}%
\special{ar 2770 830 16 16 0  6.28318530717959E+0000}%
\special{sh 1}%
\special{ar 2770 830 16 16 0  6.28318530717959E+0000}%
}}%
%
{\color[named]{Black}{%
\special{pn 4}%
\special{sh 1}%
\special{ar 4610 830 16 16 0  6.28318530717959E+0000}%
\special{sh 1}%
\special{ar 4610 830 16 16 0  6.28318530717959E+0000}%
}}%
%
{\color[named]{Black}{%
\special{pn 4}%
\special{sh 1}%
\special{ar 4850 830 16 16 0  6.28318530717959E+0000}%
\special{sh 1}%
\special{ar 4850 830 16 16 0  6.28318530717959E+0000}%
}}%
%
\put(35.7000,-13.3000){\makebox(0,0)[lb]{}}%
%
{\color[named]{Black}{%
\special{pn 4}%
\special{sh 1}%
\special{ar 3870 1330 4 4 0  6.28318530717959E+0000}%
\special{sh 1}%
\special{ar 3870 1330 4 4 0  6.28318530717959E+0000}%
}}%
%
{\color[named]{Black}{%
\special{pn 4}%
\special{pa 1640 830}%
\special{pa 1760 830}%
\special{fp}%
\special{sh 1}%
\special{pa 1760 830}%
\special{pa 1694 810}%
\special{pa 1708 830}%
\special{pa 1694 850}%
\special{pa 1760 830}%
\special{fp}%
}}%
%
{\color[named]{Black}{%
\special{pn 4}%
\special{pa 2670 830}%
\special{pa 2580 830}%
\special{fp}%
\special{sh 1}%
\special{pa 2580 830}%
\special{pa 2648 850}%
\special{pa 2634 830}%
\special{pa 2648 810}%
\special{pa 2580 830}%
\special{fp}%
}}%
%
{\color[named]{Black}{%
\special{pn 4}%
\special{pa 3560 820}%
\special{pa 3670 820}%
\special{fp}%
\special{sh 1}%
\special{pa 3670 820}%
\special{pa 3604 800}%
\special{pa 3618 820}%
\special{pa 3604 840}%
\special{pa 3670 820}%
\special{fp}%
}}%
%
{\color[named]{Black}{%
\special{pn 4}%
\special{pa 4760 820}%
\special{pa 4710 830}%
\special{fp}%
\special{sh 1}%
\special{pa 4710 830}%
\special{pa 4780 838}%
\special{pa 4762 820}%
\special{pa 4772 798}%
\special{pa 4710 830}%
\special{fp}%
}}%
%
{\color[named]{Black}{%
\special{pn 4}%
\special{pa 5180 830}%
\special{pa 5240 830}%
\special{fp}%
\special{sh 1}%
\special{pa 5240 830}%
\special{pa 5174 810}%
\special{pa 5188 830}%
\special{pa 5174 850}%
\special{pa 5240 830}%
\special{fp}%
}}%
%
{\color[named]{Black}{%
\special{pn 4}%
\special{pa 1460 1330}%
\special{pa 1530 1320}%
\special{fp}%
\special{sh 1}%
\special{pa 1530 1320}%
\special{pa 1462 1310}%
\special{pa 1478 1328}%
\special{pa 1468 1350}%
\special{pa 1530 1320}%
\special{fp}%
}}%
%
{\color[named]{Black}{%
\special{pn 4}%
\special{pa 1950 1330}%
\special{pa 1910 1330}%
\special{fp}%
\special{sh 1}%
\special{pa 1910 1330}%
\special{pa 1978 1350}%
\special{pa 1964 1330}%
\special{pa 1978 1310}%
\special{pa 1910 1330}%
\special{fp}%
}}%
%
{\color[named]{Black}{%
\special{pn 4}%
\special{pa 2660 1330}%
\special{pa 2690 1330}%
\special{fp}%
\special{sh 1}%
\special{pa 2690 1330}%
\special{pa 2624 1310}%
\special{pa 2638 1330}%
\special{pa 2624 1350}%
\special{pa 2690 1330}%
\special{fp}%
}}%
%
{\color[named]{Black}{%
\special{pn 4}%
\special{pa 3760 1330}%
\special{pa 3710 1330}%
\special{fp}%
\special{sh 1}%
\special{pa 3710 1330}%
\special{pa 3778 1350}%
\special{pa 3764 1330}%
\special{pa 3778 1310}%
\special{pa 3710 1330}%
\special{fp}%
}}%
%
{\color[named]{Black}{%
\special{pn 4}%
\special{pa 4520 1330}%
\special{pa 4610 1330}%
\special{fp}%
\special{sh 1}%
\special{pa 4610 1330}%
\special{pa 4544 1310}%
\special{pa 4558 1330}%
\special{pa 4544 1350}%
\special{pa 4610 1330}%
\special{fp}%
}}%
%
{\color[named]{Black}{%
\special{pn 4}%
\special{sh 1}%
\special{ar 3860 1330 16 16 0  6.28318530717959E+0000}%
\special{sh 1}%
\special{ar 3860 1330 16 16 0  6.28318530717959E+0000}%
}}%
\put(24.6000,-7.1000){\makebox(0,0)[lb]{$\widetilde z$}}%
\put(27.3000,-6.8000){\makebox(0,0)[lb]{$\widetilde u$}}%
\put(44.8000,-6.8000){\makebox(0,0)[lb]{$T(\widetilde z)$}}%
\put(49.1000,-6.6000){\makebox(0,0)[lb]{$T(\widetilde u)$}}%
\put(17.1000,-11.9000){\makebox(0,0)[lb]{$\widetilde x$}}%
\put(20.8000,-11.7000){\makebox(0,0)[lb]{$\widetilde y$}}%
\put(34.2000,-11.7000){\makebox(0,0)[lb]{$T(\widetilde x)$}}%
\put(39.0000,-11.4000){\makebox(0,0)[lb]{$T(\widetilde y)$}}%
\put(9.5000,-8.4000){\makebox(0,0)[lb]{$\widetilde h_n$}}%
\put(9.8000,-13.3000){\makebox(0,0)[lb]{$\widetilde f_n$}}%
%
{\color[named]{Black}{%
\special{pn 4}%
\special{sh 1}%
\special{ar 1830 1320 16 16 0  6.28318530717959E+0000}%
\special{sh 1}%
\special{ar 1830 1320 16 16 0  6.28318530717959E+0000}%
}}%
%
{\color[named]{Black}{%
\special{pn 4}%
\special{sh 1}%
\special{ar 2100 1320 16 16 0  6.28318530717959E+0000}%
\special{sh 1}%
\special{ar 2100 1320 16 16 0  6.28318530717959E+0000}%
}}%
%
{\color[named]{Black}{%
\special{pn 4}%
\special{sh 1}%
\special{ar 3600 1340 16 16 0  6.28318530717959E+0000}%
\special{sh 1}%
\special{ar 3600 1340 16 16 0  6.28318530717959E+0000}%
}}%
%
{\color[named]{Black}{%
\special{pn 8}%
\special{pa 1310 2100}%
\special{pa 5700 2110}%
\special{fp}%
}}%
%
{\color[named]{Black}{%
\special{pn 8}%
\special{pa 1340 2780}%
\special{pa 5700 2780}%
\special{fp}%
}}%
%
{\color[named]{Black}{%
\special{pn 8}%
\special{pa 1330 3490}%
\special{pa 5720 3490}%
\special{fp}%
}}%
%
{\color[named]{Black}{%
\special{pn 4}%
\special{sh 1}%
\special{ar 1830 2100 16 16 0  6.28318530717959E+0000}%
\special{sh 1}%
\special{ar 1830 2100 16 16 0  6.28318530717959E+0000}%
}}%
%
{\color[named]{Black}{%
\special{pn 4}%
\special{sh 1}%
\special{ar 2090 2100 16 16 0  6.28318530717959E+0000}%
\special{sh 1}%
\special{ar 2090 2100 16 16 0  6.28318530717959E+0000}%
}}%
%
{\color[named]{Black}{%
\special{pn 4}%
\special{sh 1}%
\special{ar 2530 2100 16 16 0  6.28318530717959E+0000}%
\special{sh 1}%
\special{ar 2530 2100 16 16 0  6.28318530717959E+0000}%
}}%
%
{\color[named]{Black}{%
\special{pn 4}%
\special{sh 1}%
\special{ar 2770 2100 16 16 0  6.28318530717959E+0000}%
\special{sh 1}%
\special{ar 2770 2100 16 16 0  6.28318530717959E+0000}%
}}%
%
{\color[named]{Black}{%
\special{pn 4}%
\special{sh 1}%
\special{ar 3610 2110 16 16 0  6.28318530717959E+0000}%
\special{sh 1}%
\special{ar 3610 2110 16 16 0  6.28318530717959E+0000}%
}}%
%
{\color[named]{Black}{%
\special{pn 4}%
\special{sh 1}%
\special{ar 3880 2110 16 16 0  6.28318530717959E+0000}%
\special{sh 1}%
\special{ar 3880 2110 16 16 0  6.28318530717959E+0000}%
}}%
%
{\color[named]{Black}{%
\special{pn 4}%
\special{sh 1}%
\special{ar 4590 2110 16 16 0  6.28318530717959E+0000}%
\special{sh 1}%
\special{ar 4590 2110 16 16 0  6.28318530717959E+0000}%
}}%
%
\put(48.5000,-21.1000){\makebox(0,0)[lb]{}}%
%
{\color[named]{Black}{%
\special{pn 4}%
\special{sh 1}%
\special{ar 4850 2100 16 16 0  6.28318530717959E+0000}%
\special{sh 1}%
\special{ar 4860 2100 16 16 0  6.28318530717959E+0000}%
}}%
%
{\color[named]{Black}{%
\special{pn 4}%
\special{sh 1}%
\special{ar 1830 2790 16 16 0  6.28318530717959E+0000}%
\special{sh 1}%
\special{ar 1830 2790 16 16 0  6.28318530717959E+0000}%
}}%
%
{\color[named]{Black}{%
\special{pn 4}%
\special{sh 1}%
\special{ar 2090 2780 16 16 0  6.28318530717959E+0000}%
\special{sh 1}%
\special{ar 2090 2780 16 16 0  6.28318530717959E+0000}%
}}%
%
{\color[named]{Black}{%
\special{pn 4}%
\special{sh 1}%
\special{ar 2530 2780 16 16 0  6.28318530717959E+0000}%
\special{sh 1}%
\special{ar 2520 2780 16 16 0  6.28318530717959E+0000}%
}}%
%
{\color[named]{Black}{%
\special{pn 4}%
\special{sh 1}%
\special{ar 2780 2780 16 16 0  6.28318530717959E+0000}%
\special{sh 1}%
\special{ar 2780 2780 16 16 0  6.28318530717959E+0000}%
}}%
%
{\color[named]{Black}{%
\special{pn 4}%
\special{sh 1}%
\special{ar 3610 2780 16 16 0  6.28318530717959E+0000}%
\special{sh 1}%
\special{ar 3620 2770 16 16 0  6.28318530717959E+0000}%
\special{sh 1}%
\special{ar 3610 2770 16 16 0  6.28318530717959E+0000}%
}}%
%
{\color[named]{Black}{%
\special{pn 4}%
\special{sh 1}%
\special{ar 3880 2780 16 16 0  6.28318530717959E+0000}%
\special{sh 1}%
\special{ar 3880 2780 16 16 0  6.28318530717959E+0000}%
}}%
%
{\color[named]{Black}{%
\special{pn 4}%
\special{sh 1}%
\special{ar 4590 2770 16 16 0  6.28318530717959E+0000}%
\special{sh 1}%
\special{ar 4580 2770 16 16 0  6.28318530717959E+0000}%
}}%
%
{\color[named]{Black}{%
\special{pn 4}%
\special{sh 1}%
\special{ar 4860 2780 16 16 0  6.28318530717959E+0000}%
\special{sh 1}%
\special{ar 4860 2780 16 16 0  6.28318530717959E+0000}%
}}%
%
{\color[named]{Black}{%
\special{pn 4}%
\special{sh 1}%
\special{ar 1850 3490 16 16 0  6.28318530717959E+0000}%
\special{sh 1}%
\special{ar 1850 3490 16 16 0  6.28318530717959E+0000}%
}}%
%
{\color[named]{Black}{%
\special{pn 4}%
\special{sh 1}%
\special{ar 2100 3490 16 16 0  6.28318530717959E+0000}%
\special{sh 1}%
\special{ar 2100 3490 16 16 0  6.28318530717959E+0000}%
}}%
%
{\color[named]{Black}{%
\special{pn 4}%
\special{sh 1}%
\special{ar 2530 3490 16 16 0  6.28318530717959E+0000}%
\special{sh 1}%
\special{ar 2530 3490 16 16 0  6.28318530717959E+0000}%
}}%
%
{\color[named]{Black}{%
\special{pn 4}%
\special{sh 1}%
\special{ar 2770 3490 16 16 0  6.28318530717959E+0000}%
\special{sh 1}%
\special{ar 2770 3490 16 16 0  6.28318530717959E+0000}%
}}%
%
{\color[named]{Black}{%
\special{pn 4}%
\special{sh 1}%
\special{ar 3620 3490 16 16 0  6.28318530717959E+0000}%
\special{sh 1}%
\special{ar 3620 3490 16 16 0  6.28318530717959E+0000}%
}}%
%
{\color[named]{Black}{%
\special{pn 4}%
\special{sh 1}%
\special{ar 3890 3490 16 16 0  6.28318530717959E+0000}%
\special{sh 1}%
\special{ar 3890 3490 16 16 0  6.28318530717959E+0000}%
}}%
%
{\color[named]{Black}{%
\special{pn 4}%
\special{sh 1}%
\special{ar 4590 3490 16 16 0  6.28318530717959E+0000}%
\special{sh 1}%
\special{ar 4590 3490 16 16 0  6.28318530717959E+0000}%
}}%
%
{\color[named]{Black}{%
\special{pn 4}%
\special{sh 1}%
\special{ar 4860 3490 16 16 0  6.28318530717959E+0000}%
\special{sh 1}%
\special{ar 4860 3490 16 16 0  6.28318530717959E+0000}%
}}%
%
\put(47.5000,-36.9000){\makebox(0,0)[lb]{}}%
\put(47.7000,-36.8000){\makebox(0,0)[lb]{$T(\widetilde u_n)$}}%
\put(27.2000,-19.6000){\makebox(0,0)[lb]{$\widetilde u_n$}}%
\put(27.1000,-37.2000){\makebox(0,0)[lb]{$\widetilde u_n$}}%
%
{\color[named]{Black}{%
\special{pn 8}%
\special{pa 2770 2280}%
\special{pa 2770 2670}%
\special{fp}%
\special{sh 1}%
\special{pa 2770 2670}%
\special{pa 2790 2604}%
\special{pa 2770 2618}%
\special{pa 2750 2604}%
\special{pa 2770 2670}%
\special{fp}%
}}%
%
{\color[named]{Black}{%
\special{pn 8}%
\special{pa 2770 2920}%
\special{pa 2786 2950}%
\special{pa 2802 2978}%
\special{pa 2820 3004}%
\special{pa 2838 3028}%
\special{pa 2860 3050}%
\special{pa 2886 3068}%
\special{pa 2912 3084}%
\special{pa 2940 3098}%
\special{pa 2970 3110}%
\special{pa 3034 3130}%
\special{pa 3066 3138}%
\special{pa 3134 3154}%
\special{pa 3166 3162}%
\special{pa 3200 3170}%
\special{pa 3232 3178}%
\special{pa 3264 3188}%
\special{pa 3294 3200}%
\special{pa 3350 3228}%
\special{pa 3374 3246}%
\special{pa 3396 3266}%
\special{pa 3416 3290}%
\special{pa 3436 3316}%
\special{pa 3468 3372}%
\special{pa 3482 3402}%
\special{pa 3490 3420}%
\special{fp}%
}}%
%
{\color[named]{Black}{%
\special{pn 8}%
\special{pa 3440 3340}%
\special{pa 3490 3440}%
\special{fp}%
\special{sh 1}%
\special{pa 3490 3440}%
\special{pa 3478 3372}%
\special{pa 3466 3392}%
\special{pa 3442 3390}%
\special{pa 3490 3440}%
\special{fp}%
}}%
%
{\color[named]{Black}{%
\special{pn 8}%
\special{pa 2850 2160}%
\special{pa 2876 2180}%
\special{pa 2928 2216}%
\special{pa 2984 2248}%
\special{pa 3010 2262}%
\special{pa 3040 2276}%
\special{pa 3068 2290}%
\special{pa 3098 2300}%
\special{pa 3126 2312}%
\special{pa 3156 2322}%
\special{pa 3188 2330}%
\special{pa 3218 2340}%
\special{pa 3248 2348}%
\special{pa 3376 2372}%
\special{pa 3440 2380}%
\special{pa 3474 2384}%
\special{pa 3506 2388}%
\special{pa 3540 2390}%
\special{pa 3574 2394}%
\special{pa 3606 2396}%
\special{pa 3674 2400}%
\special{pa 3708 2400}%
\special{pa 3810 2406}%
\special{pa 3844 2406}%
\special{pa 3878 2408}%
\special{pa 3912 2412}%
\special{pa 3944 2414}%
\special{pa 3978 2416}%
\special{pa 4042 2424}%
\special{pa 4076 2430}%
\special{pa 4106 2436}%
\special{pa 4138 2442}%
\special{pa 4170 2450}%
\special{pa 4230 2470}%
\special{pa 4286 2494}%
\special{pa 4342 2526}%
\special{pa 4368 2542}%
\special{pa 4392 2562}%
\special{pa 4416 2584}%
\special{pa 4438 2608}%
\special{pa 4456 2634}%
\special{pa 4470 2662}%
\special{pa 4482 2692}%
\special{pa 4490 2720}%
\special{fp}%
}}%
%
{\color[named]{Black}{%
\special{pn 8}%
\special{pa 4460 2650}%
\special{pa 4490 2750}%
\special{fp}%
\special{sh 1}%
\special{pa 4490 2750}%
\special{pa 4490 2680}%
\special{pa 4476 2700}%
\special{pa 4452 2692}%
\special{pa 4490 2750}%
\special{fp}%
}}%
%
{\color[named]{Black}{%
\special{pn 8}%
\special{pa 4480 2840}%
\special{pa 4498 2868}%
\special{pa 4514 2896}%
\special{pa 4532 2924}%
\special{pa 4550 2948}%
\special{pa 4570 2972}%
\special{pa 4592 2992}%
\special{pa 4616 3012}%
\special{pa 4640 3028}%
\special{pa 4666 3044}%
\special{pa 4694 3058}%
\special{pa 4722 3070}%
\special{pa 4752 3080}%
\special{pa 4782 3088}%
\special{pa 4846 3104}%
\special{pa 4878 3110}%
\special{pa 5014 3126}%
\special{pa 5050 3128}%
\special{pa 5084 3132}%
\special{pa 5156 3136}%
\special{pa 5190 3140}%
\special{pa 5226 3142}%
\special{pa 5260 3146}%
\special{pa 5294 3152}%
\special{pa 5328 3156}%
\special{pa 5392 3172}%
\special{pa 5422 3182}%
\special{pa 5450 3194}%
\special{pa 5478 3208}%
\special{pa 5502 3226}%
\special{pa 5526 3246}%
\special{pa 5546 3268}%
\special{pa 5566 3294}%
\special{pa 5584 3320}%
\special{pa 5602 3348}%
\special{pa 5618 3376}%
\special{pa 5620 3380}%
\special{fp}%
}}%
%
{\color[named]{Black}{%
\special{pn 8}%
\special{pa 5560 3310}%
\special{pa 5650 3410}%
\special{fp}%
\special{sh 1}%
\special{pa 5650 3410}%
\special{pa 5620 3348}%
\special{pa 5614 3370}%
\special{pa 5592 3374}%
\special{pa 5650 3410}%
\special{fp}%
}}%
%
{\color[named]{Black}{%
\special{pn 8}%
\special{pa 2880 2050}%
\special{pa 2830 2150}%
\special{fp}%
}}%
%
{\color[named]{Black}{%
\special{pn 8}%
\special{pa 2820 2050}%
\special{pa 2900 2140}%
\special{fp}%
}}%
\put(23.9000,-25.2000){\makebox(0,0)[lb]{$\widetilde h_n$}}%
\put(24.2000,-32.7000){\makebox(0,0)[lb]{$\widetilde f_n$}}%
\end{picture}

\end{figure}
For $n$ large, $\widetilde h_n$ admits 
a fixed point, say $\tilde u_n$, near $\tilde u$. 
Now consider the composite $\widetilde f_n\circ\widetilde h_n$. 
Clearly we have 
$\widetilde u_n<\widetilde f_n\circ\widetilde h_n(\widetilde u_n)< T(\widetilde u_n)$. On the other hand, if we choose $\widetilde u'$ very near to
$\tilde u$ so that $\widetilde u'>\widetilde u$. Then for any large $n$,
we have $\widetilde 
f_n\circ\widetilde h_n(u')>T(u')$. (See Figure 1.)
This shows $\tau(\widetilde f_n\circ\widetilde
h_n)=1$. Therefore $c(f_n,h_n)=1$ and 
$\Vert\varphi^*c\Vert=\Vert \varphi_\sharp^*c\Vert=1$, as is required.
Also it is not difficult to show that the
above inequalities also show the following.
\begin{equation}\label{open}
\mbox{For any $\varphi'\in\RR_G$ sufficiently
near to $\varphi\in\RR_G(1)$, we have $\Vert(\varphi')^*c\Vert=1$.}
\end{equation}

To show (\ref{e2.1}) for $k\geq1$, choose any $\varphi\in\RR_G(k)$, with
$\varphi_\sharp$ a $k$-fold lift of some $\psi\in\RR_G$. Clearly
$\psi\in\RR_G(1)$. Moreover the cocycle $\varphi^*c=\varphi_\sharp^*c$ is precisely
$(1/k)\psi^*c$. This shows $\Vert\varphi^*c\Vert=1/k$. \qed

\medskip
Now the openness of $\RR_G(1)$ follows from Lemma \ref{l2.2} and
(\ref{open}). The proof that the set $\bigcup_{1\leq i\leq k}\RR_G(i)$
is open is left to the reader.

\section{Outline}
Before getting into a detailed proof of Theorem \ref{locally stable},
we shall give an outline of its first two steps.
The basic idea is that a homomorphism in $eu^{-1}(l)\cap\RR_{\Pi_g}(k)$
of Theorem \ref{locally stable} has the following very special property:
There is a finite set, say $R$, of $S^1$ such that the knowledge about
how the generators of the group moves points of $R$ completely
determines the semiconjugacy class of the homomorphism.

First of all, let us explain this phenomenon in a much simpler example.
Let $\Gamma$ be the free group on two generators $A$ and $B$. Let 
$\varphi\in\RR_\Gamma$ and denote $a=\varphi(A)$ and $b=\varphi(B)$.
Assume that $\tau([\widetilde a,\widetilde b])=1$, where $\widetilde a$
(resp.\ $\widetilde b$) is an arbitrary lift of $a$ (resp.\ $b$).
Then one can show that such $\varphi$ belongs to a single semiconjugacy
class. This will be actually done in Section 4. But we can present a rough outline
here.

By the assumption $\tau([\widetilde a,\widetilde b])=1$,
 there is a fixed point $x\in S^1$ of $[a,b]$ such that
 $$x\prec b^{-1}(x)\prec a^{-1}b^{-1}(x)\prec ba^{-1}b^{-1}(x)\prec
 [a,b](x)=x.$$ See Figure 2 left. 

\begin{figure}\caption{}
\unitlength 0.1in
\begin{picture}( 60.1000, 25.4000)(  12.1000,-31.1000)
%
{\color[named]{Black}{%
\special{pn 8}%
\special{ar 2100 1850 1100 1100  1.5798051  1.5707963}%

}}%
%
{\color[named]{Black}{%
\special{pn 8}%
\special{ar 5130 1810 1100 1100  0.0000000  6.2831853}%
}}%
%
{\color[named]{Black}{%
\special{pn 4}%
\special{sh 1}%
\special{ar 2080 760 16 16 0  6.28318530717959E+0000}%
\special{sh 1}%
\special{ar 2090 760 16 16 0  6.28318530717959E+0000}%
}}%
%
{\color[named]{Black}{%
\special{pn 4}%
\special{sh 1}%
\special{ar 5160 720 16 16 0  6.28318530717959E+0000}%
\special{sh 1}%
\special{ar 5160 720 16 16 0  6.28318530717959E+0000}%
}}%
%
{\color[named]{Black}{%
\special{pn 4}%
\special{sh 1}%
\special{ar 3210 1800 16 16 0  6.28318530717959E+0000}%
\special{sh 1}%
\special{ar 3210 1800 16 16 0  6.28318530717959E+0000}%
}}%
%
{\color[named]{Black}{%
\special{pn 4}%
\special{sh 1}%
\special{ar 6240 1780 16 16 0  6.28318530717959E+0000}%
\special{sh 1}%
\special{ar 6240 1780 16 16 0  6.28318530717959E+0000}%
}}%
%
{\color[named]{Black}{%
\special{pn 4}%
\special{sh 1}%
\special{ar 4040 1800 16 16 0  6.28318530717959E+0000}%
\special{sh 1}%
\special{ar 4040 1800 16 16 0  6.28318530717959E+0000}%
}}%
%
{\color[named]{Black}{%
\special{pn 4}%
\special{sh 1}%
\special{ar 1010 1820 16 16 0  6.28318530717959E+0000}%
\special{sh 1}%
\special{ar 1010 1820 16 16 0  6.28318530717959E+0000}%
}}%
%
{\color[named]{Black}{%
\special{pn 4}%
\special{sh 1}%
\special{ar 2080 2950 16 16 0  6.28318530717959E+0000}%
\special{sh 1}%
\special{ar 2080 2950 16 16 0  6.28318530717959E+0000}%
}}%
%
{\color[named]{Black}{%
\special{pn 4}%
\special{sh 1}%
\special{ar 5180 2920 16 16 0  6.28318530717959E+0000}%
\special{sh 1}%
\special{ar 5180 2920 16 16 0  6.28318530717959E+0000}%
}}%
\put(23.9000,-18.8000){\makebox(0,0)[lb]{$x=[a,b](x)$}}%
\put(19.9000,-6.8000){\makebox(0,0)[lb]{$b^{-1}(x)$}}%
\put(20.7000,-31.3000){\makebox(0,0)[lb]{$ba^{-1}b^{-1}(x)$}}%
\put(11.1000,-19.2000){\makebox(0,0)[lb]{$a^{-1}b^{-1}(x)$}}%
\put(56.2000,-20.9000){\makebox(0,0)[lb]{$x=a(x)$}}%
\put(49.7000,-6.0000){\makebox(0,0)[lb]{$y=ab(y)$}}%
\put(34.6000,-18.4000){\makebox(0,0)[lb]{$z=b(z)$}}%
\put(50.5000,-31.4000){\makebox(0,0)[lb]{$b(y)$}}%
%
{\color[named]{Black}{%
\special{pn 4}%
\special{pa 2030 900}%
\special{pa 2014 928}%
\special{pa 1960 1012}%
\special{pa 1944 1038}%
\special{pa 1908 1094}%
\special{pa 1872 1146}%
\special{pa 1852 1174}%
\special{pa 1834 1200}%
\special{pa 1814 1226}%
\special{pa 1796 1250}%
\special{pa 1776 1276}%
\special{pa 1756 1300}%
\special{pa 1736 1326}%
\special{pa 1714 1350}%
\special{pa 1694 1372}%
\special{pa 1672 1396}%
\special{pa 1650 1418}%
\special{pa 1602 1462}%
\special{pa 1580 1482}%
\special{pa 1554 1502}%
\special{pa 1530 1522}%
\special{pa 1426 1594}%
\special{pa 1398 1610}%
\special{pa 1370 1628}%
\special{pa 1342 1644}%
\special{pa 1316 1660}%
\special{pa 1288 1676}%
\special{pa 1260 1694}%
\special{pa 1230 1708}%
\special{pa 1210 1720}%
\special{fp}%
}}%
%
{\color[named]{Black}{%
\special{pn 4}%
\special{pa 2210 920}%
\special{pa 2224 952}%
\special{pa 2238 982}%
\special{pa 2250 1014}%
\special{pa 2292 1104}%
\special{pa 2306 1136}%
\special{pa 2322 1164}%
\special{pa 2336 1194}%
\special{pa 2352 1224}%
\special{pa 2384 1280}%
\special{pa 2400 1306}%
\special{pa 2418 1334}%
\special{pa 2454 1386}%
\special{pa 2472 1410}%
\special{pa 2492 1434}%
\special{pa 2514 1458}%
\special{pa 2534 1480}%
\special{pa 2556 1502}%
\special{pa 2580 1522}%
\special{pa 2604 1540}%
\special{pa 2628 1560}%
\special{pa 2654 1576}%
\special{pa 2680 1594}%
\special{pa 2708 1610}%
\special{pa 2792 1652}%
\special{pa 2822 1666}%
\special{pa 2882 1690}%
\special{pa 2912 1704}%
\special{pa 2944 1714}%
\special{pa 2974 1726}%
\special{pa 3038 1750}%
\special{pa 3040 1750}%
\special{fp}%
}}%
%
{\color[named]{Black}{%
\special{pn 4}%
\special{pa 1280 1970}%
\special{pa 1308 1988}%
\special{pa 1334 2008}%
\special{pa 1360 2026}%
\special{pa 1388 2044}%
\special{pa 1414 2064}%
\special{pa 1440 2082}%
\special{pa 1492 2122}%
\special{pa 1518 2140}%
\special{pa 1542 2160}%
\special{pa 1566 2182}%
\special{pa 1592 2202}%
\special{pa 1614 2224}%
\special{pa 1638 2244}%
\special{pa 1660 2266}%
\special{pa 1682 2290}%
\special{pa 1704 2312}%
\special{pa 1724 2336}%
\special{pa 1744 2362}%
\special{pa 1764 2386}%
\special{pa 1836 2490}%
\special{pa 1900 2602}%
\special{pa 1914 2630}%
\special{pa 1930 2658}%
\special{pa 1944 2688}%
\special{pa 1960 2716}%
\special{pa 1988 2776}%
\special{pa 2010 2820}%
\special{fp}%
}}%
%
{\color[named]{Black}{%
\special{pn 4}%
\special{pa 3060 1930}%
\special{pa 3034 1948}%
\special{pa 3006 1966}%
\special{pa 2978 1982}%
\special{pa 2952 2000}%
\special{pa 2924 2018}%
\special{pa 2898 2036}%
\special{pa 2870 2052}%
\special{pa 2844 2070}%
\special{pa 2816 2088}%
\special{pa 2790 2108}%
\special{pa 2738 2144}%
\special{pa 2714 2164}%
\special{pa 2662 2204}%
\special{pa 2590 2264}%
\special{pa 2566 2286}%
\special{pa 2544 2308}%
\special{pa 2520 2330}%
\special{pa 2498 2352}%
\special{pa 2476 2376}%
\special{pa 2456 2400}%
\special{pa 2434 2424}%
\special{pa 2374 2496}%
\special{pa 2314 2574}%
\special{pa 2296 2598}%
\special{pa 2276 2626}%
\special{pa 2240 2678}%
\special{pa 2220 2704}%
\special{pa 2202 2730}%
\special{pa 2184 2758}%
\special{pa 2166 2784}%
\special{pa 2160 2790}%
\special{fp}%
}}%
%
{\color[named]{Black}{%
\special{pn 4}%
\special{pa 5110 860}%
\special{pa 5090 924}%
\special{pa 5080 954}%
\special{pa 5060 1018}%
\special{pa 5050 1048}%
\special{pa 5040 1080}%
\special{pa 5032 1112}%
\special{pa 5022 1142}%
\special{pa 5012 1174}%
\special{pa 5004 1206}%
\special{pa 4994 1236}%
\special{pa 4978 1300}%
\special{pa 4970 1330}%
\special{pa 4962 1362}%
\special{pa 4956 1394}%
\special{pa 4948 1424}%
\special{pa 4936 1488}%
\special{pa 4930 1518}%
\special{pa 4924 1550}%
\special{pa 4920 1582}%
\special{pa 4916 1612}%
\special{pa 4908 1676}%
\special{pa 4906 1706}%
\special{pa 4904 1738}%
\special{pa 4902 1768}%
\special{pa 4900 1800}%
\special{pa 4900 1862}%
\special{pa 4904 1926}%
\special{pa 4906 1956}%
\special{pa 4908 1988}%
\special{pa 4912 2018}%
\special{pa 4916 2050}%
\special{pa 4920 2080}%
\special{pa 4924 2112}%
\special{pa 4930 2144}%
\special{pa 4936 2174}%
\special{pa 4942 2206}%
\special{pa 4950 2236}%
\special{pa 4956 2268}%
\special{pa 4964 2298}%
\special{pa 4972 2330}%
\special{pa 4982 2362}%
\special{pa 4990 2392}%
\special{pa 4998 2424}%
\special{pa 5008 2454}%
\special{pa 5018 2486}%
\special{pa 5028 2516}%
\special{pa 5038 2548}%
\special{pa 5048 2578}%
\special{pa 5058 2610}%
\special{pa 5068 2640}%
\special{pa 5078 2672}%
\special{pa 5090 2704}%
\special{pa 5100 2734}%
\special{pa 5112 2766}%
\special{pa 5130 2820}%
\special{fp}%
}}%
%
{\color[named]{Black}{%
\special{pn 4}%
\special{pa 5250 870}%
\special{pa 5266 900}%
\special{pa 5280 932}%
\special{pa 5296 962}%
\special{pa 5338 1052}%
\special{pa 5354 1084}%
\special{pa 5368 1114}%
\special{pa 5380 1144}%
\special{pa 5408 1204}%
\special{pa 5420 1236}%
\special{pa 5432 1266}%
\special{pa 5446 1296}%
\special{pa 5456 1326}%
\special{pa 5468 1356}%
\special{pa 5478 1388}%
\special{pa 5508 1478}%
\special{pa 5516 1510}%
\special{pa 5532 1570}%
\special{pa 5544 1630}%
\special{pa 5548 1662}%
\special{pa 5556 1722}%
\special{pa 5558 1752}%
\special{pa 5560 1784}%
\special{pa 5560 1844}%
\special{pa 5556 1904}%
\special{pa 5554 1936}%
\special{pa 5546 1996}%
\special{pa 5540 2026}%
\special{pa 5534 2058}%
\special{pa 5528 2088}%
\special{pa 5512 2148}%
\special{pa 5502 2180}%
\special{pa 5494 2210}%
\special{pa 5484 2240}%
\special{pa 5472 2272}%
\special{pa 5462 2302}%
\special{pa 5438 2362}%
\special{pa 5426 2394}%
\special{pa 5402 2454}%
\special{pa 5388 2484}%
\special{pa 5374 2516}%
\special{pa 5362 2546}%
\special{pa 5334 2606}%
\special{pa 5318 2638}%
\special{pa 5290 2698}%
\special{pa 5276 2730}%
\special{pa 5270 2740}%
\special{fp}%
}}%
%
{\color[named]{Black}{%
\special{pn 4}%
\special{pa 6160 1830}%
\special{pa 6142 1858}%
\special{pa 6118 1878}%
\special{pa 6082 1886}%
\special{pa 6064 1866}%
\special{pa 6060 1832}%
\special{pa 6052 1800}%
\special{pa 6050 1770}%
\special{pa 6068 1742}%
\special{pa 6096 1722}%
\special{pa 6126 1720}%
\special{pa 6166 1732}%
\special{pa 6178 1742}%
\special{pa 6142 1740}%
\special{pa 6140 1740}%
\special{fp}%
}}%
%
{\color[named]{Black}{%
\special{pn 4}%
\special{pa 4140 1760}%
\special{pa 4174 1750}%
\special{pa 4206 1746}%
\special{pa 4236 1748}%
\special{pa 4264 1766}%
\special{pa 4284 1794}%
\special{pa 4290 1828}%
\special{pa 4282 1862}%
\special{pa 4260 1886}%
\special{pa 4230 1892}%
\special{pa 4196 1892}%
\special{pa 4166 1888}%
\special{pa 4144 1866}%
\special{pa 4140 1860}%
\special{fp}%
}}%
%
{\color[named]{Black}{%
\special{pn 4}%
\special{pa 2000 920}%
\special{pa 1930 940}%
\special{fp}%
}}%
%
{\color[named]{Black}{%
\special{pn 4}%
\special{pa 2010 910}%
\special{pa 2040 990}%
\special{fp}%
}}%
%
{\color[named]{Black}{%
\special{pn 4}%
\special{pa 3050 1760}%
\special{pa 2990 1690}%
\special{fp}%
}}%
%
{\color[named]{Black}{%
\special{pn 4}%
\special{pa 3050 1760}%
\special{pa 2990 1770}%
\special{fp}%
}}%
%
{\color[named]{Black}{%
\special{pn 4}%
\special{pa 3070 1940}%
\special{pa 2940 1950}%
\special{fp}%
}}%
%
{\color[named]{Black}{%
\special{pn 4}%
\special{pa 3040 1950}%
\special{pa 3000 2010}%
\special{fp}%
}}%
%
{\color[named]{Black}{%
\special{pn 4}%
\special{pa 2000 2830}%
\special{pa 2000 2720}%
\special{fp}%
}}%
%
{\color[named]{Black}{%
\special{pn 4}%
\special{pa 2000 2820}%
\special{pa 1950 2760}%
\special{fp}%
}}%
%
{\color[named]{Black}{%
\special{pn 4}%
\special{pa 5110 880}%
\special{pa 5020 920}%
\special{fp}%
\special{pa 5110 900}%
\special{pa 5030 930}%
\special{fp}%
}}%
%
{\color[named]{Black}{%
\special{pn 4}%
\special{pa 5120 880}%
\special{pa 5120 970}%
\special{fp}%
}}%
%
{\color[named]{Black}{%
\special{pn 4}%
\special{pa 5260 2720}%
\special{pa 5290 2610}%
\special{fp}%
}}%
%
{\color[named]{Black}{%
\special{pn 4}%
\special{pa 5280 2730}%
\special{pa 5350 2650}%
\special{fp}%
}}%
\put(14.7000,-14.0000){\makebox(0,0)[lb]{$a$}}%
\put(24.8000,-24.8000){\makebox(0,0)[lb]{$a$}}%
\put(25.4000,-13.2000){\makebox(0,0)[lb]{$b$}}%
\put(15.6000,-24.3000){\makebox(0,0)[lb]{$b$}}%
\put(56.4000,-17.8000){\makebox(0,0)[lb]{$b$}}%
\put(47.2000,-18.3000){\makebox(0,0)[lb]{$a$}}%
\put(60.1000,-16.3000){\makebox(0,0)[lb]{$a$}}%
\put(41.6000,-20.5000){\makebox(0,0)[lb]{$b$}}%
\end{picture}%

\end{figure}
The homeomorphism $a$ maps
the long interval $[ba^{-1}b^{-1}(x),a^{-1}b^{-1}(x)]$ onto a subinterval $[x,b^{-1}(x)]$.
Therefore there is a fixed point of $a$ in the open interval
$(x,b^{-1}(x))$. There is also a fixed point in 
$(a^{-1}b^{-1}(x),ba^{-1}b^{-1}(x))$. Likewise $b$ admits at least two
fixed points, one in $(b^{-1}(x),a^{-1}b^{-1}(x))$, another in
$(ba^{-1}b^{-1}(x),x)$.

Let $R$ be the set of four points in Figure 2 left, and set
$S=\{A,A^{-1},B,B^{-1}\}$. Let $R^2=\bigcup_{s\in S}\varphi(s)R$. Then $R^2$
contains $R$, and has 8 more points. The configuration of
$R^2$ in $S^1$ is determined uniquely. Likewise if we set 
$R^3=\bigcup_{s\in S}\varphi(s)R^2$, then its configuration is also unique. 
See Figure 3. 
\begin{figure}\caption{}
\unitlength 0.1in
\begin{picture}( 40.3600, 22.6700)( 14.7400,-29.2000)
%
{\color[named]{Black}{%
\special{pn 8}%
\special{ar 2640 1780 1068 1068  2.1745457  2.1899577}%
}}%
%
{\color[named]{Black}{%
\special{pn 8}%
\special{ar 2530 1850 1056 1056  2.4962118  2.4298863}%
}}%
%
{\color[named]{Black}{%
\special{pn 8}%
\special{ar 2520 2930 2926 2926  5.3910515  5.4058143}%
}}%
%
{\color[named]{Black}{%
\special{pn 8}%
\special{ar 3780 2890 2100 2100  4.7029553  4.7171283}%
}}%
%
{\color[named]{Black}{%
\special{pn 8}%
\special{pa 3780 780}%
\special{pa 3812 788}%
\special{pa 3844 794}%
\special{pa 3874 802}%
\special{pa 3906 808}%
\special{pa 3938 816}%
\special{pa 3968 824}%
\special{pa 4000 832}%
\special{pa 4032 838}%
\special{pa 4062 846}%
\special{pa 4094 856}%
\special{pa 4124 864}%
\special{pa 4156 872}%
\special{pa 4186 882}%
\special{pa 4216 890}%
\special{pa 4248 900}%
\special{pa 4308 920}%
\special{pa 4338 932}%
\special{pa 4368 942}%
\special{pa 4398 954}%
\special{pa 4426 966}%
\special{pa 4486 994}%
\special{pa 4542 1022}%
\special{pa 4572 1036}%
\special{pa 4628 1068}%
\special{pa 4656 1086}%
\special{pa 4682 1104}%
\special{pa 4710 1120}%
\special{pa 4736 1140}%
\special{pa 4762 1158}%
\special{pa 4840 1218}%
\special{pa 4864 1240}%
\special{pa 4890 1260}%
\special{pa 4914 1282}%
\special{pa 4960 1328}%
\special{pa 4982 1352}%
\special{pa 5004 1374}%
\special{pa 5026 1398}%
\special{pa 5048 1424}%
\special{pa 5068 1448}%
\special{pa 5128 1526}%
\special{pa 5164 1578}%
\special{pa 5228 1690}%
\special{pa 5242 1718}%
\special{pa 5256 1748}%
\special{pa 5270 1776}%
\special{pa 5284 1806}%
\special{pa 5320 1896}%
\special{pa 5340 1956}%
\special{pa 5350 1988}%
\special{pa 5370 2048}%
\special{pa 5386 2112}%
\special{pa 5394 2142}%
\special{pa 5400 2174}%
\special{pa 5408 2206}%
\special{pa 5426 2302}%
\special{pa 5430 2334}%
\special{pa 5436 2366}%
\special{pa 5440 2398}%
\special{pa 5444 2428}%
\special{pa 5456 2524}%
\special{pa 5458 2556}%
\special{pa 5462 2588}%
\special{pa 5480 2876}%
\special{pa 5480 2880}%
\special{fp}%
}}%
%
{\color[named]{Black}{%
\special{pn 4}%
\special{sh 1}%
\special{ar 3590 1890 16 16 0  6.28318530717959E+0000}%
\special{sh 1}%
\special{ar 3590 1890 16 16 0  6.28318530717959E+0000}%
}}%
%
{\color[named]{Black}{%
\special{pn 4}%
\special{sh 1}%
\special{ar 2540 790 16 16 0  6.28318530717959E+0000}%
\special{sh 1}%
\special{ar 2540 790 16 16 0  6.28318530717959E+0000}%
}}%
%
{\color[named]{Black}{%
\special{pn 4}%
\special{sh 1}%
\special{ar 1480 1890 16 16 0  6.28318530717959E+0000}%
\special{sh 1}%
\special{ar 1480 1890 16 16 0  6.28318530717959E+0000}%
}}%
%
{\color[named]{Black}{%
\special{pn 4}%
\special{sh 1}%
\special{ar 2540 2920 4 4 0  6.28318530717959E+0000}%
\special{sh 1}%
\special{ar 2540 2920 4 4 0  6.28318530717959E+0000}%
}}%
%
{\color[named]{Black}{%
\special{pn 4}%
\special{pa 3510 1340}%
\special{pa 3420 1400}%
\special{fp}%
}}%
%
{\color[named]{Black}{%
\special{pn 4}%
\special{pa 3110 910}%
\special{pa 3040 1010}%
\special{fp}%
}}%
%
{\color[named]{Black}{%
\special{pn 4}%
\special{pa 1940 890}%
\special{pa 2020 1000}%
\special{fp}%
}}%
%
{\color[named]{Black}{%
\special{pn 4}%
\special{pa 1540 1360}%
\special{pa 1650 1410}%
\special{fp}%
}}%
%
{\color[named]{Black}{%
\special{pn 4}%
\special{pa 1680 2320}%
\special{pa 1540 2420}%
\special{fp}%
}}%
%
{\color[named]{Black}{%
\special{pn 4}%
\special{pa 2180 2790}%
\special{pa 2130 2880}%
\special{fp}%
}}%
%
{\color[named]{Black}{%
\special{pn 4}%
\special{pa 3070 2670}%
\special{pa 3180 2780}%
\special{fp}%
}}%
%
{\color[named]{Black}{%
\special{pn 4}%
\special{sh 1}%
\special{ar 3790 770 16 16 0  6.28318530717959E+0000}%
\special{sh 1}%
\special{ar 3810 790 16 16 0  6.28318530717959E+0000}%
\special{sh 1}%
\special{ar 3810 790 16 16 0  6.28318530717959E+0000}%
}}%
%
{\color[named]{Black}{%
\special{pn 4}%
\special{sh 1}%
\special{ar 5470 2810 16 16 0  6.28318530717959E+0000}%
\special{sh 1}%
\special{ar 5470 2810 16 16 0  6.28318530717959E+0000}%
}}%
%
{\color[named]{Black}{%
\special{pn 20}%
\special{pa 4730 1060}%
\special{pa 4640 1150}%
\special{fp}%
}}%
%
{\color[named]{Black}{%
\special{pn 20}%
\special{pa 5340 1800}%
\special{pa 5230 1890}%
\special{fp}%
}}%
%
{\color[named]{Black}{%
\special{pn 4}%
\special{pa 4170 790}%
\special{pa 4110 920}%
\special{fp}%
}}%
%
{\color[named]{Black}{%
\special{pn 4}%
\special{pa 4500 930}%
\special{pa 4430 1040}%
\special{fp}%
}}%
%
{\color[named]{Black}{%
\special{pn 4}%
\special{pa 5040 1300}%
\special{pa 4910 1410}%
\special{fp}%
}}%
%
{\color[named]{Black}{%
\special{pn 4}%
\special{pa 5200 1560}%
\special{pa 5110 1630}%
\special{fp}%
}}%
%
{\color[named]{Black}{%
\special{pn 4}%
\special{pa 5450 2070}%
\special{pa 5320 2130}%
\special{fp}%
}}%
%
{\color[named]{Black}{%
\special{pn 4}%
\special{pa 5510 2440}%
\special{pa 5390 2440}%
\special{fp}%
}}%
%
{\color[named]{Black}{%
\special{pn 8}%
\special{pa 3360 2300}%
\special{pa 3550 2440}%
\special{fp}%
}}%
%
{\color[named]{Black}{%
\special{pn 8}%
\special{pa 1680 2490}%
\special{pa 1720 2550}%
\special{fp}%
}}%
%
{\color[named]{Black}{%
\special{pn 4}%
\special{sh 1}%
\special{ar 2550 2900 16 16 0  6.28318530717959E+0000}%
\special{sh 1}%
\special{ar 2550 2900 16 16 0  6.28318530717959E+0000}%
}}%
\end{picture}%

\end{figure}
The left depicts $R^2$ and the right a part of $R^3$.
This way, we can determine the configuration of the whole orbit 
$\varphi(\Gamma)x$, which, according to Proposition \ref{equivariant},
 implies that the semiconjugacy class of
$\varphi$ is uniquely determined.
The actual proof can be organized as an induction.

Here is another example of this kind. See Figure 2 right. This is also
a homomorphism $\varphi$ from the free group on two generators $A$ and $B$,
and we denote $a=\varphi(A)$ and $b=\varphi(B)$. The homeomorphism $a$ (resp.\
$b$) has a fixed point $x$ (resp.\ $z$), and we have $y=ab(y)$ for the
point $y$ in the figure. 
Clearly $c(a,b)=1$ and any homomorphism with $c(a,b)=1$ has a
configuration as in Figure 2 right.
Again one can show that such $\varphi$ belongs to a single semiconjugacy
class. That is,
if we let $R$ be the set of four points $x$, $y$,
$z$ and $b(y)$, then the same thing holds with this $R$.

 What is good about these partitions $R$ is the following. Let $\psi$ be
any $k$-fold lift of $\varphi$. Then the pull back image  $\pi_k^{-1}(R)$
has the same property: it determines the
semiconjugacy class of the homomorphism $\psi$.

What is not good is that this kind of partitions are difficult to find out. 
To show Theorem \ref{locally stable}, we need something more.

\bigskip
Let us consider a Fuchsian representation $\varphi\in\RR_{\Pi_2}$ of
the surface group $\Pi_2$ of genus $2$ such that $eu(\varphi)=2$. One
can
assume the elements
$a_\nu=\varphi(A_\nu)$ and $b_\nu=\varphi(B_\nu)$ ($\nu=1,2$) are the
hyperbolic motions in Figure 4 left. 
\begin{figure}\caption{}
\unitlength 0.1in
\begin{picture}( 50.2000, 22.8000)( 13.2000,-29.8000)
%
{\color[named]{Black}{%
\special{pn 8}%
\special{ar 2460 1840 1140 1140  1.8655264  1.8008825}%
}}%
%
{\color[named]{Black}{%
\special{pn 8}%
\special{ar 5140 1860 1108 1108  0.0000000  6.2831853}%
}}%
%
{\color[named]{Black}{%
\special{pn 8}%
\special{pa 2460 850}%
\special{pa 2442 880}%
\special{pa 2388 964}%
\special{pa 2348 1016}%
\special{pa 2328 1040}%
\special{pa 2306 1062}%
\special{pa 2284 1082}%
\special{pa 2260 1100}%
\special{pa 2236 1116}%
\special{pa 2210 1130}%
\special{pa 2154 1150}%
\special{pa 2124 1158}%
\special{pa 2094 1164}%
\special{pa 2062 1168}%
\special{pa 1998 1172}%
\special{pa 1964 1172}%
\special{pa 1862 1166}%
\special{pa 1826 1164}%
\special{pa 1790 1160}%
\special{pa 1756 1156}%
\special{pa 1720 1152}%
\special{pa 1710 1150}%
\special{fp}%
}}%
%
{\color[named]{Black}{%
\special{pn 8}%
\special{pa 1730 1160}%
\special{pa 1730 1296}%
\special{pa 1728 1330}%
\special{pa 1726 1362}%
\special{pa 1718 1426}%
\special{pa 1712 1458}%
\special{pa 1704 1488}%
\special{pa 1684 1544}%
\special{pa 1670 1572}%
\special{pa 1654 1598}%
\special{pa 1636 1622}%
\special{pa 1618 1648}%
\special{pa 1574 1692}%
\special{pa 1526 1736}%
\special{pa 1500 1756}%
\special{pa 1472 1778}%
\special{pa 1446 1798}%
\special{pa 1418 1818}%
\special{pa 1400 1830}%
\special{fp}%
}}%
%
{\color[named]{Black}{%
\special{pn 8}%
\special{pa 1410 1800}%
\special{pa 1440 1820}%
\special{pa 1468 1840}%
\special{pa 1498 1860}%
\special{pa 1526 1880}%
\special{pa 1552 1902}%
\special{pa 1578 1922}%
\special{pa 1604 1944}%
\special{pa 1628 1966}%
\special{pa 1648 1990}%
\special{pa 1668 2012}%
\special{pa 1686 2036}%
\special{pa 1702 2062}%
\special{pa 1716 2088}%
\special{pa 1726 2114}%
\special{pa 1734 2142}%
\special{pa 1742 2202}%
\special{pa 1742 2232}%
\special{pa 1738 2296}%
\special{pa 1734 2328}%
\special{pa 1728 2360}%
\special{pa 1696 2496}%
\special{pa 1686 2530}%
\special{pa 1680 2550}%
\special{fp}%
}}%
%
{\color[named]{Black}{%
\special{pn 8}%
\special{pa 1680 2540}%
\special{pa 1718 2538}%
\special{pa 1754 2536}%
\special{pa 1792 2532}%
\special{pa 1828 2532}%
\special{pa 1864 2530}%
\special{pa 1898 2530}%
\special{pa 1934 2532}%
\special{pa 1966 2534}%
\special{pa 1998 2538}%
\special{pa 2028 2544}%
\special{pa 2058 2552}%
\special{pa 2086 2562}%
\special{pa 2112 2574}%
\special{pa 2136 2590}%
\special{pa 2158 2608}%
\special{pa 2176 2628}%
\special{pa 2194 2652}%
\special{pa 2210 2676}%
\special{pa 2224 2704}%
\special{pa 2238 2734}%
\special{pa 2250 2764}%
\special{pa 2260 2796}%
\special{pa 2290 2898}%
\special{pa 2290 2900}%
\special{fp}%
}}%
%
{\color[named]{Black}{%
\special{pn 8}%
\special{pa 2300 2890}%
\special{pa 2344 2834}%
\special{pa 2364 2806}%
\special{pa 2408 2750}%
\special{pa 2430 2724}%
\special{pa 2452 2700}%
\special{pa 2498 2654}%
\special{pa 2522 2634}%
\special{pa 2548 2616}%
\special{pa 2572 2600}%
\special{pa 2624 2576}%
\special{pa 2652 2568}%
\special{pa 2680 2562}%
\special{pa 2708 2558}%
\special{pa 2768 2558}%
\special{pa 2800 2562}%
\special{pa 2832 2568}%
\special{pa 2862 2574}%
\special{pa 2896 2582}%
\special{pa 2960 2602}%
\special{pa 2994 2612}%
\special{pa 3028 2624}%
\special{pa 3060 2636}%
\special{pa 3094 2648}%
\special{pa 3100 2650}%
\special{fp}%
}}%
%
{\color[named]{Black}{%
\special{pn 8}%
\special{pa 3070 2650}%
\special{pa 3072 2614}%
\special{pa 3074 2580}%
\special{pa 3078 2544}%
\special{pa 3080 2508}%
\special{pa 3092 2406}%
\special{pa 3098 2374}%
\special{pa 3104 2340}%
\special{pa 3112 2310}%
\special{pa 3132 2250}%
\special{pa 3144 2222}%
\special{pa 3156 2196}%
\special{pa 3172 2170}%
\special{pa 3188 2146}%
\special{pa 3228 2102}%
\special{pa 3276 2066}%
\special{pa 3302 2050}%
\special{pa 3330 2034}%
\special{pa 3420 1992}%
\special{pa 3516 1956}%
\special{pa 3530 1950}%
\special{fp}%
}}%
%
{\color[named]{Black}{%
\special{pn 8}%
\special{pa 3520 1940}%
\special{pa 3494 1916}%
\special{pa 3470 1892}%
\special{pa 3444 1870}%
\special{pa 3420 1846}%
\special{pa 3396 1820}%
\special{pa 3372 1796}%
\special{pa 3350 1772}%
\special{pa 3328 1746}%
\special{pa 3308 1722}%
\special{pa 3288 1696}%
\special{pa 3256 1644}%
\special{pa 3242 1616}%
\special{pa 3230 1590}%
\special{pa 3218 1562}%
\special{pa 3210 1532}%
\special{pa 3206 1502}%
\special{pa 3202 1474}%
\special{pa 3200 1442}%
\special{pa 3200 1412}%
\special{pa 3202 1380}%
\special{pa 3210 1316}%
\special{pa 3216 1284}%
\special{pa 3222 1250}%
\special{pa 3230 1218}%
\special{pa 3238 1184}%
\special{pa 3246 1152}%
\special{pa 3250 1130}%
\special{fp}%
}}%
%
{\color[named]{Black}{%
\special{pn 8}%
\special{pa 3250 1160}%
\special{pa 3160 1160}%
\special{pa 3128 1158}%
\special{pa 3098 1158}%
\special{pa 3066 1156}%
\special{pa 3036 1154}%
\special{pa 3004 1150}%
\special{pa 2972 1148}%
\special{pa 2940 1142}%
\special{pa 2908 1138}%
\special{pa 2876 1132}%
\special{pa 2842 1124}%
\special{pa 2810 1116}%
\special{pa 2776 1106}%
\special{pa 2740 1094}%
\special{pa 2706 1082}%
\special{pa 2672 1066}%
\special{pa 2638 1052}%
\special{pa 2606 1034}%
\special{pa 2576 1014}%
\special{pa 2550 994}%
\special{pa 2526 972}%
\special{pa 2504 948}%
\special{pa 2488 920}%
\special{pa 2478 892}%
\special{pa 2472 862}%
\special{pa 2470 840}%
\special{fp}%
}}%
%
{\color[named]{Black}{%
\special{pn 4}%
\special{pa 2230 1230}%
\special{pa 2254 1254}%
\special{pa 2276 1278}%
\special{pa 2300 1302}%
\special{pa 2322 1326}%
\special{pa 2346 1350}%
\special{pa 2442 1438}%
\special{pa 2466 1458}%
\special{pa 2492 1478}%
\special{pa 2544 1514}%
\special{pa 2570 1530}%
\special{pa 2598 1544}%
\special{pa 2624 1558}%
\special{pa 2654 1570}%
\special{pa 2682 1582}%
\special{pa 2712 1590}%
\special{pa 2742 1600}%
\special{pa 2802 1612}%
\special{pa 2834 1618}%
\special{pa 2930 1630}%
\special{pa 2964 1632}%
\special{pa 2996 1634}%
\special{pa 3030 1636}%
\special{pa 3062 1638}%
\special{pa 3096 1638}%
\special{pa 3128 1640}%
\special{pa 3150 1640}%
\special{fp}%
}}%
%
{\color[named]{Black}{%
\special{pn 4}%
\special{pa 2800 1210}%
\special{pa 2748 1250}%
\special{pa 2722 1268}%
\special{pa 2696 1288}%
\special{pa 2668 1308}%
\special{pa 2616 1344}%
\special{pa 2588 1362}%
\special{pa 2562 1380}%
\special{pa 2450 1444}%
\special{pa 2422 1458}%
\special{pa 2394 1470}%
\special{pa 2366 1484}%
\special{pa 2336 1494}%
\special{pa 2306 1506}%
\special{pa 2276 1514}%
\special{pa 2246 1524}%
\special{pa 2216 1532}%
\special{pa 2184 1538}%
\special{pa 2154 1544}%
\special{pa 2122 1550}%
\special{pa 2090 1554}%
\special{pa 2058 1560}%
\special{pa 2026 1562}%
\special{pa 1962 1570}%
\special{pa 1898 1574}%
\special{pa 1864 1576}%
\special{pa 1832 1580}%
\special{pa 1810 1580}%
\special{fp}%
}}%
%
{\color[named]{Black}{%
\special{pn 4}%
\special{pa 1830 2110}%
\special{pa 1858 2108}%
\special{pa 1884 2106}%
\special{pa 1912 2104}%
\special{pa 1968 2104}%
\special{pa 1996 2106}%
\special{pa 2024 2110}%
\special{pa 2054 2114}%
\special{pa 2114 2130}%
\special{pa 2146 2142}%
\special{pa 2178 2158}%
\special{pa 2212 2176}%
\special{pa 2246 2196}%
\special{pa 2282 2222}%
\special{pa 2318 2250}%
\special{pa 2356 2284}%
\special{pa 2396 2320}%
\special{pa 2432 2356}%
\special{pa 2468 2390}%
\special{pa 2500 2422}%
\special{pa 2526 2450}%
\special{pa 2544 2470}%
\special{pa 2556 2480}%
\special{pa 2558 2480}%
\special{pa 2548 2466}%
\special{pa 2528 2442}%
\special{pa 2520 2430}%
\special{fp}%
}}%
%
{\color[named]{Black}{%
\special{pn 4}%
\special{pa 2110 2420}%
\special{pa 2166 2380}%
\special{pa 2192 2360}%
\special{pa 2248 2320}%
\special{pa 2274 2302}%
\special{pa 2302 2282}%
\special{pa 2330 2266}%
\special{pa 2358 2248}%
\special{pa 2386 2232}%
\special{pa 2414 2218}%
\special{pa 2444 2204}%
\special{pa 2472 2192}%
\special{pa 2502 2180}%
\special{pa 2530 2170}%
\special{pa 2560 2162}%
\special{pa 2620 2150}%
\special{pa 2680 2146}%
\special{pa 2712 2146}%
\special{pa 2742 2148}%
\special{pa 2838 2160}%
\special{pa 2870 2166}%
\special{pa 2902 2174}%
\special{pa 2934 2180}%
\special{pa 2998 2196}%
\special{pa 3032 2206}%
\special{pa 3050 2210}%
\special{fp}%
}}%
%
{\color[named]{Black}{%
\special{pn 4}%
\special{pa 1910 1580}%
\special{pa 1820 1570}%
\special{fp}%
\special{sh 1}%
\special{pa 1820 1570}%
\special{pa 1884 1598}%
\special{pa 1874 1576}%
\special{pa 1888 1558}%
\special{pa 1820 1570}%
\special{fp}%
}}%
%
{\color[named]{Black}{%
\special{pn 4}%
\special{pa 3040 1640}%
\special{pa 3160 1650}%
\special{fp}%
\special{sh 1}%
\special{pa 3160 1650}%
\special{pa 3096 1626}%
\special{pa 3108 1646}%
\special{pa 3092 1664}%
\special{pa 3160 1650}%
\special{fp}%
}}%
%
{\color[named]{Black}{%
\special{pn 4}%
\special{pa 1920 2080}%
\special{pa 1880 2110}%
\special{fp}%
\special{sh 1}%
\special{pa 1880 2110}%
\special{pa 1946 2086}%
\special{pa 1924 2078}%
\special{pa 1922 2054}%
\special{pa 1880 2110}%
\special{fp}%
}}%
%
{\color[named]{Black}{%
\special{pn 4}%
\special{pa 2950 2180}%
\special{pa 3050 2230}%
\special{fp}%
\special{sh 1}%
\special{pa 3050 2230}%
\special{pa 3000 2182}%
\special{pa 3002 2206}%
\special{pa 2982 2218}%
\special{pa 3050 2230}%
\special{fp}%
}}%
%
{\color[named]{Black}{%
\special{pn 4}%
\special{pa 1950 2090}%
\special{pa 1850 2090}%
\special{fp}%
\special{sh 1}%
\special{pa 1850 2090}%
\special{pa 1918 2110}%
\special{pa 1904 2090}%
\special{pa 1918 2070}%
\special{pa 1850 2090}%
\special{fp}%
}}%
%
{\color[named]{Black}{%
\special{pn 4}%
\special{pa 4590 910}%
\special{pa 4616 936}%
\special{pa 4642 960}%
\special{pa 4692 1010}%
\special{pa 4770 1082}%
\special{pa 4848 1148}%
\special{pa 4926 1208}%
\special{pa 4954 1226}%
\special{pa 5006 1258}%
\special{pa 5034 1272}%
\special{pa 5060 1286}%
\special{pa 5088 1298}%
\special{pa 5116 1308}%
\special{pa 5144 1316}%
\special{pa 5170 1324}%
\special{pa 5200 1328}%
\special{pa 5228 1332}%
\special{pa 5256 1334}%
\special{pa 5314 1334}%
\special{pa 5344 1330}%
\special{pa 5372 1326}%
\special{pa 5432 1314}%
\special{pa 5462 1304}%
\special{pa 5492 1296}%
\special{pa 5522 1284}%
\special{pa 5552 1274}%
\special{pa 5584 1260}%
\special{pa 5614 1248}%
\special{pa 5644 1234}%
\special{pa 5676 1218}%
\special{pa 5706 1204}%
\special{pa 5738 1188}%
\special{pa 5768 1170}%
\special{pa 5800 1154}%
\special{pa 5830 1136}%
\special{pa 5862 1120}%
\special{pa 5892 1102}%
\special{pa 5924 1084}%
\special{pa 5930 1080}%
\special{fp}%
}}%
%
{\color[named]{Black}{%
\special{pn 4}%
\special{pa 5480 810}%
\special{pa 5440 866}%
\special{pa 5418 892}%
\special{pa 5378 948}%
\special{pa 5356 974}%
\special{pa 5336 1000}%
\special{pa 5314 1026}%
\special{pa 5294 1052}%
\special{pa 5250 1104}%
\special{pa 5184 1176}%
\special{pa 5162 1198}%
\special{pa 5114 1242}%
\special{pa 5092 1262}%
\special{pa 5068 1282}%
\special{pa 5042 1300}%
\special{pa 5018 1318}%
\special{pa 4992 1336}%
\special{pa 4966 1350}%
\special{pa 4940 1366}%
\special{pa 4914 1380}%
\special{pa 4858 1404}%
\special{pa 4802 1424}%
\special{pa 4772 1432}%
\special{pa 4742 1438}%
\special{pa 4712 1446}%
\special{pa 4682 1450}%
\special{pa 4650 1454}%
\special{pa 4620 1458}%
\special{pa 4588 1462}%
\special{pa 4524 1466}%
\special{pa 4492 1466}%
\special{pa 4458 1468}%
\special{pa 4426 1468}%
\special{pa 4392 1466}%
\special{pa 4358 1466}%
\special{pa 4326 1464}%
\special{pa 4292 1462}%
\special{pa 4258 1462}%
\special{pa 4224 1460}%
\special{pa 4190 1456}%
\special{pa 4122 1452}%
\special{pa 4100 1450}%
\special{fp}%
}}%
%
{\color[named]{Black}{%
\special{pn 4}%
\special{pa 4120 2290}%
\special{pa 4154 2292}%
\special{pa 4186 2292}%
\special{pa 4218 2294}%
\special{pa 4252 2294}%
\special{pa 4284 2296}%
\special{pa 4318 2296}%
\special{pa 4414 2302}%
\special{pa 4446 2306}%
\special{pa 4478 2308}%
\special{pa 4606 2324}%
\special{pa 4636 2330}%
\special{pa 4668 2336}%
\special{pa 4698 2342}%
\special{pa 4788 2366}%
\special{pa 4848 2386}%
\special{pa 4876 2398}%
\special{pa 4906 2410}%
\special{pa 4934 2422}%
\special{pa 5018 2464}%
\special{pa 5074 2496}%
\special{pa 5100 2512}%
\special{pa 5128 2530}%
\special{pa 5154 2546}%
\special{pa 5180 2564}%
\special{pa 5208 2582}%
\special{pa 5234 2602}%
\special{pa 5260 2620}%
\special{pa 5364 2700}%
\special{pa 5388 2720}%
\special{pa 5414 2742}%
\special{pa 5440 2762}%
\special{pa 5466 2784}%
\special{pa 5490 2806}%
\special{pa 5516 2826}%
\special{pa 5542 2848}%
\special{pa 5566 2870}%
\special{pa 5590 2890}%
\special{fp}%
}}%
%
{\color[named]{Black}{%
\special{pn 4}%
\special{pa 4680 2870}%
\special{pa 4704 2846}%
\special{pa 4752 2794}%
\special{pa 4774 2770}%
\special{pa 4798 2744}%
\special{pa 4822 2720}%
\special{pa 4846 2694}%
\special{pa 4942 2598}%
\special{pa 4966 2576}%
\special{pa 4992 2552}%
\special{pa 5016 2530}%
\special{pa 5040 2510}%
\special{pa 5066 2488}%
\special{pa 5090 2468}%
\special{pa 5142 2428}%
\special{pa 5220 2374}%
\special{pa 5272 2342}%
\special{pa 5300 2328}%
\special{pa 5326 2314}%
\special{pa 5354 2300}%
\special{pa 5382 2288}%
\special{pa 5438 2268}%
\special{pa 5468 2260}%
\special{pa 5496 2252}%
\special{pa 5526 2244}%
\special{pa 5556 2240}%
\special{pa 5586 2234}%
\special{pa 5618 2232}%
\special{pa 5648 2228}%
\special{pa 5680 2228}%
\special{pa 5710 2226}%
\special{pa 5774 2226}%
\special{pa 5870 2232}%
\special{pa 5904 2236}%
\special{pa 5936 2238}%
\special{pa 5970 2242}%
\special{pa 6002 2246}%
\special{pa 6036 2250}%
\special{pa 6068 2256}%
\special{pa 6102 2260}%
\special{pa 6134 2266}%
\special{pa 6168 2270}%
\special{pa 6170 2270}%
\special{fp}%
}}%
%
{\color[named]{Black}{%
\special{pn 4}%
\special{pa 4040 1840}%
\special{pa 6250 1830}%
\special{fp}%
}}%
%
{\color[named]{Black}{%
\special{pn 8}%
\special{pa 5410 1320}%
\special{pa 5480 1300}%
\special{fp}%
\special{sh 1}%
\special{pa 5480 1300}%
\special{pa 5410 1300}%
\special{pa 5430 1316}%
\special{pa 5422 1338}%
\special{pa 5480 1300}%
\special{fp}%
}}%
%
{\color[named]{Black}{%
\special{pn 8}%
\special{pa 4730 1430}%
\special{pa 4550 1470}%
\special{fp}%
\special{sh 1}%
\special{pa 4550 1470}%
\special{pa 4620 1476}%
\special{pa 4602 1458}%
\special{pa 4612 1436}%
\special{pa 4550 1470}%
\special{fp}%
}}%
%
{\color[named]{Black}{%
\special{pn 8}%
\special{pa 5050 1830}%
\special{pa 5210 1830}%
\special{fp}%
\special{sh 1}%
\special{pa 5210 1830}%
\special{pa 5144 1810}%
\special{pa 5158 1830}%
\special{pa 5144 1850}%
\special{pa 5210 1830}%
\special{fp}%
}}%
%
{\color[named]{Black}{%
\special{pn 8}%
\special{pa 4820 2380}%
\special{pa 4660 2320}%
\special{fp}%
}}%
%
{\color[named]{Black}{%
\special{pn 8}%
\special{pa 5360 2290}%
\special{pa 5520 2230}%
\special{fp}%
\special{sh 1}%
\special{pa 5520 2230}%
\special{pa 5452 2236}%
\special{pa 5470 2250}%
\special{pa 5466 2272}%
\special{pa 5520 2230}%
\special{fp}%
}}%
%
{\color[named]{Black}{%
\special{pn 8}%
\special{pa 4740 2350}%
\special{pa 4680 2330}%
\special{fp}%
\special{sh 1}%
\special{pa 4680 2330}%
\special{pa 4738 2370}%
\special{pa 4732 2348}%
\special{pa 4750 2332}%
\special{pa 4680 2330}%
\special{fp}%
}}%
%
{\color[named]{Black}{%
\special{pn 4}%
\special{sh 1}%
\special{ar 6250 1830 16 16 0  6.28318530717959E+0000}%
\special{sh 1}%
\special{ar 6250 1830 16 16 0  6.28318530717959E+0000}%
}}%
%
{\color[named]{Black}{%
\special{pn 20}%
\special{pa 4000 1820}%
\special{pa 4090 1880}%
\special{fp}%
}}%
%
{\color[named]{Black}{%
\special{pn 20}%
\special{pa 4070 1800}%
\special{pa 4000 1880}%
\special{fp}%
}}%
\put(38.4000,-19.9000){\makebox(0,0)[lb]{$y$}}%
\put(63.4000,-19.1000){\makebox(0,0)[lb]{$x$}}%
\put(54.9000,-14.8000){\makebox(0,0)[lb]{$a_1$}}%
\put(46.2000,-16.3000){\makebox(0,0)[lb]{$b_1$}}%
\put(53.5000,-22.0000){\makebox(0,0)[lb]{$b_2$}}%
\put(46.0000,-22.1000){\makebox(0,0)[lb]{$a_2$}}%
\put(20.3000,-17.6000){\makebox(0,0)[lb]{$b_1$}}%
\put(26.3000,-17.5000){\makebox(0,0)[lb]{$a_1$}}%
\put(20.5000,-20.6000){\makebox(0,0)[lb]{$a_2$}}%
\put(26.2000,-20.4000){\makebox(0,0)[lb]{$b_2$}}%
\put(51.7000,-17.6000){\makebox(0,0)[lb]{$[a_1,b_1]$}}%
\end{picture}%

\end{figure}
The axes of $a_\nu$, $b_\nu$ and
$[a_1,b_1]=[b_2,a_2]$ are depicted in Figure 4 right. Let $x$ and $y$
be the fixed points of $[a_1,b_1]$. See Figure 5 for parts of orbits of
$x$ and $y$. 

\begin{figure}\caption{}
\unitlength 0.1in
\begin{picture}( 43.0000, 47.8000)( 28.9000,-85.3000)
%
{\color[named]{Black}{%
\special{pn 8}%
\special{ar 5100 6260 2090 2090  0.0000000  6.2831853}%
}}%
%
{\color[named]{Black}{%
\special{pn 4}%
\special{sh 1}%
\special{ar 7180 6180 16 16 0  6.28318530717959E+0000}%
\special{sh 1}%
\special{ar 7180 6180 16 16 0  6.28318530717959E+0000}%
}}%
%
{\color[named]{Black}{%
\special{pn 4}%
\special{sh 1}%
\special{ar 6470 4710 16 16 0  6.28318530717959E+0000}%
\special{sh 1}%
\special{ar 6470 4710 16 16 0  6.28318530717959E+0000}%
}}%
%
{\color[named]{Black}{%
\special{pn 4}%
\special{sh 1}%
\special{ar 3990 4490 16 16 0  6.28318530717959E+0000}%
\special{sh 1}%
\special{ar 3990 4490 16 16 0  6.28318530717959E+0000}%
}}%
%
{\color[named]{Black}{%
\special{pn 4}%
\special{sh 1}%
\special{ar 3200 5390 16 16 0  6.28318530717959E+0000}%
\special{sh 1}%
\special{ar 3200 5390 16 16 0  6.28318530717959E+0000}%
}}%
%
{\color[named]{Black}{%
\special{pn 4}%
\special{sh 1}%
\special{ar 3480 7590 16 16 0  6.28318530717959E+0000}%
\special{sh 1}%
\special{ar 3480 7590 16 16 0  6.28318530717959E+0000}%
}}%
%
{\color[named]{Black}{%
\special{pn 4}%
\special{sh 1}%
\special{ar 5190 8360 16 16 0  6.28318530717959E+0000}%
\special{sh 1}%
\special{ar 5190 8360 16 16 0  6.28318530717959E+0000}%
}}%
%
{\color[named]{Black}{%
\special{pn 4}%
\special{sh 1}%
\special{ar 6690 7640 16 16 0  6.28318530717959E+0000}%
\special{sh 1}%
\special{ar 6690 7640 16 16 0  6.28318530717959E+0000}%
}}%
%
{\color[named]{Black}{%
\special{pn 4}%
\special{pa 6780 4910}%
\special{pa 6690 5010}%
\special{fp}%
}}%
%
{\color[named]{Black}{%
\special{pn 4}%
\special{pa 4400 4220}%
\special{pa 4470 4340}%
\special{fp}%
}}%
%
{\color[named]{Black}{%
\special{pn 4}%
\special{pa 3340 5070}%
\special{pa 3420 5140}%
\special{fp}%
}}%
%
{\color[named]{Black}{%
\special{pn 4}%
\special{pa 2890 6250}%
\special{pa 3090 6250}%
\special{fp}%
}}%
%
{\color[named]{Black}{%
\special{pn 4}%
\special{pa 3730 7930}%
\special{pa 3820 7820}%
\special{fp}%
}}%
%
{\color[named]{Black}{%
\special{pn 4}%
\special{pa 5640 8170}%
\special{pa 5720 8320}%
\special{fp}%
}}%
%
{\color[named]{Black}{%
\special{pn 4}%
\special{pa 3790 4860}%
\special{pa 3820 4876}%
\special{pa 3848 4892}%
\special{pa 3876 4906}%
\special{pa 3906 4922}%
\special{pa 3934 4936}%
\special{pa 3964 4952}%
\special{pa 3992 4968}%
\special{pa 4020 4982}%
\special{pa 4050 4998}%
\special{pa 4078 5012}%
\special{pa 4108 5028}%
\special{pa 4136 5042}%
\special{pa 4166 5058}%
\special{pa 4194 5072}%
\special{pa 4224 5086}%
\special{pa 4252 5102}%
\special{pa 4282 5116}%
\special{pa 4310 5130}%
\special{pa 4340 5144}%
\special{pa 4368 5158}%
\special{pa 4398 5172}%
\special{pa 4426 5186}%
\special{pa 4456 5200}%
\special{pa 4484 5214}%
\special{pa 4514 5228}%
\special{pa 4542 5242}%
\special{pa 4572 5254}%
\special{pa 4602 5268}%
\special{pa 4630 5282}%
\special{pa 4690 5306}%
\special{pa 4720 5320}%
\special{pa 4748 5332}%
\special{pa 4838 5368}%
\special{pa 4868 5378}%
\special{pa 4896 5390}%
\special{pa 4926 5402}%
\special{pa 4986 5422}%
\special{pa 5016 5434}%
\special{pa 5076 5454}%
\special{pa 5106 5462}%
\special{pa 5166 5482}%
\special{pa 5198 5490}%
\special{pa 5228 5500}%
\special{pa 5288 5516}%
\special{pa 5320 5524}%
\special{pa 5350 5530}%
\special{pa 5380 5538}%
\special{pa 5412 5546}%
\special{pa 5472 5558}%
\special{pa 5536 5570}%
\special{pa 5566 5576}%
\special{pa 5598 5580}%
\special{pa 5628 5586}%
\special{pa 5692 5594}%
\special{pa 5722 5598}%
\special{pa 5850 5614}%
\special{pa 5880 5616}%
\special{pa 5912 5618}%
\special{pa 5944 5622}%
\special{pa 6168 5636}%
\special{pa 6200 5636}%
\special{pa 6232 5638}%
\special{pa 6264 5638}%
\special{pa 6296 5640}%
\special{pa 6362 5640}%
\special{pa 6394 5642}%
\special{pa 6848 5642}%
\special{pa 6880 5640}%
\special{pa 6910 5640}%
\special{fp}%
}}%
%
{\color[named]{Black}{%
\special{pn 4}%
\special{pa 5770 4410}%
\special{pa 5750 4438}%
\special{pa 5730 4464}%
\special{pa 5710 4492}%
\special{pa 5690 4518}%
\special{pa 5670 4546}%
\special{pa 5650 4572}%
\special{pa 5628 4600}%
\special{pa 5588 4652}%
\special{pa 5568 4680}%
\special{pa 5548 4706}%
\special{pa 5526 4732}%
\special{pa 5506 4760}%
\special{pa 5486 4786}%
\special{pa 5464 4812}%
\special{pa 5424 4864}%
\special{pa 5402 4888}%
\special{pa 5382 4914}%
\special{pa 5360 4940}%
\special{pa 5338 4964}%
\special{pa 5318 4990}%
\special{pa 5296 5014}%
\special{pa 5274 5040}%
\special{pa 5252 5064}%
\special{pa 5232 5088}%
\special{pa 5188 5136}%
\special{pa 5164 5158}%
\special{pa 5142 5182}%
\special{pa 5120 5204}%
\special{pa 5098 5228}%
\special{pa 5074 5250}%
\special{pa 5052 5272}%
\special{pa 5028 5294}%
\special{pa 5006 5314}%
\special{pa 4982 5336}%
\special{pa 4862 5436}%
\special{pa 4836 5454}%
\special{pa 4812 5474}%
\special{pa 4786 5492}%
\special{pa 4762 5510}%
\special{pa 4736 5526}%
\special{pa 4710 5544}%
\special{pa 4632 5592}%
\special{pa 4604 5606}%
\special{pa 4578 5622}%
\special{pa 4522 5650}%
\special{pa 4494 5662}%
\special{pa 4466 5676}%
\special{pa 4410 5700}%
\special{pa 4382 5710}%
\special{pa 4352 5722}%
\special{pa 4322 5732}%
\special{pa 4294 5742}%
\special{pa 4234 5762}%
\special{pa 4204 5770}%
\special{pa 4174 5780}%
\special{pa 4142 5788}%
\special{pa 4082 5804}%
\special{pa 4050 5810}%
\special{pa 4020 5818}%
\special{pa 3988 5824}%
\special{pa 3956 5832}%
\special{pa 3924 5838}%
\special{pa 3894 5844}%
\special{pa 3830 5856}%
\special{pa 3798 5860}%
\special{pa 3764 5866}%
\special{pa 3732 5870}%
\special{pa 3700 5876}%
\special{pa 3636 5884}%
\special{pa 3602 5890}%
\special{pa 3570 5894}%
\special{pa 3536 5898}%
\special{pa 3472 5906}%
\special{pa 3438 5910}%
\special{pa 3406 5914}%
\special{pa 3372 5916}%
\special{pa 3340 5920}%
\special{pa 3272 5928}%
\special{pa 3250 5930}%
\special{fp}%
}}%
%
{\color[named]{Black}{%
\special{pn 4}%
\special{pa 4450 8040}%
\special{pa 4510 7956}%
\special{pa 4530 7930}%
\special{pa 4610 7818}%
\special{pa 4630 7792}%
\special{pa 4650 7764}%
\special{pa 4672 7736}%
\special{pa 4712 7684}%
\special{pa 4732 7656}%
\special{pa 4754 7630}%
\special{pa 4794 7578}%
\special{pa 4816 7552}%
\special{pa 4838 7528}%
\special{pa 4858 7502}%
\special{pa 4924 7430}%
\special{pa 4944 7406}%
\special{pa 4990 7360}%
\special{pa 5012 7336}%
\special{pa 5034 7314}%
\special{pa 5058 7292}%
\special{pa 5080 7272}%
\special{pa 5104 7250}%
\special{pa 5128 7230}%
\special{pa 5150 7210}%
\special{pa 5174 7190}%
\special{pa 5200 7172}%
\special{pa 5224 7152}%
\special{pa 5248 7134}%
\special{pa 5274 7118}%
\special{pa 5298 7100}%
\special{pa 5324 7084}%
\special{pa 5350 7070}%
\special{pa 5376 7054}%
\special{pa 5402 7040}%
\special{pa 5430 7026}%
\special{pa 5456 7014}%
\special{pa 5484 7000}%
\special{pa 5512 6988}%
\special{pa 5596 6958}%
\special{pa 5626 6948}%
\special{pa 5654 6940}%
\special{pa 5744 6916}%
\special{pa 5834 6898}%
\special{pa 5866 6892}%
\special{pa 5896 6888}%
\special{pa 5960 6880}%
\special{pa 5990 6876}%
\special{pa 6022 6874}%
\special{pa 6054 6870}%
\special{pa 6086 6868}%
\special{pa 6120 6866}%
\special{pa 6152 6864}%
\special{pa 6184 6864}%
\special{pa 6218 6862}%
\special{pa 6250 6862}%
\special{pa 6284 6860}%
\special{pa 6384 6860}%
\special{pa 6418 6862}%
\special{pa 6484 6862}%
\special{pa 6552 6866}%
\special{pa 6586 6866}%
\special{pa 6654 6870}%
\special{pa 6688 6870}%
\special{pa 6756 6874}%
\special{pa 6792 6876}%
\special{pa 6860 6880}%
\special{pa 6870 6880}%
\special{fp}%
}}%
%
{\color[named]{Black}{%
\special{pn 4}%
\special{pa 6170 7810}%
\special{pa 6118 7770}%
\special{pa 6090 7748}%
\special{pa 6012 7688}%
\special{pa 5984 7668}%
\special{pa 5932 7628}%
\special{pa 5904 7608}%
\special{pa 5852 7568}%
\special{pa 5824 7548}%
\special{pa 5798 7528}%
\special{pa 5770 7510}%
\special{pa 5744 7490}%
\special{pa 5716 7472}%
\special{pa 5690 7452}%
\special{pa 5662 7434}%
\special{pa 5636 7416}%
\special{pa 5608 7398}%
\special{pa 5582 7380}%
\special{pa 5554 7362}%
\special{pa 5526 7346}%
\special{pa 5500 7328}%
\special{pa 5388 7264}%
\special{pa 5360 7250}%
\special{pa 5334 7234}%
\special{pa 5306 7220}%
\special{pa 5276 7206}%
\special{pa 5220 7178}%
\special{pa 5164 7154}%
\special{pa 5134 7142}%
\special{pa 5078 7118}%
\special{pa 5018 7098}%
\special{pa 4990 7088}%
\special{pa 4930 7072}%
\special{pa 4902 7064}%
\special{pa 4872 7056}%
\special{pa 4842 7050}%
\special{pa 4812 7042}%
\special{pa 4782 7036}%
\special{pa 4750 7032}%
\special{pa 4720 7028}%
\special{pa 4690 7022}%
\special{pa 4660 7020}%
\special{pa 4628 7016}%
\special{pa 4598 7012}%
\special{pa 4566 7010}%
\special{pa 4536 7008}%
\special{pa 4504 7008}%
\special{pa 4472 7006}%
\special{pa 4314 7006}%
\special{pa 4250 7010}%
\special{pa 4218 7010}%
\special{pa 4186 7014}%
\special{pa 4122 7018}%
\special{pa 4090 7022}%
\special{pa 4058 7024}%
\special{pa 4026 7028}%
\special{pa 3992 7032}%
\special{pa 3928 7040}%
\special{pa 3894 7044}%
\special{pa 3862 7048}%
\special{pa 3830 7054}%
\special{pa 3796 7058}%
\special{pa 3764 7064}%
\special{pa 3730 7068}%
\special{pa 3698 7074}%
\special{pa 3664 7078}%
\special{pa 3632 7084}%
\special{pa 3598 7090}%
\special{pa 3566 7096}%
\special{pa 3532 7102}%
\special{pa 3500 7108}%
\special{pa 3466 7114}%
\special{pa 3434 7120}%
\special{pa 3400 7126}%
\special{pa 3370 7130}%
\special{fp}%
}}%
%
{\color[named]{Black}{%
\special{pn 4}%
\special{pa 6820 7300}%
\special{pa 6990 7370}%
\special{fp}%
}}%
%
{\color[named]{Black}{%
\special{pn 20}%
\special{pa 6800 5650}%
\special{pa 6910 5630}%
\special{fp}%
\special{sh 1}%
\special{pa 6910 5630}%
\special{pa 6842 5622}%
\special{pa 6858 5640}%
\special{pa 6848 5662}%
\special{pa 6910 5630}%
\special{fp}%
}}%
%
{\color[named]{Black}{%
\special{pn 13}%
\special{pa 3370 5900}%
\special{pa 3260 5920}%
\special{fp}%
\special{sh 1}%
\special{pa 3260 5920}%
\special{pa 3330 5928}%
\special{pa 3312 5910}%
\special{pa 3322 5888}%
\special{pa 3260 5920}%
\special{fp}%
}}%
%
{\color[named]{Black}{%
\special{pn 13}%
\special{pa 6780 6880}%
\special{pa 6910 6880}%
\special{fp}%
\special{sh 1}%
\special{pa 6910 6880}%
\special{pa 6844 6860}%
\special{pa 6858 6880}%
\special{pa 6844 6900}%
\special{pa 6910 6880}%
\special{fp}%
}}%
%
{\color[named]{Black}{%
\special{pn 13}%
\special{pa 3470 7110}%
\special{pa 3400 7120}%
\special{fp}%
\special{sh 1}%
\special{pa 3400 7120}%
\special{pa 3470 7130}%
\special{pa 3454 7112}%
\special{pa 3464 7092}%
\special{pa 3400 7120}%
\special{fp}%
}}%
\put(32.0000,-63.6000){\makebox(0,0)[lb]{$y=[a_1,b_1](y)=[a_2,b_2](y)$}}%
\put(56.6000,-54.8000){\makebox(0,0)[lb]{$a_1$}}%
\put(40.1000,-57.3000){\makebox(0,0)[lb]{$b_1$}}%
\put(60.0000,-70.9000){\makebox(0,0)[lb]{$b_2$}}%
\put(41.3000,-73.2000){\makebox(0,0)[lb]{$a_2$}}%
\put(66.3000,-46.8000){\makebox(0,0)[lb]{$b_1^{-1}(x)$}}%
\put(69.0000,-50.0000){\makebox(0,0)[lb]{$b_1^{-1}(y)$}}%
\put(44.1000,-40.8000){\makebox(0,0)[lb]{$a_1^{-1}b_1^{-1}(y)$}}%
\put(36.6000,-43.7000){\makebox(0,0)[lb]{$a_1^{-1}b_1^{-1}(x)$}}%
\put(34.7000,-52.8000){\makebox(0,0)[lb]{$b_1a_1^{-1}b_1^{-1}(y)$}}%
\put(32.5000,-55.7000){\makebox(0,0)[lb]{$b_1a_1^{-1}b_1^{-1}(x)$}}%
\put(70.7000,-74.6000){\makebox(0,0)[lb]{$b_2a_2^{-1}b_2^{-1}(y)$}}%
\put(67.8000,-77.4000){\makebox(0,0)[lb]{$b_2a_2^{-1}b_2^{-1}(x)$}}%
\put(58.5000,-84.9000){\makebox(0,0)[lb]{$a_2^{-1}b_2^{-1}(y)$}}%
\put(47.8000,-86.6000){\makebox(0,0)[lb]{$a_2^{-1}b_2^{-1}(x)$}}%
\put(36.5000,-82.4000){\makebox(0,0)[lb]{$b_2^{-1}(y)$}}%
\put(30.7000,-78.3000){\makebox(0,0)[lb]{$b_2^{-1}(x)$}}%
\put(55.1000,-62.8000){\makebox(0,0)[lb]{$[a_1,b_1](x)=[a_2,b_2](x)=x$}}%
\end{picture}%

\end{figure}
The set $R$ of fourteen points there is enough to determine
the semiconjugacy class of the homomorphism $\varphi$. In fact, the
configuration of $R$ immediately implies that $eu(\varphi)=2$, and by
\cite{Matsumoto2}, the semiconjugacy class is unique.  
However 
when we consider a $2$-fold lift $\psi$ of $\varphi$, 
it is not clear if the inverse image $\pi_2^{-1}(R)$ actually
determines the semiconjugacy class or not.
To cope with the problem, we need an algorythm to determine
the orbits of $x$ and $y$, which can be inherited to a $k$-fold cover. 
{\em But this is not according to the word length
of the elements of $\Pi_2$.}

\medskip
Consider the amalgamated product $$\Pi_2=\Gamma_1*_\Lambda\Gamma_2,$$
where $\Gamma_\nu$ is the subgroup generated by $A_\nu$ and $B_\nu$ 
and $\Lambda$ generated by $[A_1,B_1]=[B_2,A_2]$.
First we consider the homomorphism $\varphi_\nu=\varphi\vert_{\Gamma_\nu}$. This is a
homomorphism from the free group $\Gamma_\nu$ on two generators
such that $\tau([\widetilde a_\nu,\widetilde b_\nu])=1$, and the previous
observation works. However notice that one can define
 the set $R$ of four points in Figure 2 in two different ways:
one from the orbit of $x$, the other $y$.
It is more natural and more convenient
to consider disjoint four intervals (instead
of points). For $\Gamma_1$, they are $E_1=[y,x]$ and its iterates  
in Figure 6. The complement of the four intervals is denoted by $P_1$.
The stabilizer (in $\Gamma_1$) of $E_1$ is $\Lambda$, and the limit set
of the Fuchsian group $\Gamma_1$ is contained in $P_1$.
Likewise in the right figure, the four intervals are
$E_2=[x,y]$ and its iterates. The complement is denoted by $P_2$. 

For $\gamma_\nu\in\Gamma_\nu$, the configuration of the
orbits $\varphi_\nu(\gamma_\nu)(x)$ and $\varphi_\nu(\gamma_\nu)(y)$
is determined just by the data in Figure 6
inductively on the word length of $\gamma_\nu$, as we have explained. 
They are contained in $P_\nu$.
\begin{figure}\caption{}
\unitlength 0.1in
\begin{picture}( 49.8800, 24.0000)( 13.5000,-34.8000)
%
{\color[named]{Black}{%
\special{pn 8}%
\special{ar 2560 2290 1150 1150  0.0000000  6.2831853}%
}}%
%
{\color[named]{Black}{%
\special{pn 4}%
\special{sh 1}%
\special{ar 3710 2280 8 8 0  6.28318530717959E+0000}%
\special{sh 1}%
\special{ar 3710 2280 8 8 0  6.28318530717959E+0000}%
}}%
%
{\color[named]{Black}{%
\special{pn 4}%
\special{sh 1}%
\special{ar 3350 1460 16 16 0  6.28318530717959E+0000}%
\special{sh 1}%
\special{ar 3350 1460 16 16 0  6.28318530717959E+0000}%
}}%
%
{\color[named]{Black}{%
\special{pn 4}%
\special{sh 1}%
\special{ar 2100 1240 16 16 0  6.28318530717959E+0000}%
\special{sh 1}%
\special{ar 2100 1240 16 16 0  6.28318530717959E+0000}%
}}%
%
{\color[named]{Black}{%
\special{pn 4}%
\special{sh 1}%
\special{ar 1500 1850 16 16 0  6.28318530717959E+0000}%
\special{sh 1}%
\special{ar 1500 1850 16 16 0  6.28318530717959E+0000}%
}}%
%
{\color[named]{Black}{%
\special{pn 4}%
\special{pa 3520 1760}%
\special{pa 3610 1710}%
\special{fp}%
}}%
%
{\color[named]{Black}{%
\special{pn 4}%
\special{pa 2570 1200}%
\special{pa 2570 1080}%
\special{fp}%
}}%
%
{\color[named]{Black}{%
\special{pn 4}%
\special{pa 1640 1520}%
\special{pa 1730 1600}%
\special{fp}%
}}%
%
{\color[named]{Black}{%
\special{pn 4}%
\special{pa 1350 2290}%
\special{pa 1480 2290}%
\special{fp}%
}}%
%
{\color[named]{Black}{%
\special{pn 20}%
\special{ar 2550 2280 1150 1150  5.7962861  6.2831853}%
}}%
%
{\color[named]{Black}{%
\special{pn 20}%
\special{ar 2560 2290 1146 1146  4.7210844  5.4731008}%
}}%
%
{\color[named]{Black}{%
\special{pn 20}%
\special{ar 2560 2280 1142 1142  3.8329296  4.2879370}%
}}%
%
{\color[named]{Black}{%
\special{pn 20}%
\special{ar 2560 2290 1144 1144  3.1502881  3.5269691}%
}}%
%
{\color[named]{Black}{%
\special{pn 8}%
\special{ar 5180 2290 1150 1150  0.0000000  6.2831853}%
}}%
%
{\color[named]{Black}{%
\special{pn 4}%
\special{sh 1}%
\special{ar 6330 2280 16 16 0  6.28318530717959E+0000}%
\special{sh 1}%
\special{ar 6330 2280 16 16 0  6.28318530717959E+0000}%
}}%
%
{\color[named]{Black}{%
\special{pn 4}%
\special{sh 1}%
\special{ar 6020 3090 16 16 0  6.28318530717959E+0000}%
\special{sh 1}%
\special{ar 6020 3090 16 16 0  6.28318530717959E+0000}%
}}%
%
{\color[named]{Black}{%
\special{pn 4}%
\special{sh 1}%
\special{ar 5180 3440 16 16 0  6.28318530717959E+0000}%
\special{sh 1}%
\special{ar 5180 3440 16 16 0  6.28318530717959E+0000}%
}}%
%
{\color[named]{Black}{%
\special{pn 4}%
\special{sh 1}%
\special{ar 4290 3010 16 16 0  6.28318530717959E+0000}%
\special{sh 1}%
\special{ar 4290 3010 16 16 0  6.28318530717959E+0000}%
}}%
%
{\color[named]{Black}{%
\special{pn 4}%
\special{pa 6110 2860}%
\special{pa 6200 2920}%
\special{fp}%
}}%
%
{\color[named]{Black}{%
\special{pn 4}%
\special{pa 5480 3360}%
\special{pa 5520 3480}%
\special{fp}%
}}%
%
{\color[named]{Black}{%
\special{pn 4}%
\special{pa 4470 3120}%
\special{pa 4390 3210}%
\special{fp}%
}}%
%
{\color[named]{Black}{%
\special{pn 4}%
\special{pa 3940 2290}%
\special{pa 4120 2290}%
\special{fp}%
}}%
%
{\color[named]{Black}{%
\special{pn 20}%
\special{ar 5160 2280 1140 1140  2.4490940  3.1415927}%
}}%
%
{\color[named]{Black}{%
\special{pn 20}%
\special{ar 5170 2280 1150 1150  1.5622495  2.2699889}%
}}%
%
{\color[named]{Black}{%
\special{pn 20}%
\special{ar 5170 2280 1158 1158  0.7607118  1.2818209}%
}}%
%
{\color[named]{Black}{%
\special{pn 20}%
\special{ar 5170 2280 1168 1168  6.2831853  0.5595068}%
}}%
\put(16.0000,-23.6000){\makebox(0,0)[lb]{$y$}}%
\put(34.8000,-23.6000){\makebox(0,0)[lb]{$x$}}%
\put(41.8000,-23.6000){\makebox(0,0)[lb]{$y$}}%
\put(61.4000,-23.4000){\makebox(0,0)[lb]{$x$}}%
\put(24.8000,-19.7000){\makebox(0,0)[lb]{$P_1$}}%
\put(50.5000,-27.9000){\makebox(0,0)[lb]{$P_2$}}%
%
{\color[named]{Black}{%
\special{pn 4}%
\special{pa 2880 1910}%
\special{pa 3540 1990}%
\special{fp}%
\special{sh 1}%
\special{pa 3540 1990}%
\special{pa 3476 1962}%
\special{pa 3488 1984}%
\special{pa 3472 2002}%
\special{pa 3540 1990}%
\special{fp}%
}}%
%
{\color[named]{Black}{%
\special{pn 4}%
\special{pa 2660 1760}%
\special{pa 2970 1360}%
\special{fp}%
\special{sh 1}%
\special{pa 2970 1360}%
\special{pa 2914 1400}%
\special{pa 2938 1402}%
\special{pa 2946 1426}%
\special{pa 2970 1360}%
\special{fp}%
}}%
%
{\color[named]{Black}{%
\special{pn 4}%
\special{pa 2350 1770}%
\special{pa 2040 1470}%
\special{fp}%
\special{sh 1}%
\special{pa 2040 1470}%
\special{pa 2074 1532}%
\special{pa 2078 1508}%
\special{pa 2102 1502}%
\special{pa 2040 1470}%
\special{fp}%
}}%
%
{\color[named]{Black}{%
\special{pn 4}%
\special{pa 2330 1950}%
\special{pa 1630 2070}%
\special{fp}%
\special{sh 1}%
\special{pa 1630 2070}%
\special{pa 1700 2078}%
\special{pa 1684 2062}%
\special{pa 1692 2040}%
\special{pa 1630 2070}%
\special{fp}%
}}%
%
{\color[named]{Black}{%
\special{pn 4}%
\special{pa 5440 2720}%
\special{pa 6120 2610}%
\special{fp}%
\special{sh 1}%
\special{pa 6120 2610}%
\special{pa 6052 2602}%
\special{pa 6068 2620}%
\special{pa 6058 2640}%
\special{pa 6120 2610}%
\special{fp}%
}}%
%
{\color[named]{Black}{%
\special{pn 4}%
\special{pa 5310 2840}%
\special{pa 5700 3190}%
\special{fp}%
\special{sh 1}%
\special{pa 5700 3190}%
\special{pa 5664 3132}%
\special{pa 5660 3154}%
\special{pa 5638 3160}%
\special{pa 5700 3190}%
\special{fp}%
}}%
%
{\color[named]{Black}{%
\special{pn 4}%
\special{pa 5010 2870}%
\special{pa 4740 3220}%
\special{fp}%
\special{sh 1}%
\special{pa 4740 3220}%
\special{pa 4798 3180}%
\special{pa 4774 3178}%
\special{pa 4766 3156}%
\special{pa 4740 3220}%
\special{fp}%
}}%
%
{\color[named]{Black}{%
\special{pn 4}%
\special{pa 4910 2760}%
\special{pa 4230 2670}%
\special{fp}%
\special{sh 1}%
\special{pa 4230 2670}%
\special{pa 4294 2700}%
\special{pa 4284 2678}%
\special{pa 4300 2660}%
\special{pa 4230 2670}%
\special{fp}%
}}%
\put(23.6000,-33.0000){\makebox(0,0)[lb]{$E_1$}}%
\put(49.4000,-13.8000){\makebox(0,0)[lb]{$E_2$}}%
\end{picture}%

\end{figure}
The actual proof is given in Section 4,
where we call such subsets $P_\nu$ basic partitions. The complementary
intervals $E_\nu$ is called the entrance of $\Lambda$ to $P_\nu$.
As we explained, the stabilizer of $E_\nu$ (in $\Gamma_\nu$) is $\Lambda$.
The entrances $E_1$ and $E_2$ satisfy the conditions; 
$E_1\cup E_2=S^1$ and $\Int E_1\cap \Int E_2=\emptyset$. They are said to
be combinable. Now the whole orbits $\varphi(g)(x)$ and $\varphi(g)(y)$ for
$g\in\Pi_2$ can be
determined just by this combinability condition. This part, reminiscent of the Maskit
combination theorem \cite{Maskit} in Kleinian groups, is shown in
Section 5.
What is good for this construction is that the whole process
 can be passed to a 2-fold
lift of $\varphi$. 

Moreover the set $R$ of fourteen points in Figure 5 are robust, in the sense
that any homomorphism near to $\varphi$ has the same configuration
as $R$. Furthermore, if we consider a 2-fold lift $\psi$ of $\varphi$,
The set $\pi_2^{-1}(R)$ is also robust for $\psi$.
 This part, shown in Section
7, concludes the proof of the local stability
(Theorem \ref{locally stable}) for $g=2$.

For $g\geq3$, the group $\Pi_g$ is represented as the fundamental group
of a tree of groups. Each vertex of the tree has valency either 1 or
3. For a valency 3 vertex, we have a homomorphism $\varphi\in\RR_\Gamma$,
where $\Gamma$ is the free group on two generators $A$ and $B$.
The homeomorphism $a=\varphi(A)$ and $b=\varphi(B)$ has the property that
$c(a,b)=1$. This implies that $\varphi$ admits a
configuration in Figure 2 right. For this we consider a basic partition
$P$ as in Figure 7. 
\begin{figure}\caption{}
\unitlength 0.1in
\begin{picture}( 36.8000, 30.7000)( 15.7000,-36.7000)
%
{\color[named]{Black}{%
\special{pn 8}%
\special{ar 3710 2340 1380 1380  0.0000000  6.2831853}%
}}%
%
{\color[named]{Black}{%
\special{pn 4}%
\special{sh 1}%
\special{ar 4460 1190 16 16 0  6.28318530717959E+0000}%
\special{sh 1}%
\special{ar 4460 1190 16 16 0  6.28318530717959E+0000}%
}}%
%
{\color[named]{Black}{%
\special{pn 4}%
\special{sh 1}%
\special{ar 2620 1510 16 16 0  6.28318530717959E+0000}%
\special{sh 1}%
\special{ar 2620 1510 16 16 0  6.28318530717959E+0000}%
}}%
%
{\color[named]{Black}{%
\special{pn 4}%
\special{sh 1}%
\special{ar 2710 3280 16 16 0  6.28318530717959E+0000}%
\special{sh 1}%
\special{ar 2710 3280 16 16 0  6.28318530717959E+0000}%
}}%
%
{\color[named]{Black}{%
\special{pn 4}%
\special{sh 1}%
\special{ar 4820 3150 16 16 0  6.28318530717959E+0000}%
\special{sh 1}%
\special{ar 4820 3150 16 16 0  6.28318530717959E+0000}%
}}%
%
{\color[named]{Black}{%
\special{pn 4}%
\special{pa 4950 1870}%
\special{pa 5070 1800}%
\special{fp}%
}}%
%
{\color[named]{Black}{%
\special{pn 4}%
\special{pa 3600 1030}%
\special{pa 3600 900}%
\special{fp}%
}}%
%
{\color[named]{Black}{%
\special{pn 4}%
\special{pa 2260 2350}%
\special{pa 2410 2350}%
\special{fp}%
}}%
%
{\color[named]{Black}{%
\special{pn 4}%
\special{pa 3600 3650}%
\special{pa 3600 3770}%
\special{fp}%
}}%
%
{\color[named]{Black}{%
\special{pn 20}%
\special{ar 3730 2370 1376 1376  5.2763481  5.8878352}%
}}%
%
{\color[named]{Black}{%
\special{pn 20}%
\special{ar 3720 2370 1406 1406  3.8107657  4.6268837}%
}}%
%
{\color[named]{Black}{%
\special{pn 20}%
\special{ar 3720 2370 1380 1380  2.4131277  3.1560844}%
}}%
%
{\color[named]{Black}{%
\special{pn 20}%
\special{ar 3720 2370 1354 1354  0.6254850  1.6520982}%
}}%
\put(52.5000,-25.6000){\makebox(0,0)[lb]{$E_1=a(E_1)$}}%
\put(36.8000,-8.3000){\makebox(0,0)[lb]{$E_2=a(E_4)=ab(E_2)$}}%
\put(15.7000,-20.4000){\makebox(0,0)[lb]{$E_3=b(E_3)$}}%
\put(18.3000,-37.3000){\makebox(0,0)[lb]{$E_4=b(E_2)=ba(E_4)$}}%
\put(36.8000,-24.0000){\makebox(0,0)[lb]{$P$}}%
%
{\color[named]{Black}{%
\special{pn 4}%
\special{pa 3930 2230}%
\special{pa 4660 1600}%
\special{fp}%
\special{sh 1}%
\special{pa 4660 1600}%
\special{pa 4596 1628}%
\special{pa 4620 1636}%
\special{pa 4624 1660}%
\special{pa 4660 1600}%
\special{fp}%
}}%
%
{\color[named]{Black}{%
\special{pn 4}%
\special{pa 3580 2360}%
\special{pa 2530 2820}%
\special{fp}%
\special{sh 1}%
\special{pa 2530 2820}%
\special{pa 2600 2812}%
\special{pa 2580 2800}%
\special{pa 2584 2776}%
\special{pa 2530 2820}%
\special{fp}%
}}%
%
{\color[named]{Black}{%
\special{pn 4}%
\special{pa 3790 2490}%
\special{pa 4220 3400}%
\special{fp}%
\special{sh 1}%
\special{pa 4220 3400}%
\special{pa 4210 3332}%
\special{pa 4198 3352}%
\special{pa 4174 3348}%
\special{pa 4220 3400}%
\special{fp}%
}}%
%
{\color[named]{Black}{%
\special{pn 4}%
\special{pa 3700 2200}%
\special{pa 3170 1260}%
\special{fp}%
\special{sh 1}%
\special{pa 3170 1260}%
\special{pa 3186 1328}%
\special{pa 3196 1306}%
\special{pa 3220 1308}%
\special{pa 3170 1260}%
\special{fp}%
}}%
\end{picture}%

\end{figure}
The complementary region consists of four intervals
$E_1$--$E_4$. The stabilizer of $E_1$ is the subgroup $\langle a\rangle$,
and we say that $E_1$ is the entrance of $\langle a\rangle$ to $P$.
Likewise $E_2$, $E_3$ and $E_4$ are entraces to $P$ of the subgroups
$\langle ab\rangle$, $\langle b\rangle$ and $\langle ba\rangle$, respectively.
Compare with Figure 2 right.

\section{Basic partitions}

Let $\Gamma$ be a group with a prescribed finite symmetric generating set $S$.

\begin{definition} \label{d1}
A subset $P$ of $S^1$ is called a {\em basic partition} (BP) for
 $\varphi\in\RR_\Gamma$, if it satisfies the following conditions. 

(1) $P$ is a union of finitely many disjoint closed intervals. 

(2) For any $I\sqsubset P$, there
 exists a unique element $s_I\in S$ such that
$$
\varphi(s_I)I=\bigcup_{i=1}^lI_i \cup \bigcup_{i=1}^{l-1}J_i,$$
 where $I_i\sqsubset P, J_i\sqsubset P_\sharp$
are distinct intervals and $l\geq 2$. (See Notations \ref{n1}.)

(3) For any $I\sqsubset P$ and $s\in S\setminus\{s_I\}$, 
$\varphi(s)(I)$  is a proper subset of some $I'\sqsubset P$.

(4) For any $J\sqsubset P_\sharp$ and $s\in S$, either $\varphi(s)J\sqsubset P_\sharp$ or $\varphi(s)J\subset {\rm Int}(P)$.
\end{definition}

\begin{example}\label{ex1}
The set $P_\nu$ ($\nu=1,2$) in Figure 6 is an
example of BP for homomorphisms $\varphi_\nu=\varphi\vert_{\Gamma_\nu}$. 
The set $P$ in Figure 7 is also a BP.
\end{example}

\begin{definition}\label{d2}
For a BP $P$ for $\varphi\in\RR_\Gamma$ and $l\geq2$, define inductively
$P^l=\bigcap_{s\in S\cup\{e\}}\varphi(s)P^{l-1}$, where $P^1=P$. Also define
$P^\infty=\bigcap_{l\in\N}P^l$.
\end{definition}

Thus $\{P^l\}_{l\in\N}$ is a decreasing sequence of compact subsets,
each consisting of finitely many closed intervals.
In Example \ref{ex1}, if the corresponding homomorphism is onto a
Shottky group, then $P^\infty$ coincides with the limit set.
In general, $P^\infty$ is a closed perfect set.

Let us see how $P^2$ is obtained
from $P$. By (2) and (3) of Definition \ref{d1}, we have
$$P^2=\bigcup_{I\sqsubset P}\varphi(s_I)^{-1}(P\cap s_I(I)).$$
That is, any interval $I\sqsubset P$ is divided uniquely
as:
$$
I=\bigcup_{i=1}^l\varphi(s_I)^{-1}(I_i) \cup
\bigcup_{i=1}^{l-1}\varphi(s_I)^{-1}(J_i),$$
where 
$\varphi(s_I^{-1})(I_i)\sqsubset P^2,\ \varphi(s_I^{-1})(J_i)\sqsubset P_\sharp^2=(P^2)_\sharp$.
Any $I'\sqsubset P^2$ is of the above form
$I'=\varphi(s_I)^{-1}(I_i)$, and
$\varphi(s_{I})$ maps $I'$ onto $I_i\sqsubset P$. For any other $s$,
 $\varphi(s)$ maps $I'$ onto a proper subset of
 some $I''\sqsubset P^2$. 
On the other hand, $P^2_\sharp$ is obtained from $P_\sharp$
by adding new intervals of the above form $\varphi(s_I)^{-1}(J_i)$. 
A component of $P^2_\sharp$ is called {\em level 1} if it is contained
in $P_\sharp$, and {\em level 2} otherwise.
Any level 1 component
is mapped by any $\varphi(s)$ onto
a component of $P^2_\sharp$, either to level 1  or to  level 2.
As for a level 2 component, we have the following.
\begin{enumerate}
\item A level 2 component 
$\varphi(s_I)^{-1}(J_i)$ is mapped by $\varphi(s_I)$ onto a level 1 component
$J_i$, and is mapped by
any other $\varphi(s)$ 
onto an interval contained in the interior of $P^2$. Especially no level
      2 component is mapped onto a level 2 component.
\end{enumerate}

By these considerations, we have the following lemma.

\begin{lemma}\label{l1}
For a BP $P$ for $\varphi\in\RR_\Gamma$ and $l\geq2$, $P^l$ is a BP for
 $\varphi$. \qed
\end{lemma}

Let $P$ (resp.\ $P'$) be a BP for $\varphi\in\RR_\Gamma$ (resp. $\varphi'$).
Recall that $P_*=P\cap P_\sharp$ from Notation \ref{n1}.

\begin{definition}
A COP (cyclic order preserving) bijection $\xi:P_*\to P'_*$
is called a {\em BP equivalence} if for any 
 $x,y\in P_*$ and $s\in S$, we have 

$\bullet$ $[x,y]\sqsubset P$ if and only if
 $[\xi(x),\xi(y)]\sqsubset P'$ and 

$\bullet$ $y=\varphi(s)x$ if and only if
$\xi(y)=\varphi'(s)\xi(x)$.
\end{definition}

\begin{lemma} \label{l5}
Let $P$ (resp.\ $P'$) be a BP for $\varphi\in\RR_\Gamma$ (resp.\ $\varphi'$). 
Then a BP equivalence $\xi:P_*\to P'_*$
extends uniquely to a BP equivalence $\xi^2:P^2_*\to P'^2_*$.
\end{lemma}

\bd  For any $x\in P^2_*\setminus P_*$, there exists a unique element
 $s\in S$ such that
$\varphi(s)x\in P_*$. Define
 $\xi^2(x)=\varphi'(s)^{-1}\circ \xi\circ\varphi(s)(x)$.
It is easy to show that $\xi^2$ is in fact a BP equivalence.
\qed

\medskip
Notice that $P^\infty=\bigcap_{l\in\N}P^l$ is a perfect closed set,
$P^\infty_\sharp=(P^\infty)_\sharp$ consists of countably many
disjoint closed intervals, and $P^\infty_*=P^\infty\cap P^\infty_\sharp$
is a countable set. All the three sets are $\varphi(\Gamma)$-invariant.

The next theorem says that if $P$ is a BP for
$\varphi\in\RR_\Gamma$, then the
semiconjugacy class of the homomorphism $\varphi$ is determined by
the simple dynamics of $S$ on $P$.
A semiconjugacy class is
in fact determined by how one or several orbits are located in $S^1$
(Proposition \ref{equivariant}).

\begin{theorem}\label{t1}
Let $P$ and $P'$ be BP's for $\varphi\in\RR_\Gamma$ and
$\varphi'\in\RR_\Gamma$. Then a BP equivalence $\xi:P_*\to P'_*$
extends uniquely to a $(\varphi,\varphi')$-equivariant COP bijection 
$\xi^\infty: P^\infty_*\to (P')^\infty_*$. 
\end{theorem}

\medskip
\bd This follows from inductive applications of Lemma \ref{l5}. \qed

\medskip

The next lemma plays a key role when we study a $k$-fold lift of a
homomorphism. The easy proof is omitted.

\begin{lemma}\label{l8}
Let $P$ be a BP for $\varphi\in\RR_\Gamma$ and $\psi$ a $k$-fold lift of
 $\varphi$. Then $\pi_k^{-1}(P)$ is a BP for $\psi$. \qed
\end{lemma}

The lemma joined with Theorem \ref{t1} says that if $\psi$ is a $k$-fold
lift of $\varphi$ which admits a BP $P$, then the semiconjugacy class
of $\psi$ is determined by the dynamics of $\psi(S)$ on $\pi_k^{-1}(P)$.

\bigskip
For future purpose, we need to continue to study
more  about BP's. Especially we have to show that
the stabilizer (defined later) of an interval $J\sqsubset P_\sharp$ can
be determined by a simple algorythm for a certain class of BP's.

\begin{definition}
For any $J\sqsubset P^\infty_\sharp$, define the {\em level} of
$J$, $\lev(J)\in \N$, by
$\lev(J)=l$ if and only if $J\subset P^l_\sharp\setminus P^{l-1}_\sharp$.
\end{definition}

\begin{lemma}\label{l2}
Let $P$ be a BP for $\varphi\in\RR_\Gamma$.
If $J\sqsubset P^\infty_\sharp$ satisfies $\lev(J)=l$ for some
$l\geq2$, then there is a unique element $s\in S$ such that
$\lev(\varphi(s)J)=l-1$, and for any other $s\in S$, $\lev(\varphi(s)J)=l+1$.
\end{lemma}

\bd For $l=2$, this follows from (1) placed just before Lemma \ref{l1}.
The general case can be shown by  an easy induction. \qed

\begin{definition}\label{right}
A labelled directed graph $G(P)$ associated with a BP $P$ for
 $\varphi\in\RR_\Gamma$ is defined as follows. The vertices of $G(P)$ are
components of $P_\sharp$. There is a directed edge from  $J_1$ to $J_2$
 with label $s\in S$ (written $J_1\stackrel{s}\rightarrow J_2$) if
$s=s_I$, where $I$ is the component of $P$ right adjacent to $J_1$,
and $\varphi(s)(J_1)=J_2$.
\end{definition}

\begin{example}\label{ex2}
The graph $G(P_1)$ and $G(P_2)$ of the BP's in Figure 6 consists of one
 cycle,
while the graph $G(P)$ for Figure 7 consists of 3 cycles. 
\end{example}

Notice that for any vertex $J$ of $G(P)$, there is exactly one edge leaving
$J$. However there may be a vertex at which no edges arrive.

\begin{definition}\label{l3}
A BP $P$ for $\varphi\in\RR_\Gamma$ is called {\em pure}, if the graph $G(P)$ consists of
disjoint cycles. We allow a period one cycle formed by one vertex
and one edge.
\end{definition}

In fact, the pureness does not change if we replace
``right adjacent'' by ``left adjacent'' in Definition \ref{right}, 
although the direction or labelling of the graph may change.
For any BP $P$, $P^2$ can never be pure. The BP's in Examples \ref{ex2} are pure.

\begin{definition}
For $\varphi\in\RR_\Gamma$ and a subset $A$ of $S^1$, the {\em
 stabilizer of $A$ with respect to $\varphi$}, denoted by $\Stab_\varphi(A)$,
is defined by
$$
\Stab_\varphi(A)=\{\gamma\in\Gamma\mid\varphi(\gamma)(A)=A\}.
$$
\end{definition}

\begin{lemma}\label{l4} Let $P$ be a pure BP for $\varphi\in\RR_\Gamma$. Then we
have the following.

(1) The group $\Gamma$ is free with symmetrized free generating set $S$
 and $\varphi$ is injective.

(2) For any $J\sqsubset P_\sharp$, the stabilizer $\Stab_\varphi(J)$ is
 generated by  an element written as a cyclically reduced word of $S$.

(3) For any $J\sqsubset P^\infty_\sharp$ with $\lev(J)=l$ ($l\geq2$),
$\Stab_\varphi(J)$ is generated by an element which has a nonreducing
 representation 
$\alpha\beta\alpha^{-1}$ by reduced words of $S$ such that the word length of $\alpha$ is
$l-1$ and $\beta$ is cyclically reduced.
\end{lemma}

\bd For any $J\sqsubset P_\sharp$, assume $\varphi(\gamma)(J)=J$ for some
$\gamma\in\Gamma\setminus\{e\}$. Write $\gamma$ as a reduced word in $S$:
$\gamma=s_m\cdots s_2s_1$. For any $1\leq i\leq m$, let
$J_i=\varphi(s_i\cdots s_1)(J)$.  
Then we have $\lev(J_i)=1$ for any $i$, that is, $J_i$ is a vertex
of $G(P)$. In fact, if $\lev(J_i)$ would take
the maximal value $l\geq2$ at some $i$,
then by Lemma \ref{l2}, we have
$\lev(J_{i-1})=\lev(J_{i+1})=l-1$ and $s_{i+1}=s_i^{-1}$, contrary
to the assumption that the word is reduced.
Again since the word is reduced and $P$ is pure,
we have either of the following.
$$J\stackrel{s_1}\to J_1\stackrel{s_2}\to\cdots\stackrel{s_m}\to J_m=J
\ \mbox{ or }\ J\stackrel{s_1}\leftarrow
J_1\stackrel{s_2}\leftarrow\cdots\stackrel{s_m}\leftarrow J_m=J.$$
Let $I_i\sqsubset P$ be an interval right adjacent to $J_i$.
Then in the former case $\varphi(s_{i+1})$ is always expanding on $I_i$,
that is,  $s_{i+1}=s_{I_i}$. This shows that $\varphi(\gamma)$ cannot be the
identity. The same is true in the latter case.
Points (1) and (2) follows from this, while it is easy to
derive (3) from (2). \qed

\medskip
Finally we shall prepare some terminologies and facts needed for the next section.
Let $\Lambda$ be an infinite cyclic subgroup of $\Gamma$ and
$\varphi\in\RR_\Gamma$.  

\begin{definition}
Given a closed subset $X$ of $S^1$, the set
$$
E_\varphi^\Lambda(X)=\bigcup\{J\sqsubset X_\sharp\mid \{e\}\neq\Stab_\varphi(J)\subset \Lambda\},$$
is called the {\em entrance} of $\Lambda$ to $X$
with respect to $\varphi$.
\end{definition}

\begin{definition}
A pair of closed subsets $(Q,E)$ is called a {\em
 $(\Gamma,\Lambda)$-pair for $\varphi$} 
if $Q$ is a $\varphi(\Gamma)$-invariant closed perfect
 set, $E=E_\varphi^\Lambda(Q)$, and $E$ is a finite 
disjoint union of closed intervals.
\end{definition}

\begin{lemma}\label{l9}
Let $P$ be a pure BP for $\varphi\in\RR_\Gamma$ and $\Lambda$ an infinite cyclic
 subgroup of $\Gamma$. Assume $E_\varphi^\Lambda(P)$ is nonempty. Then
$(P^\infty,E_\varphi^\Lambda(P))$ is a $(\Gamma,\Lambda)$-pair.
\end{lemma}

\bd We only need to show that
$E^\Lambda_\varphi(P)=E^\Lambda_\varphi(P^\infty)$. 
That is, if $J\sqsubset P^\infty_\sharp$ and $\Stab_\varphi(J)\subset
\Lambda$, then $\lev(J)=1$.
But this is clear from Lemma
\ref{l4}. \qed

\section{Combinations}

This section is divided into three subsections. In the first, we
are concerned with a single homomorphism, while in the second, with
a pair of homomorphisms.

\medskip
\noindent {\bf 1.} Throughout this subsection, we assume the following.

\begin{assumption}\label{a0}
(a) The group $G$ is written as an amalgamated product
$$G=\Gamma_1*_\Lambda\Gamma_2,$$
 where $\Lambda$ is an infinite cyclic subgroup.

(b) $\varphi\in\RR_G$, and $\varphi_\nu=\varphi\vert_{\Gamma_\nu}$ is
 injective for
$\nu=1,2$. 

(c) $(Q_\nu,E_\nu)$ is a $(\Gamma_\nu,\Lambda)$-pair for $\varphi_\nu$, $\nu=1,2$.
\end{assumption}

Denote $\Gamma_\nu^*=\Gamma_\nu\setminus\Lambda$. We make extensive use of the following
partition of the group $G$.

$$
G=\bigsqcup_{k\geq0}G^{k}, \mbox{ where }$$
\begin{equation}\label{e1}
G^0=\Lambda,\ G^{1}=\Gamma^*_1\sqcup\Gamma^*_2,\ G^2=\Gamma^*_1\Gamma^*_2\sqcup\Gamma^*_2\Gamma^*_1,\
G^3=\Gamma^*_1\Gamma^*_2\Gamma^*_1\sqcup\Gamma^*_2\Gamma^*_1\Gamma^*_2,\
\cdots.
\end{equation}

\begin{definition}
The pairs $(Q_1,E_1)$ and $(Q_2,E_2)$ are called {\em combinable for
 $\varphi$} 
if
$E_1$ and $E_2$ alternate in $S^1$, that is, $E_1\cup E_2=S^1$ and
$\Int(E_1)\cap\Int(E_2)=\emptyset$. In this case we denote
$E_*=\partial E_1=\partial E_2$.
\end{definition}

We also assume the following in this subsection.

\begin{assumption}
$\QQ=((Q_1,E_1),(Q_2,E_2))$ is a combinable pair for $\varphi$.
\end{assumption}

We define an (undirected) graph $(\VV(\QQ),\EE(\QQ))$ of the combinable
pair $\QQ$ as follows. 
$$
\VV(\QQ)=\{\varphi(g)Q_\nu\mid g\in G,\ \nu=1,2\},$$
$$
\EE(\QQ)=\{\{v,v'\}\mid v,v'\in\VV(\QQ),\ v\neq v',\  v\cap
v'\neq\emptyset\}.
$$
The group $G$ acts naturally on the graph $(\VV(\QQ),\EE(\QQ))$ as graph
automorphisms via the
homomorphism $\varphi$. 
The rest of this subsection is devoted to the study of properties 
of the graph $(\VV(\QQ),\EE(\QQ))$. Especially we show that the graph
$(\VV(\QQ),\EE(\QQ))$ is in fact a tree. (It is isomorphic to the
Bass-Serre tree associated to the amalgamated product
$G=\Gamma_1*_\Lambda\Gamma_2$.) 

For $v,w\in\VV(\QQ)$, we denote $v\sim w$ if $\{v,w\}\in\EE(\QQ)$, and
say
that $v$ and $w$ are adjacent.
The indexing set for $Q_\nu$ is the group $\Z/2\Z$, thus for example
$Q_3=Q_1$, while the indexing set for a group element is $\Z$, thus
in general $\gamma_3\neq\gamma_1$.

\begin{lemma}\label{l10} We have $Q_\nu\sim Q_{\nu+1}$ and
$Q_\nu\sim \varphi(\gamma_\nu)Q_{\nu+1}$
for any
 $\gamma_\nu\in\Gamma_\nu^*$.
Conversely if $ Q_\nu\sim v$, then either
 $v=Q_{\nu+1}$ or $v=\varphi(\gamma_\nu)Q_{\nu+1}$ for some
 $\gamma_\nu\in\Gamma_\nu^*$. 
\end{lemma}

\bd 
Since $Q_1\cap Q_2=E_*\neq\emptyset$, we have $Q_1\sim Q_2$.
Since $Q_\nu$ is invariant by $\varphi(\Gamma_\nu)$. we have
$Q_\nu\cap\varphi(\gamma_{\nu})Q_{\nu+1}=\varphi(\gamma_{\nu})(Q_1\cap Q_2)\neq\emptyset$ for $\gamma_{\nu}\in\Gamma_{\nu}^*$. 
That is, $Q_\nu\sim\varphi(\gamma_\nu)Q_{\nu+1}$.

In the rest, we shall show that all the other vertices are not adjacent
to $Q_\nu$. First we prepare some fundamental facts. See Figure 8.
\begin{figure}\caption{\small The subsets $Q_\nu$ should have countably many complementary intervals.
Only some of them are drawn in the figure.}
\unitlength 0.1in
\begin{picture}( 45.7000, 13.8000)( 27.4000,-45.7000)
%
{\color[named]{Black}{%
\special{pn 8}%
\special{pa 2740 3820}%
\special{pa 7310 3820}%
\special{fp}%
}}%
%
{\color[named]{Black}{%
\special{pn 8}%
\special{pa 3320 3950}%
\special{pa 3330 3950}%
\special{fp}%
}}%
%
{\color[named]{Black}{%
\special{pn 8}%
\special{pa 3200 3710}%
\special{pa 3190 3940}%
\special{fp}%
}}%
%
{\color[named]{Black}{%
\special{pn 8}%
\special{pa 5040 3710}%
\special{pa 5040 3940}%
\special{fp}%
}}%
%
{\color[named]{Black}{%
\special{pn 8}%
\special{pa 6750 3720}%
\special{pa 6750 3940}%
\special{fp}%
}}%
%
{\color[named]{Black}{%
\special{pn 20}%
\special{pa 3180 3810}%
\special{pa 3420 3810}%
\special{fp}%
}}%
%
{\color[named]{Black}{%
\special{pn 20}%
\special{pa 3780 3820}%
\special{pa 3970 3820}%
\special{fp}%
}}%
%
{\color[named]{Black}{%
\special{pn 20}%
\special{pa 4320 3820}%
\special{pa 4520 3820}%
\special{fp}%
}}%
%
{\color[named]{Black}{%
\special{pn 20}%
\special{pa 4830 3820}%
\special{pa 5210 3820}%
\special{fp}%
}}%
%
{\color[named]{Black}{%
\special{pn 20}%
\special{pa 5570 3820}%
\special{pa 5820 3820}%
\special{fp}%
}}%
%
{\color[named]{Black}{%
\special{pn 20}%
\special{pa 6140 3820}%
\special{pa 6320 3820}%
\special{fp}%
}}%
%
{\color[named]{Black}{%
\special{pn 20}%
\special{pa 6590 3820}%
\special{pa 6740 3820}%
\special{fp}%
}}%
\put(38.3000,-34.4000){\makebox(0,0)[lb]{$Q_{\nu+1}$}}%
\put(56.3000,-34.2000){\makebox(0,0)[lb]{$Q_\nu$}}%
%
{\color[named]{Black}{%
\special{pn 4}%
\special{pa 3910 3530}%
\special{pa 3360 3740}%
\special{fp}%
\special{sh 1}%
\special{pa 3360 3740}%
\special{pa 3430 3736}%
\special{pa 3410 3722}%
\special{pa 3416 3698}%
\special{pa 3360 3740}%
\special{fp}%
}}%
%
{\color[named]{Black}{%
\special{pn 4}%
\special{pa 4010 3540}%
\special{pa 3910 3770}%
\special{fp}%
\special{sh 1}%
\special{pa 3910 3770}%
\special{pa 3956 3718}%
\special{pa 3932 3722}%
\special{pa 3918 3702}%
\special{pa 3910 3770}%
\special{fp}%
}}%
%
{\color[named]{Black}{%
\special{pn 4}%
\special{pa 4240 3490}%
\special{pa 4390 3750}%
\special{fp}%
\special{sh 1}%
\special{pa 4390 3750}%
\special{pa 4374 3682}%
\special{pa 4364 3704}%
\special{pa 4340 3702}%
\special{pa 4390 3750}%
\special{fp}%
}}%
%
{\color[named]{Black}{%
\special{pn 4}%
\special{pa 4470 3510}%
\special{pa 4900 3750}%
\special{fp}%
\special{sh 1}%
\special{pa 4900 3750}%
\special{pa 4852 3700}%
\special{pa 4854 3724}%
\special{pa 4832 3736}%
\special{pa 4900 3750}%
\special{fp}%
}}%
%
{\color[named]{Black}{%
\special{pn 4}%
\special{pa 5620 3500}%
\special{pa 5200 3740}%
\special{fp}%
\special{sh 1}%
\special{pa 5200 3740}%
\special{pa 5268 3724}%
\special{pa 5246 3714}%
\special{pa 5248 3690}%
\special{pa 5200 3740}%
\special{fp}%
}}%
%
{\color[named]{Black}{%
\special{pn 4}%
\special{pa 5750 3530}%
\special{pa 5700 3750}%
\special{fp}%
\special{sh 1}%
\special{pa 5700 3750}%
\special{pa 5734 3690}%
\special{pa 5712 3698}%
\special{pa 5696 3682}%
\special{pa 5700 3750}%
\special{fp}%
}}%
%
{\color[named]{Black}{%
\special{pn 4}%
\special{pa 5950 3490}%
\special{pa 6220 3750}%
\special{fp}%
\special{sh 1}%
\special{pa 6220 3750}%
\special{pa 6186 3690}%
\special{pa 6182 3714}%
\special{pa 6158 3718}%
\special{pa 6220 3750}%
\special{fp}%
}}%
%
{\color[named]{Black}{%
\special{pn 4}%
\special{pa 6080 3500}%
\special{pa 6650 3750}%
\special{fp}%
\special{sh 1}%
\special{pa 6650 3750}%
\special{pa 6598 3706}%
\special{pa 6602 3730}%
\special{pa 6582 3742}%
\special{pa 6650 3750}%
\special{fp}%
}}%
\put(39.3000,-44.1000){\makebox(0,0)[lb]{$E_\nu$}}%
\put(55.7000,-44.1000){\makebox(0,0)[lb]{$E_{\nu+1}$}}%
%
{\color[named]{Black}{%
\special{pn 4}%
\special{pa 3810 4350}%
\special{pa 3778 4344}%
\special{pa 3744 4336}%
\special{pa 3648 4312}%
\special{pa 3584 4292}%
\special{pa 3554 4282}%
\special{pa 3524 4270}%
\special{pa 3496 4258}%
\special{pa 3440 4230}%
\special{pa 3388 4194}%
\special{pa 3364 4176}%
\special{pa 3320 4132}%
\special{pa 3300 4108}%
\special{pa 3240 4030}%
\special{pa 3222 4002}%
\special{pa 3204 3976}%
\special{pa 3200 3970}%
\special{fp}%
}}%
%
{\color[named]{Black}{%
\special{pn 4}%
\special{pa 4390 4350}%
\special{pa 4426 4348}%
\special{pa 4562 4332}%
\special{pa 4658 4314}%
\special{pa 4690 4306}%
\special{pa 4750 4286}%
\special{pa 4778 4272}%
\special{pa 4804 4258}%
\special{pa 4830 4242}%
\special{pa 4854 4226}%
\special{pa 4878 4206}%
\special{pa 4900 4184}%
\special{pa 4920 4160}%
\special{pa 4940 4134}%
\special{pa 4976 4082}%
\special{pa 4994 4054}%
\special{pa 5010 4024}%
\special{pa 5028 3996}%
\special{pa 5030 3990}%
\special{fp}%
}}%
%
{\color[named]{Black}{%
\special{pn 4}%
\special{pa 5480 4360}%
\special{pa 5454 4342}%
\special{pa 5426 4324}%
\special{pa 5348 4270}%
\special{pa 5322 4250}%
\special{pa 5298 4230}%
\special{pa 5272 4210}%
\special{pa 5250 4188}%
\special{pa 5226 4166}%
\special{pa 5204 4144}%
\special{pa 5182 4120}%
\special{pa 5162 4096}%
\special{pa 5142 4070}%
\special{pa 5122 4046}%
\special{pa 5104 4020}%
\special{pa 5084 3994}%
\special{pa 5066 3968}%
\special{pa 5060 3960}%
\special{fp}%
}}%
%
{\color[named]{Black}{%
\special{pn 4}%
\special{pa 6290 4360}%
\special{pa 6318 4358}%
\special{pa 6344 4356}%
\special{pa 6370 4352}%
\special{pa 6398 4346}%
\special{pa 6426 4338}%
\special{pa 6452 4326}%
\special{pa 6482 4312}%
\special{pa 6510 4296}%
\special{pa 6540 4272}%
\special{pa 6570 4246}%
\special{pa 6602 4212}%
\special{pa 6634 4174}%
\special{pa 6666 4130}%
\special{pa 6696 4086}%
\special{pa 6722 4044}%
\special{pa 6742 4008}%
\special{pa 6754 3984}%
\special{pa 6754 3972}%
\special{pa 6744 3978}%
\special{pa 6730 3990}%
\special{fp}%
}}%
%
{\color[named]{Black}{%
\special{pn 8}%
\special{pa 3180 4440}%
\special{pa 3580 3890}%
\special{fp}%
\special{sh 1}%
\special{pa 3580 3890}%
\special{pa 3526 3932}%
\special{pa 3550 3934}%
\special{pa 3558 3956}%
\special{pa 3580 3890}%
\special{fp}%
}}%
%
{\color[named]{Black}{%
\special{pn 8}%
\special{pa 5210 4410}%
\special{pa 5410 3890}%
\special{fp}%
\special{sh 1}%
\special{pa 5410 3890}%
\special{pa 5368 3946}%
\special{pa 5392 3940}%
\special{pa 5406 3960}%
\special{pa 5410 3890}%
\special{fp}%
}}%
\put(29.3000,-46.0000){\makebox(0,0)[lb]{$\varphi(\gamma_{\nu+1})(E_{\nu+1})$}}%
\put(51.3000,-45.9000){\makebox(0,0)[lb]{$\varphi(\gamma_\nu)(E_\nu)$}}%
\end{picture}%

\end{figure}
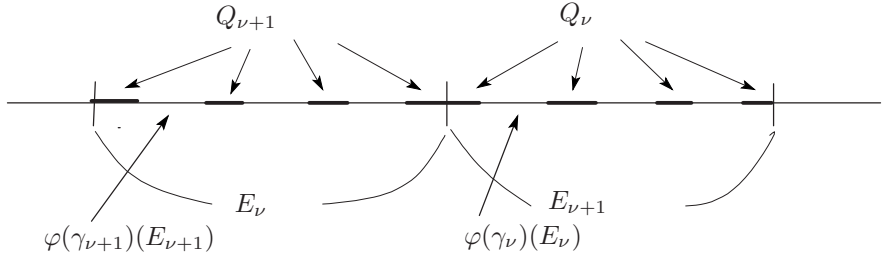
\begin{eqnarray}\label{e2}
& Q_\nu\subset E_{\nu+1} \mbox{ and }Q_\nu\cap\Int E_\nu=\emptyset,\\
& \varphi(\gamma_\nu)E_\nu\subset\Int E_{\nu+1}
\mbox{ and } \Int\varphi(\gamma_\nu)E_\nu\cap Q_\nu=\emptyset \mbox{ for any }\gamma_\nu\in\Gamma_\nu^*.\label{e2.5}
\end{eqnarray}
For (\ref{e2.5}), recall that $Q_\nu$ is assumed to be perfect.

Now for $\gamma_{\nu+1}\in\Gamma_{\nu+1}^*$,
the vertex $\varphi(\gamma_{\nu+1})Q_\nu$ is not adjacent to $Q_\nu$, since
$$\varphi(\gamma_{\nu+1})Q_\nu\subset\varphi(\gamma_{\nu+1})E_{\nu+1}\subset
\Int E_\nu.
$$

We shall show by induction on $k$ that if $k\geq2$ and
$\gamma_{\nu+i}\in\Gamma^*_{\nu+i}$ ($1\leq i\leq k$), then
\begin{equation}\label{e3}
\varphi(\gamma_{\nu+k}\cdots\gamma_{\nu+1})Q_\nu\subset\Int \varphi(\gamma_{\nu+k})E_{\nu+k}.
\end{equation}
This shows that $\varphi(\gamma_{\nu+k}\cdots\gamma_{\nu+1})Q_\nu$ is
adjacent neither to $Q_1$ nor to $Q_2$, by virtue of
(\ref{e2.5}).

To show (\ref{e3}) for $k=2$, notice by (\ref{e2})
$$
\varphi(\gamma_{\nu+2}\gamma_{\nu+1})Q_\nu\subset\varphi(\gamma_{\nu+2}\gamma_{\nu+1})E_{\nu+1}\subset\Int\varphi(\gamma_{\nu+2})
E_{\nu+2}.$$
For the inductive step, 
$$\varphi(\gamma_{\nu+k+1}\gamma_{\nu+k}\cdots\gamma_{\nu+1})Q_\nu\subset\varphi(\gamma_{\nu+k+1})
\Int\varphi(\gamma_{\nu+k})E_{\nu+k}
\subset\Int\varphi(\gamma_{\nu+k+1})E_{\nu+k+1}.$$
\qed

\begin{remark}\label{r}
The above proof shows that any vertex $\varphi(g)Q_\nu$ is distinct from
$Q_\nu$ unless $g\in\Gamma_\nu$.
\end{remark}

\begin{figure}
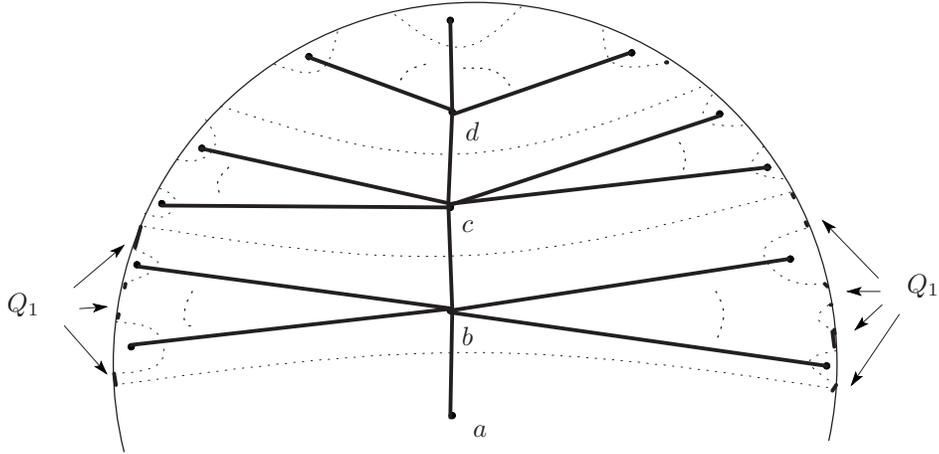
\caption{\small The vertices $a$--$d$ are $a=Q_2$, $b=Q_1$,
 $c=\gamma_1Q_2$ and $d=\gamma_1\gamma_2Q_1$. The actual set $Q_1$ is
 depicted on the circle.}
\unitlength 0.1in
%

\end{figure}

See Figure 9 for the graph $(\VV(\QQ),\EE(\QQ))$.

\begin{lemma}\label{l11}
For any interval $J\sqsubset (Q_\nu)_\sharp$ which is 
distinct from $\varphi(\gamma_\nu)E_\nu$ for any $\gamma_\nu\in \Gamma_\nu$,
we have $\Int J\cap v=\emptyset$ for any $v\in\VV(\QQ)$
and  $\Stab_\varphi(J)=\Stab_{\varphi_\nu}(J)$.  
\end{lemma}

\bd Any vertex other than $Q_\nu$ contained in $E_{\nu+1}$ is contained in
$\varphi(\gamma_\nu)E_\nu$ for some $\gamma_\nu\in \Gamma_\nu^*$, by virtue
of (\ref{e3}), showing the first statement.  
For the last statement, choose an arbitrary element $g\in\Stab_\varphi(J)$.
Then $g$ leaves $\partial J$ invariant.
The set $\partial J$ is contained in $Q_\nu$ and
disjoint from any other vertex of $\VV(\QQ)$.
Therefore $g$ stabilizes the vertex $Q_\nu$ in the action of $G$ on the 
graph. This shows $g\in\Gamma_\nu$ by Remark \ref{r}. \qed

\medskip

\medskip
Let us continue the study of the graph $(\VV(\QQ),\EE(\QQ))$.

\begin{lemma}\label{l12}
Let $v,w\in\VV(\QQ)$. If $v\sim w$ and
 $v=\varphi(\gamma_{\nu+k}\cdots\gamma_{\nu+1})Q_\nu$ for some $k\geq1$ and
 $\gamma_{\nu+i}\in\Gamma_{\nu+i}^*$, then either
 $w=\varphi(\gamma_{\nu+k}\cdots\gamma_{\nu+2})Q_{\nu+1}$ or
$w=\varphi(\gamma_{\nu+k}\cdots\gamma_{\nu+1}\gamma_{\nu})Q_{\nu-1}$
for some $\gamma_\nu\in\Gamma_\nu^*$, and moreover 
$v\cap w$ is contained in the $\varphi(G)$-orbit of $E_*$.
\end{lemma}

\bd Recall that the group $G$ acts on the graph $(\VV(\QQ),\EE(\QQ))$
as graph automorphisms. Thus if
$w\sim\varphi(\gamma_{\nu+k}\cdots\gamma_{\nu+1})Q_\nu$, then
$\varphi(\gamma_{\nu+k}\cdots\gamma_{\nu+1})^{-1}w\sim Q_\nu$.
Therefore either
$\varphi(\gamma_{\nu+k}\cdots\gamma_{\nu+1})^{-1}w$ is equal to $Q_{\nu+1}$
or $\varphi(\gamma_\nu)Q_{\nu-1}$ for some $\gamma_\nu\in\Gamma^*_\nu$.
Since $\varphi(\gamma_{\nu+1})Q_{\nu+1}=Q_{\nu+1}$, this shows the first part. 
An immediate consequence is that $G$ acts transitively on the set of
edges $\EE(\QQ)$.
That is, there is $g\in G$ which maps
 $E_*=Q_1\cap Q_2$ onto $v\cap w$, showing the second part. \qed

\begin{definition}\label{distance}
Any vertex $v$ of the graph is written as
$v=\varphi(\gamma_{\nu+k}\cdots\gamma_{\nu+1})Q_\nu$ for
$\gamma_{\nu+i}\in\Gamma_{\nu+i}^*$. The number $k$ is unique,
and is called the {\em distance}
of $v$.
\end{definition}

\begin{lemma}\label{l13}
(1) We have $\Stab_\varphi(Q_\nu)=\Gamma_\nu$.

(2) The graph $(\VV(\QQ),\EE(\QQ))$ is a tree and $\varphi$ is injective.

(3) For $g\in G$, we have $\varphi(g)E_*\cap E_*\neq\emptyset$ if and only if $g\in\Lambda$.
\end{lemma}

\bd
Point (1) is a rephrasing of Remark  \ref{r}.
Lemma \ref{l12} says that any vertex of distance $k$ ($k\geq1$)
is only adjacent to vertices of distance $k-1$ and $k+1$, and the
 vertex of distance $k-1$ is unique. This shows that
the graph $(\VV(\QQ),\EE(\QQ))$ is a tree.
To show that $\varphi$ is injective,
assume $\varphi(g)={\rm id}$ for some $g\in G$. Then $g$ acts trivially
on the graph. Especially $g$ leaves the vertex $Q_\nu$ invariant. That is, 
$g\in\Gamma_\nu$ by (1). By the assumption that
$\varphi_\nu=\varphi\vert_{\Gamma_\nu}$ is injective, we get $g=e$, as is
required. 
The if part of (3) is clear. To show the converse, notice that 
$E_*= Q_1\cap Q_2$. By Lemma \ref{l10}, if
$\varphi(g)Q_\nu\cap Q_{\nu+1}\neq\emptyset$, then $g\in\Gamma_{\nu+1}$.
This holds for each $\nu=1,2$, and thus $g\in\Gamma_1\cap\Gamma_2=\Lambda$. \qed

\medskip
Further discussions are necessary for the development of the
next subsection.

\begin{definition}\label{J}
Let $\JJ_0$ be the family of connected components of $E_1$ and 
 $E_2$, and  for $n\geq1$, let
 $$\JJ_n=\{\varphi(\gamma_{\nu+n-1}\cdots\gamma_{\nu})I\mid \gamma_{\nu+i}\in\Gamma^*_{\nu+i},\ I\sqsubset E_\nu,\,\nu=1,2\}.$$
\end{definition}

\begin{lemma}\label{l14}

(1) For any $J\in\JJ_n$ ($n\geq2$), there are $J_{i}\in\JJ_{i}$
 ($1\leq i\leq n-1$) and $\nu$ such
 that
$$
J\subset  \Int J_{n-1}\subset J_{n-1}\subset\cdots\subset \Int J_1\subset
J_1\sqsubset (Q_\nu)_\sharp,$$

(2) Any two intervals $J,J'\in\JJ_n$, $n\geq1$, satisfy either $J=J'$ or
$J\cap J'=\emptyset$.
\end{lemma}

\bd 
Point (1) is shown inductively using (\ref{e2.5}). 
Point (2) is clear for $n=1$. (See Figure 8.)
The general case can be shown by an induction on $n$ based upon (1). \qed

\medskip
For a subset $K$ of $G$ and $X$ of $S^1$, denote $\varphi(K)X=\cup_{g\in
K}\varphi(g)X$. 
\begin{definition}\label{X}
Define a subset $X_n$ of $S^1$ by $X_0=E_*$ and for $n\geq1$,
 $X_{n}=\varphi(G^n)X_0$. Let $X=\bigcup_nX_n$. (For the definition of
 $G^n$, see (\ref{e1}).)
\end{definition}

The following easy lemma is useful to clarify an argument in the next
subsection. 

\begin{lemma}\label{added}

(1) $X_n=\bigcup_{I\in\JJ_n}\partial I$.

(2) $X_n\cap X_m=\emptyset$ if $n\neq m$.

(3) $X=\varphi(G)E_*$.

(4) For any $v\in\VV(\QQ)$ of distance $n$ ($n\geq0$), $v\cap X_m\neq\emptyset$ if
 and only if $m=n$ or $m=n+1$. \qed

\end{lemma}

\medskip
The following lemma will be used in Section 6 where we consider successive combinations.

\begin{lemma}\label{next section}
Assume there is a subset $E'\subset (Q_1)_\sharp$ such that $(Q_1,E')$ is a
 $(\Gamma_1,\Lambda')$-pair for $\varphi_1$, where $\Lambda'$ is an 
infinite cyclic subgroup of
$\Gamma_1$ such that $\Lambda'\cap\gamma_1\Lambda\gamma_1^{-1}=\{e\}$ for
 any $\gamma_1\in\Gamma_1$. 
Then $({\rm Cl}(\varphi(G)(Q_1\cup Q_2)),E')$ is
a $(G,\Lambda')$-pair for $\varphi$.
\end{lemma}

\bd It is clear that
$$Z:={\rm Cl}(\varphi(G)(Q_1\cup Q_2))={\rm Cl}(\bigcup_{v\in\VV(\QQ)}v)
$$
is a $\varphi(G)$-invariant closed perfect  set. 
So what is left is to show that $E'= E^{\Lambda'}_\varphi(Z)$,
where by definition $E'=E^{\Lambda'}_{\varphi_1}(Q_1)$.
The assumption on $\Lambda'$ implies that 
$E'\cap \varphi(\gamma_\nu)E_\nu=\emptyset$ for any $\nu$ and
$\gamma_\nu\in\Gamma_\nu^*$.
By Lemma \ref{l11}, we
have $E'\subset E^{\Lambda'}_\varphi(Z)$. 
To show the converse, assume $J\sqsubset E_\varphi^{\Lambda'}(Z)$. If 
$J\sqsubset (Q_1\cup Q_2)_\sharp$, then clearly we have
$J\sqsubset E'$. Otherwise $J$ must be contained in
 $\varphi(\gamma_\nu)I_\nu\in\JJ_1$ for some $I_\nu\sqsubset E_\nu$ and $\gamma_\nu\in\Gamma_\nu^*$. Since $J\sqsubset E_\varphi^{\Lambda'}(Z)$, there is
$g\in \Lambda'\setminus\{e\}\subset\Gamma_1^*$ such that $\varphi(g)J=J$. 
Then $\varphi(g)\varphi(\gamma_\nu)I_\nu\cap\varphi(\gamma_\nu)I_\nu\neq\emptyset$.
If $\nu=2$, then $\varphi(\gamma_2)I_2\subset \Int E_1$,
while $\varphi(g)\varphi(\gamma_2)I_2\subset \Int E_2$. A contradiction.
If $\nu=1$,
 $\varphi(g)\varphi(\gamma_1)I_1\in \JJ_1$ since
$g\in\Gamma_1^*$. Then by Lemma \ref{l14} (2), 
 $\varphi(g\gamma_1)I_1=\varphi(\gamma_1)I_1$, and
$\gamma_1^{-1} g\gamma_1\in\Lambda$ by Lemma \ref{l13} (3).
But this is contrary to the assumption on $\Lambda'$.
\qed

\bigskip
\noindent
{\bf 2.} In this subsection, we assume the following.

\begin{assumption}\label{a1}

 Let $\nu=1,2$ and $i=1,2$.

(a) The group $G$ is just as in Assumption \ref{a0}.

(b) Let $\varphi^i\in\RR_G$, and assume
 $\varphi^i_\nu=\varphi^i\vert_{\Gamma_\nu}$ is injective.

(c) Let $(Q^i_\nu,E^i_\nu)$  be a
$(\Gamma_\nu,\Lambda)$-pair for $\varphi^i_\nu$.

(d) The pair $\QQ^i=((Q^i_1,E^i_1),(Q^i_2,E^i_2))$ 
 is combinable for $\varphi^i$.

(e) There is a COP bijection $\xi:Q^1_{1,*}\cup Q^1_{2,*}\to Q^2_{1,*}\cup Q^2_{2,*}$ such that $\xi(Q^1_{\nu,*})=Q^2_{\nu,*}$
and the restrictions $\xi_\nu=\xi\vert_{Q^1_{\nu,*}}:Q^1_{\nu,*}\to Q^2_{\nu,*}$
is $(\varphi_\nu^1,\varphi_\nu^2)$-equivariant.
\end{assumption}

\medskip
Our purpose is to show that $\xi$ extends to a $(\varphi^1,\varphi^2)$-equivariant COP
bijection from the saturation $\varphi^1(G)(Q^1_{1,*}\cup Q^1_{2,*})$ to
$\varphi^2(G)(Q^2_{1,*}\cup Q^2_{2,*})$ (Theorem \ref{t2}).
The proof is by two steps: the first step is the following Lemma.
Let $\JJ^i_n$, $X^i_n$ and $X^i$ be defined as in Definitions
\ref{J} and  \ref{X} for $\varphi^i$.

\begin{lemma}\label{l16}
The map $\xi$ extends to a  COP bijection
$$\hat\xi:Q^1_{1,*}\cup Q^1_{2,*}\cup X^1\to Q^2_{1,*}\cup Q^2_{2,*}\cup X^2$$
which is $(\varphi^1,\varphi^2)$-equivariant as a map from $X^1$ to $X^2$.
\end{lemma}

\bd Recall that $X^i=\varphi_i(G)E^i_*$. The map $\xi$ extends to
$X^1$ by the $(\varphi^1,\varphi^2)$-equivariance. Namely, given $x\in X^1$, choose $g\in G$ and
$x_0\in E^1_*$ such that $x=\varphi^1(g)x_0$, and define
$\hat\xi(x)=\varphi^2(g)\xi(x_0)$. The map $\hat\xi$ is a well defined bijection since by Lemma
\ref{l13} (3), $\Stab_{\varphi^i}(E^i_*)\subset\Lambda$, and
$\xi\vert_{E^1_*}$ is
$(\varphi^1\vert_{\Lambda},\varphi^2\vert_{\Lambda})$-equivariant. Notice also
that $\hat\xi$ coincides with the original $\xi$ on $X^1_1\subset
Q_{1,*}^1\cup Q_{2.*}^1$ by the $(\varphi^1_\nu,\varphi^2_\nu)$-equivariance
of $\xi_\nu$, and $X^1_n$ ($n\geq2$) is disjoint from
$Q_{1,*}^1\cup Q_{2.*}^1$.
Therefore we only need to show that $\hat\xi$ is COP.

We shall prove that $\hat\xi$ is COP on 
$Q^1_{1,*}\cup Q^1_{2,*}\cup \bigcup_{0\leq i\leq n}X^1_i$ by an
induction on $n$. This is sufficient since $X^1=\bigcup_nX^1_n$.
For $n=1$, this is true by the assumption since 
$X^1_1\subset
Q^1_{1,*}\cup Q^1_{2,*}$. To show it for $n+1$,
choose an arbitrary open interval 
$$\Int I\sqsubset S^1\setminus(Q^1_{1,*}\cup Q^1_{2,*}\cup \bigcup_{0\leq i\leq n}X^1_i)$$ such that $\Int I\cap X^1_{n+1}\neq\emptyset.$
Clearly we only have to show that $\hat\xi$ is COP on 
$$I\cap(Q^1_{1,*}\cup Q^1_{2,*}\cup \bigcup_{0\leq i\leq n+1}X^1_i),$$ where  $I$ is the closure 
of $\Int I$. 
Now any point of $\Int I\cap X^1_{n+1}$ is an endpoint of some interval
of $\JJ_{n+1}^1$, and by Lemma \ref{l14},
we have $I\in\JJ^1_n$. This shows $$I\cap(Q^1_{1,*}\cup Q^1_{2,*}\cup
\bigcup_{0\leq i\leq n+1}X^1_i)=I\cap(X_n^1\cup X_{n+1}^1).$$
Furthermore $I=\varphi^1(g)J$ for some $J\in\JJ^1_0$ and
$g\in G^n$.

Finally since we have defined
$$\hat\xi\vert_{I\cap (X^1_n\cup X^1_{n+1})}=(\varphi^2(g)\vert_{\xi(J)\cap
(X^2_0\cup X^2_1)})\circ(\xi\vert_{J\cap (X^1_0\cup X^1_1})
\circ(\varphi^1(g^{-1})\vert_{I\cap (X_n^1\cup X^1_{n+1})}),$$
 and all the maps on the RHS is COP,
the map $\hat\xi\vert_{I\cap (X_n^1\cup X^1_{n+1})}$ is COP, as is required. \qed

\begin{theorem}\label{t2}
Under Assumption \ref{a1}, the COP bijection $$\xi:Q^1_{1,*}\cup Q^1_{2,*}\to Q^2_{1,*}\cup Q^2_{2,*}$$ extends uniquely to a $(\varphi^1,\varphi^2)$-equivariant
COP bijection $$\hat\xi:\varphi^1(G)(Q^1_{1,*}\cup Q^1_{2,*})\to\varphi^2(G)(Q^2_{1,*}\cup Q^2_{2,*}).$$
\end{theorem}

\bd Recall that 
$$
\varphi^1(G)(Q^1_1\cup Q^1_2)=\bigcup\{v\mid v\in\VV(\QQ^1)\},$$
 where $v=\varphi^1(g)Q^1_\nu$ for some $g\in G$ and $\nu$.
Denote $v_*=\varphi^1(g)Q^1_{\nu,*}$. Define $\hat\xi$ on each $v_*$ by the
$(\varphi^1,\varphi^2)$-equivariance. This is well defined because $\xi$ is
$(\varphi^1_\nu,\varphi^2_\nu)$-equivariant on $Q^1_{\nu,*}$ and
$\Stab_{\varphi^i}(Q^i_\nu)=\Gamma_\nu$ by Lemma \ref{l13} (1).
Of course the map $\hat\xi$ is COP on each $v_*$.
The map $\hat\xi$ coincides with the one defined in Lemma \ref{l16}
on $v_*\cap X^1$. The proof is complete by Lemma \ref{l16}.
\qed

\bigskip
\noindent
{\bf 3.} Let $\nu=1,2$ and $i=1,2$. Assume the following.

(a) The group $G$ is just as in Assumption \ref{a0}.

(b) Let $\varphi^i\in\RR_G$, and denote $\varphi^i_\nu=\varphi^i\vert_{\Gamma_\nu}$. 

(c) Let $P_\nu^i$ is a pure BP for $\varphi_\nu^i$, with $E_\nu^i$ the
entrance of $\Lambda$ to $P_\nu^i$.

(d) The pairs $(P_1^i,E_1^i)$ and $(P_2^i,E_2^i)$ are combinable in the
sense that $E_1^i$ and $E_2^i$ are alternating in $S^1$.

(e) There is a COP bijection 
$\xi: P^1_{1,*}\cup P^1_{2,*}\to P^2_{1,*}\cup P^2_{2,*}$ such that
$\xi\vert_{P^1_{\nu,*}}$ is a BP equivalence from $P^1_{\nu,*}$ 
onto $P^2_{\nu,*}$.

Joining Theorems \ref{t1} and \ref{t2}, we get the following.

\begin{theorem} \label{t2.5} 
Under the above assumption, $\varphi^1$ and $\varphi^2$ are semiconjugate.
\qed
\end{theorem}

\medskip
Notice that the set $R$ of fourteen points in Figure 5 is equal to
$P^1_{1,*}\cup P^1_{2,*}$ for the homomorphism (here denoted
$\varphi^1$) in $\RR_{\Pi_2}$ with $eu(\varphi^1)=2$.
Thus the above theorem says that any homomorphisms which admit
the same configuration $R$ are mutually semiconjugate.
This, together with the robustness of $R$ (discussed in Section 7),
implies the local stability of $\varphi^1$.
Furthermore a $2$-fold lift of $\varphi^1$ is also shown to be locally stable.

\section{Trees of groups}

\begin{definition} {\em A tree of groups} is a finite tree $\TT=(\VVV,\EEE)$
 such that
\\
(1) to each vertex $v\in\VVV$ (resp.\ edge $e\in\EEE$) is associated a
 group $\Gamma_v$ (resp.\ $\Lambda_e$), 
\\
(2) and if $v\in V$ is an end point of $e\in \EEE$, then a monomorphism
 $\iota_{e,v}:\Lambda_e\to\Gamma_v$ is assigned. 

The {\em fundamental group} $G(\TT)$ of a tree $\TT$ of groups is the group
 generated by $\Gamma_v$ and $\Lambda_e$ ($v\in\VVV,e\in\EEE$) subject
 to the relation $\lambda=\gamma$ whenever $\lambda\in \Lambda_e$,
$\gamma\in\Gamma_v$, $v$ is an end point of $e$, and
 $\iota_{e,v}(\lambda)=\gamma$.
\end{definition}

\begin{example}\label{ex4}
Consider the closed oriented surface $\Sigma_g$ of genus $g$.
Divide $\Sigma_g$ by circles into once puctured tori and pairs of
pants. Embed a tree in $\Sigma_g$ as in Figure 10 top. Then the fundamental groups
 $\Gamma_i$ of subsurfaces $\Sigma_i$ and the fundamental groups $\Lambda_j$
of circles $C_j$ are considered to be subgroups of the fundamental
 group $\Pi_g$ of the total surface, the base points being taken on the tree.
This yields a tree of groups as in Figure 10 bottom whose fundamental
 group is isomorphic to $\Pi_g$.
\end{example}

\begin{figure}\caption{}
\unitlength 0.1in
%

\end{figure}

Throughout this section we work under the following assumption.

\begin{assumption}\label{a3}
\noindent
(a) The group $G=G(\TT)$ is the fundamental group of a tree
 $\TT=(\VVV,\EEE)$ of groups. 
\\
(b) The vertex group $\Gamma_v$ admits a finite symmetric generating set
 $S_v$, and the edge group $\Lambda_e$ is infinite cyclic.
\\
(c) If $e$ and $e'$ are distinct edges starting at a vertex $v$,
then $\Lambda_e\cap\lambda_v\Lambda_{e'}\lambda_v^{-1}=\{e\}$ for any
$\lambda_v\in\Gamma_v$.
\\
(d) There are two homomorphisms $\varphi^i\in\RR_G$, $i=1,2$. We denote
$\varphi^i_v$ the restriction of $\varphi^i$ to the vertex group $\Gamma_v$.
\\
(e) For each vertex group $\Gamma_v$, there is a pure BP $P^i_v$ for $\varphi^i_v$
with repect to the generating set $S_v$.
\\
(f) If $v$ is an end point of $e$, then there is an entrance,
denoted $E^i_{v,e}$, 
of $\Lambda_e$ to $P^i_v$ with respect to $\varphi^i_v$. Put
 $Q_v^i=P_v^{i,\infty}$. 

\begin{quote}Then $(Q_v^i,E^i_{v,e})$ is a $(\Gamma_v,\Lambda_e)$-pair for
 $\varphi^i_v$ by Lemma \ref{l9}.
\end{quote}

\noindent
(g) If $v$ and $v'$ are two end points of $e$, then $(Q_v^i,E^i_{v,e})$ 
and $(Q_{v'}^i,E^i_{v',e})$ form a combinable pair. Denote the finite
 set $E^i_{e,*}=E^i_{v,e}\cap E^i_{v',e}$. 

\begin{quote}
The set $P^i_*=\bigcup_{v\in\VVV}P^i_{v,*}$ is called the {\em basic
configuration} (abbreviated BC) of $G=G(\TT)$ for $\varphi^i$. 
 A COP bijection 
$\xi:P^1_*\to P^2_*$ is called a {\em BC equivalence} if
 $\xi(P^1_{v,*})=P^2_{v,*}$ and $\xi\vert_{P^1_{v,*}}$ is a
 BP equivalence from $P^1_{v,*}$ to $P^2_{v,*}$ for each
$v\in\VVV$. 
\end{quote}

\noindent
(h) There is a BC equivalence $\xi:P_{*}^1\to P_{*}^2$.
\end{assumption}

For our purpose, the following example of BC is the most important.

\begin{example}\label{ex5}
Cosider a Fuchsian representation of the surface group $\Pi_5$ in
Figure 10.
Choose a lift $\widetilde T$ of the tree $T$ embedded in the surface
 to the universal covering space $\DD$. See Figure
 11. The lift of the curve $C_j$ to $\DD$ which intersects
$\widetilde T$ is denoted by $\widetilde C_j$.
The edge group $\Lambda_j$ is the stabilizer of $\widetilde C_j$.
As for the vertex group $\Gamma_i$, its convex core (of the limit set)
is contained in the
 region $\Sigma'_i$ depicted in Figure 11.

For a vertex of valency 1, the vertex group is generated by two
hyperbolic motions $a$ and $b$ such that $\tau([\widetilde a,\widetilde b])=1$. So it has 
a BP as in Figure 6. For a vertex of valency 3, generators $a,b$ of
 the vertex group satisfies $c(a,b)=1$, and it has a BP as in Figure 7.
The BP $P_1$ (resp.\ $P_3$) corresponding to the vertex group $\Gamma_1$
 (resp.\ $\Gamma_3$) is depicted in Figure 12.
The BC of $\Pi_5$ consists of 
50 points and satisfies all the conditions of Definition \ref{a3}. 
\end{example}

\begin{figure}\caption{}
\unitlength 0.1in
%
\end{figure}
\begin{figure}\caption{}
\unitlength 0.1in
%
%
\end{figure}

The following lemma is straightforward.

\begin{lemma}\label{last-lemma}
If $P_*^1$ is a BC for $\varphi^1\in\RR_G$, where
$G=G(\TT)$ is the fundamental group of a tree $\TT$, and if
$\psi^1$ is a $k$-fold lift of $\varphi^1$ for some $k\geq2$, then
$\pi_k^{-1}(P^1_*)$ is a BC for $\psi^1$. \qed
\end{lemma}

Before stating the main theorem of this section, we prepare a lemma.
By Theorem \ref{t1}, the BP equivalence $\xi\vert_{P^1_{v,*}}:P^1_{v,*}\to P_{v,*}^2$ extends
to a $(\varphi_v^1,\varphi^2_v)$-equivariant COP bijection 
$\hat\xi_v:Q_{v,*}^1\to Q_{v,*}^2$ for each vertex $v$. 
Notice that $Q_{v,*}^i=\varphi^i(\Gamma_v)P^i_{v,*}$.
\begin{lemma}\label{disjoint} 
There is a COP bijection $\hat\xi:\bigcup_vQ^1_{v,*}\to \bigcup_vQ^2_{v,*}$
such that $\hat\xi\vert_{Q^1_{v,*}}=\hat\xi_v$.
\end{lemma}

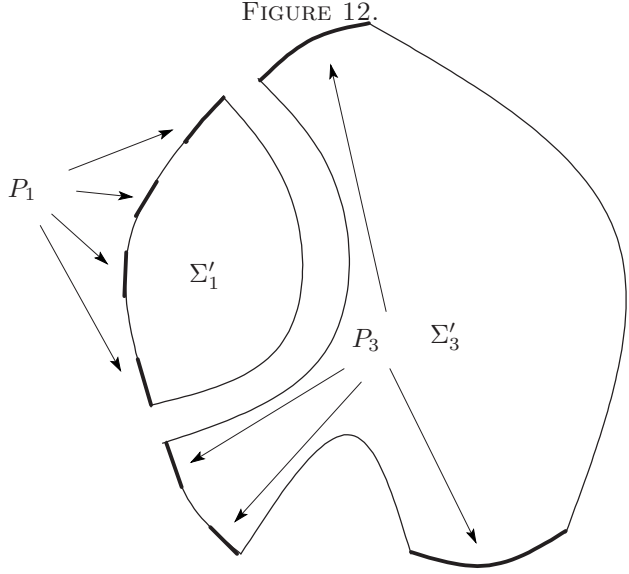
\bd  
If $v,v'$ are distinct vertices, then
 $P_v^i\cap P^i_{v'}\subset P_*^i$.
In fact, if $v,v'$ are adjacent, this follows from (g). If not,
$P^i_{v'}$ is contained in $\Int E^i_{v,e}$, where $e$ is the edge that starts at $v$ and tends toward the
direction of $v'$, which implies $P_v^i\cap P^i_{v'}=\emptyset$. Since
$Q_{v,*}^i\subset P^i_v$, the lemma follows from the fact that both
$\xi:P^1_*\to P^2_*$ and $\hat\xi_v:Q_{v,*}^1\to Q_{v,*}^2$ are COP bijections.
\qed

\begin{theorem}\label{t3}
The BC-equivalence $\xi:P^1_*\to P^2_*$ extends to a $(\varphi^1,\varphi^2)$-equivariant COP
 bijection $\hat\xi:\varphi^1(G)P^1_*\to\varphi^2(G)P^2_*$.
\end{theorem}

\bd 
The proof is by an induction on the number $n$ of vertices of $\TT$.
If $n=2$, this is just Theorem \ref{t2}.
Given $\TT$, delete a vertex $v$ of valency 1 and the edge $e$ which
starts at $v$.
Denote the resultant subtree by $\TT'$ and the other end point
of $e$ by $v'$. 
Then the group $G=G(\TT)$ can be written as an amalgamated product:
$$G=G(\TT')*_{\Lambda_e}\Gamma_v.$$ 
Let $$Q'^{i}=\varphi^i(G(\TT'))(\bigcup_{v\in \TT'}Q^i_v).$$
 Then $(Q'^{i},E^i_{v',e})$ is shown to be
a $(G(\TT'),\Lambda_e)$-pair  by virtue 
of Assumption \ref{a3} (c) and
successive use of Lemma \ref{next section}. 
Clearly the pair $(Q'^{i},E^i_{v',e})$ is combinable with the
$(\Gamma_v,\Lambda_e)$-pair $(Q^i_v,E^i_{v,e})$.
On the other hand, by the induction hypothesis,
$\xi$ has an $G(\TT')$-equivariant extension $\xi':Q'^1_*\to Q'^2_*$.
Moreover $\xi'$ and $\xi_v$ satisfy point (e) of Assumption \ref{a1}.
The proof is complete by Theorem \ref{t2}. \qed

\section{Robust basic configurations}

Again let $G=G(\TT)$ be the fundamental group of a tree $\TT$ of groups.
Assume that $\varphi^1\in\RR_G$ satisfies Assumption \ref{a3} for $\nu=1$, and let
$P^1_*$ be the associated BC.
Recall that for each vertex $v$ of $\TT$ and $l\geq2$,
$(P_v^1)^l$ is the BP for $\Gamma_v$ derived from the BP $P_v^1$.
(Definition \ref{d2}).
Denote $(P^1)^l_*=\bigcup_{v}(P_v^1)^l_*$.

For each point $x\in P^1_*$, the stabilizer $\Stab_{\varphi^1}(x)$ is
infinite cyclic
by Lemma \ref{l4}, Lemma \ref{l13} and a repeated use of Lemma \ref{l11}.
Denote by $x^+_l$ (resp.\ $x^-_l$) the point in $(P^1)^l_*$ right
(resp.\ left) adjacent to $x$.

\begin{definition}\label{robust}
The BC $P^1_*$ is called {\em robust} if for any point
$x\in P^1_*$ and any big
 $l$, one of the generators of $\varphi^1(\Stab_{\varphi^1}(x))$ maps
 the interval $[x^-_l,x^+_l]$ into a proper subinterval of it.
\end{definition}

\begin{lemma}
For a homomorphism $\varphi^1\in\RR_{\Pi_g}$ with $eu(\varphi^1)=2g-2$ ($g\geq2$),
the BC given by Examples \ref{ex4} and \ref{ex5} is robust.
\end{lemma}

\bd If we choose a Fuchsian representation as a model of $\varphi^1$, then
any point of the BC is a fixed point of a hyperbolic motion.
Any representation $\varphi^1$ with
$eu(\varphi^1)=2g-2$ is semiconjugate to the Fuchsian representation by
\cite{Matsumoto2}, showing the lemma.
\qed

\medskip
Finally we have the following theorem.

\begin{theorem}\label{t4} Assume that $\varphi^1$ admits a robust
BC $P^1_*$.
Then there is a neighbourhood $\mathcal U$ of $\varphi^1$ in $\RR_G$ such that if
$\varphi^2\in \mathcal U$, $\varphi^2$ admits a BC
 $P^2_*$
and a BC equivalence $\xi:P^1_*\to P^2_*$. 
\end{theorem}

\bd Choose $l$ large enough so that the condition of Definition
\ref{robust} is met by all the points $x$ in $P^1_*$ and that the intervals
$[x^-_l,x^+_l]$'s are disjoint. Let $g_x$ be the generator of
$\Stab_{\varphi^1}(x)$ such that 
$$\varphi^1(g_x)[x^-_l,x^+_l]\subset\Int[x^-_l,x^+_l].$$
Choose a neighbourhood $\mathcal U$ of
$\varphi^{1}$ so that for any $\varphi^2\in\mathcal U$ and $x\in
P^1_*$, 
we have
$$\varphi^2(g_x)[x^-_l,x^+_l]\subset\Int[x^-_l,x^+_l].$$
Let $\xi(x)$ be the leftmost point in
$\Fix(\varphi^2(g_x))\cap[x^-_l,x^+_l]$.
Then the set $$P_*^2=\{\xi(x)\mid x\in P_*^1\}$$ forms a BC for
$\varphi^2$, 
and the map $\xi$ is a BC equivalence. In fact, it is easy to see
that for any vertex $v$,
$$P^2_{v,*}=\{\xi(x)\mid x\in P^1_{v,*}\}$$ is a BP for $\varphi^2\vert_{\Gamma_v}$,
 because we have assumed that $P^1_{v,*}$ is a pure BP. \qed

\medskip
Joining this theorem with Lemma \ref{last-lemma}
and Theorem \ref{t3}, we get the following corollary, which conclude the
proof of Theorem \ref{locally stable}.

\begin{corollary}\label{c1}
If $\varphi^1\in\RR_G$ admits a robust BC, and
 $\psi^1\in\RR_G$ is a $k$-fold lift  of $\varphi^1$ ($k\geq 1$),
 then $\psi^1$ is locally stable.
\end{corollary}

\bd Let $P^1_*$ be a robust BC for $\varphi^1$. Then clearly
$\pi_k^{-1}(P_*^1)$ is a robust BC for $\psi^1$.
By Theorem \ref{t4}, there is a  BC $\widetilde P^2_*$ for any
 $\psi^2$ sufficiently near to $\psi^1$ 
and a BC equivalence $\xi:\pi_k^{-1}(P_*^1)\to\widetilde P^2_*$.
By Theorem \ref{t3}, the BC equivalence $\xi$ extends
to a $(\psi^1,\psi^2)$-equivariant 
COP map $\hat\xi:\psi^1(G)(\pi_k^{-1}(P_*^1))\to\psi^2(G)(\widetilde P^2_*)$.
The map $\hat\xi$ extends to a $(\psi^1,\psi^2)$-equivariant
semiconjugacy. 
\qed 

\section{Appendix A: The proof of Proposition \ref{first}}

We shall show that the semiconjugacy as defined in Definition
\ref{semiconjugacy} is an equivalence relation in
$\RR_G\setminus\RR_G^*$. All that need proof is the reflexiveness.
Let $\varphi^1,\varphi^2\in\RR_G\setminus \RR_G^*$. Assume there is a
degree one monotone map $h:S^1\to S^1$ such that 
\begin{equation}\label{h}
\varphi^2(g)\circ h=h\circ \varphi^1(g),\ \ \forall g\in G.
\end{equation}
Since $\varphi^i\in\RR_G\setminus \RR_G^*$,  $h$ is not
a constant map.
 Let $\widetilde h:\R\to\R$
be a lift of $h$ as in Definition \ref{dom}.
Notice that such a lift $\tilde h$ is unique up to the composition with
$T^n$, since the map
$h$ is nonconstant. (This is why we divide the definition of
semiconjugacy into two parts.)
Define $\widetilde h^\diamond:\R\to\R$ by
$$
\widetilde h^\diamond(y)=\inf\{x\in\R\mid \tilde h(x)=y\}.$$
Clearly $\widetilde h^\diamond$ commutes with $T$, and there is a degree one
monotone map $h^\diamond:S^1\to S^1$ such that 
$h^\diamond\circ
\pi=\pi\circ\widetilde h^\diamond$. The well-definedness of $h^\diamond$
is
garanteed by the uniqueness of $\widetilde h$.
Moreover if $h$, $h'$ and $h'\circ h$ are nonconstant monotone degree one maps, then
we have
$$
(h'\circ h)^\diamond=h^\diamond\circ(h')^\diamond.$$ 
Thus (\ref{h}) implies that
$$h^\diamond\circ\varphi^2(g^{-1})=\varphi^1(g^{-1})\circ h^\diamond,
$$
completing the proof.

\section{Appendix B: The  proof of Theorem \ref{Ghys}}
We assume that $\varphi\in\RR_G$ is type 1 and minimal, and will show that
$\varphi$ is proximal, the other implication being obvious. Call a closed interval $I\subset S^1$
$\varphi$-{\em contactible} if $\inf_{g\in G}\abs{\varphi(g)I}=0$.
First of all we have the following easy fact.
\begin{quote}
(1) For any $g\in G$ and any closed interval $I$, $I$ is $\varphi$-contractible
if and only if $\varphi(g)I$ is
$\varphi$-contractible. \qed
\end{quote}
Next let us show:
\begin{quote}
(2) There is $\delta>0$ such that if $\abs{I}<\delta$, then $I$ is
 $\varphi$-contractible.
\end{quote}

\bd Since $\varphi$ is not type 0, there is a nontrivial homeomorphism
$\varphi(g)$ which admits a fixed point. This shows that there is
a $\varphi$-contractible interval $J$.  
Since $\varphi$ is minimal,
the family 
$$
{\mathcal J}=\{\varphi(g)\Int J\mid g\in G\}$$
must cover $S^1$.
Now the Lebesgue number $\delta$ of the open covering $\mathcal J$
works.
\qed  

\medskip
Define a map $\widetilde U:\R\to\R$ by
$$\widetilde U(\widetilde x)=\sup\{\widetilde y\in(\widetilde x,\infty)\mid\pi([\widetilde
x,\widetilde y])\ \mbox{ is $\varphi$-contractible}\}.
$$
We have the following easy properties. 
\begin{quote}
(3) $\widetilde x+\delta\leq\widetilde U(\widetilde x)\leq\widetilde x+1$.
\qed
\end{quote}
\begin{quote} (4) $\widetilde U$ is monotone nondecreasing. \qed
\end{quote}
Also (1) implies the following.
\begin{quote}
(5)
For any $g\in G$ and a lift $\widetilde{\varphi(g)}$ of $\varphi(g)$ to
 $\R$, $$\widetilde{\varphi(g)}\circ\widetilde U=\widetilde U\circ\widetilde{\varphi(g)}.$$
Especially, $\widetilde U\circ T=T\circ\widetilde U$. \qed
\end{quote}
\begin{quote}
(6) The map $\widetilde U$ is injective.
\end{quote}

\bd Assume on the contrary that there is $\widetilde y\in\R$ such that
${\rm Cl}(\widetilde U^{-1}(\widetilde y))=[\widetilde x_0,\widetilde
x_1]$ is an interval. By the minimality of $\varphi$, there is a lift
$\widetilde{\varphi(g)}$ such that $\widetilde{\varphi(g)}(\widetilde
x_1)\in(\widetilde x_0,\widetilde x_1)$.
Then there is $\widetilde x_2\in(\widetilde x_0,\widetilde x_1)$ such
that $\widetilde{\varphi(g)}(\widetilde x_2)\in(\widetilde x_0,\widetilde
x_1)$
and $\widetilde{\varphi(g)}^{-1}(\widetilde x_2)\in(\widetilde x_1,\infty)$.
Now
$$\widetilde{\varphi(g)}(\widetilde
y)=\widetilde{\varphi(g)}\circ\widetilde U(\widetilde x_2)=\widetilde U
\circ\widetilde{\varphi(g)}(\widetilde x_2)=\widetilde y.$$
This shows
$$\widetilde U\circ\widetilde{\varphi(g)}^{-1}(\widetilde x_2)
=\widetilde{\varphi(g)}^{-1}\circ\widetilde U(\widetilde x_2)
=\widetilde{\varphi(g)}^{-1}(\widetilde y)=\widetilde y.
$$
This contradicts with the fact that
$\widetilde{\varphi(g)}^{-1}(\widetilde x_2)\not\in[\widetilde
x_0,\widetilde x_1]={\rm Cl}(\widetilde U^{-1}(\widetilde y))$. \qed

\begin{quote}
(7) $\widetilde U$ is bijective.
\end{quote}

\bd
Define $\widetilde V:\R\to\R$ by
$$\widetilde V(\widetilde x)=\inf\{\widetilde y\in(-\infty,\widetilde
x)\mid\pi([\widetilde y,
\widetilde x])\ \mbox{ is $\varphi$-contractible}\}.
$$
For any point $\widetilde x\in\R$,
and any point $\widetilde x_1$ in $(\widetilde x,\widetilde
U(\widetilde x))$, (6) implies that $\widetilde U(\widetilde
x)<\widetilde U(\widetilde x_1)$. This shows that the interval
$\pi([\widetilde x_1,\widetilde U(\widetilde x)])$ is $\varphi$-contractible.
Since $\widetilde x_1$ is an arbitrary point of $(\widetilde x,\widetilde
U(\widetilde x))$, this shows that $\widetilde V(\widetilde U(\widetilde
x))\leq \widetilde x$. Again by (6), we have in fact
$$\widetilde V(\widetilde U(\widetilde
x))= \widetilde x.
$$
The same argument shows that $\widetilde U\circ\widetilde V={\rm Id}$.
\qed

\medskip
By (4) and (7), $\widetilde U$ is a homeomorphism. By (5), there is
$U\in\HH$ such that $\pi\circ\widetilde U=U\circ\pi$. Also by (5),
$U$ commutes with any element of $\varphi(G)$.
Finally let us show:
\begin{quote}
(8) There is $k\in\N$ such that $U^k={\rm Id}$.
\end{quote}

\bd If $\Fix(U^k)$ is nonempty for some $k\in\N$, then $\Fix(U^k)$
must be invariant by $\varphi(G)$, since $U^k$ commutes with 
any element of $\varphi(G)$. That is, $\Fix(U^k)=S^1$, showing (8). 
If not, the rotation number of $U$ must be irrational, and there 
is a unique minimal set $X$ of $U$. Since $X$ is unique and since
$U$ commutes with any element of $\varphi(G)$, $X$ must be left
invariant
 by any element of $\varphi(G)$. Since $\varphi$ is
minimal, this implies $X=S^1$. That is, $U$ is topologically conjugate to
an irrational rotation. But then $\varphi(G)$ must be abelian, and
$\varphi$ must be of type 0. A contradiction. \qed

\medskip
To conclude, since $\varphi$ is assumed to be of type 1, we have $k=1$.
But by (3), this implies $\widetilde U=T$. That is,
$\varphi$ is proximal.

\end{document}